\documentstyle{amsppt}
\magnification=\magstep1
\pageheight{23truecm}
\NoBlackBoxes

\def\id{\operatorname{id}}
\def\cat{\operatorname{cat}}

\def\lk{\operatorname{lk}}
\def\invlim{\operatorname{invlim}}

\def\J{{\bold{J}}}

\def\Z{{\Bbb Z}}
\def\e{{ e}}

\def\R{{\Bbb R}}
\def\Sb{{\Bbb S}}
\def\P{{\Bbb P}}
\def\K{{\Bbb K}}
\def\T{{\Bbb T}}
\def\B{{\Bbb B}}
\def\ssm{{\smallsmile\,}}

\def\noi{\noindent}

\def\Ext{\operatorname{Ext}}

\def\Hom{\operatorname{Hom}}

\def\dim{\operatorname{dim}}
\def\im{\operatorname{im}}
\def\id{\operatorname{id}}
\def\diam{\operatorname{diam}}

\def\rank{\operatorname{rank}}

\def\ss{\smallskip}
\def\ms{\medskip}
\def\bs{\bigskip}
\def\kn{\mathaccent"7017 }
\def\kw{\mathbin{\raise 1pt \hbox {$\scriptstyle\rlap
{$\scriptstyle \sqcap $} \sqcup $}}}
\catcode`@=11
\def\dwiestrz#1#2{\,\vcenter{\m@th\ialign{##\crcr
      $\hfil\scriptstyle{\ #1\ }\hfil$\crcr
      \noalign{\kern2\p@\nointerlineskip}
      \rightarrowfill\crcr
      \noalign{\kern\p@\nointerlineskip}
      \leftarrowfill\crcr\noalign{\kern2\p@\nointerlineskip}
      $\hfil\scriptstyle{\ #2\ }\hfil$\crcr}}\,}
\catcode`@=13

\rightheadtext{Embedding compacta into products of curves}
\leftheadtext{A. Koyama, J. Krasinkiewicz and S. Spie{\D z}}

\topmatter

%\date
%\enddate

\title
Embedding compacta into products of curves \vskip-3mm
\endtitle

\author
Akira Koyama$^{(*)}$, J{\' o}zef Krasinkiewicz and Stanis{\l}aw
Spie{\D z} \vskip2mm
\endauthor

\address
Department of Mathematics, Faculty of Science, Suruga, Shizuoka
University, Shizuoka, 422-8529, Japan
\endaddress
\email sakoyam$\@$ipc.shizuoka.ac.jp
\endemail

\address
The Institute of Mathematics, Polish Academy of Sciences, ul.{\'
S}niadeckich 8, 00-950, Warsaw, Poland \endaddress
\address
Institute of Mathematics and Informatics, University of Opole, ul.
Oleska 48, 45-052 Opole, Poland
\endaddress
\email jokra$\@$impan.gov.pl
\endemail

\address
The Institute of Mathematics, Polish Academy of Sciences, ul.{\'
S}niadeckich 8, 00-950, Warsaw, Poland
\endaddress
\email
spiez$\@$impan.gov.pl
\endemail

\abstract We present some results on $n$-dimensional compacta
embeddable into $n$-dimensional Cartesian products of compacta. We
pay special attention to compacta embeddable into products of
1-dimensional compacta. Most of our basic results are proven under
the assumption that the compacta $X$ admit essential maps into the
$n$-sphere $\Sb^n$ (equivalently, the {\v Cech} cohomology
$H^n(X)\neq 0$). Our investigations have been inspired by some
results in this direction established by Borsuk, Cauty, Dydak,
Koyama and Kuperberg. The results of the present paper may be
viewed as an extension of the theory developed so far by these
authors.

First, we prove that if $X$ is an $n$-dimensional compactum with
$H^n(X) \neq 0$ that embeds in a product of $n$ curves then there
exists an algebraically essential map $X\to\T^n$ into the
$n$-torus. Then we show that the same is true if $X$ embeds in the
$n$th symmetric product of a curve. The existence of such a
mapping implies the existence of elements $a_1,\cdots,a_n\in
H^1(X)$ whose cup product $a_1\ssm\cdots\ssm a_n$ is non-zero.
Consequently, $\rank H^1(X)\geq n$ and $\cat X >n$. In particular,
$\Sb^n$, $n\geq 2$, is not embeddable in the $n$th symmetric
product of any curve. Next, we introduce some new classes of
$n$-dimensional continua and show that embeddability of locally
connected {\it quasi $n$-manifolds} into products of $n$ curves
also implies $\rank H^1(X)\geq n$. It follows that some
2-dimensional contractible polyhedra are not embeddable in
products of two curves. On the other hand, we show that any
collapsible 2-dimensional polyhedron can be embedded in a product
of two trees. We answer a question posed by Cauty proving that
closed surfaces embeddable in products of two curves can be also
embedded in products of two graphs. We prove that no closed
surface $\neq\T^2$ lying in a product of two curves is a retract
of that product.
\endabstract

\keywords
Embeddings, surfaces, cohomology groups, symmetric products
\endkeywords

\subjclass 54E45, 55M10, 55U25, 57N05, 57N35
\endsubjclass

\thanks
$^{(*)}$ This paper has been originated during a visit of the
first author in the Institute of Mathematics, Polish Academy of
Sciences, under the exchange program between Japan Society for the
Promotion of Sciences and the Polish Academy of Sciences.
\endthanks

\endtopmatter

\vfill \break

\centerline{\eightbf CONTENTS}

\ms\centerline{\eightrm CHAPTER 1}

\centerline{\eightbf Introduction}

\item{}{\sevenit Problems to Chapter 1 \hfill{9}}

\ms\centerline{\eightrm CHAPTER 2}

\centerline{\eightbf Embeddability into products of curves and the
first cohomology}

\item{\sevenrm 2A.} {\sevenit Mapping products of compacta into
corresponding products of spheres \hfill{11}}

\vskip-1mm\item{\sevenrm 2B.} {\sevenit Mapping into products
involving 1-dimensional factors \hfill{14}}

\ms\centerline{\eightrm CHAPTER 3}

\centerline{\eightbf Embeddability into products and the
cohomology ring}

{\sevenit
\item{\sevenrm 3A.} {\sevenit Algebraically non-trivial mappings into
products of compacta and non-zero cup products of \hfill
 \vskip-1mm
spherical elements \hfill{17}}

\vskip-1mm\item{\sevenrm 3B.} {\sevenit Categories of spaces
\hfill{21}}

\ms\centerline{\eightrm CHAPTER 4}

\centerline{\eightbf Embedding compacta into symmetric product of
curves}

\item{\sevenrm 4A.} {\sevenit Symmetric products as functors
\hfill{23}}

\vskip-1mm \item{\sevenrm 4B.} {\sevenit Symmetric products of a
pointed space \hfill{24}}

\vskip-1mm\item{\sevenrm 4C.} {\sevenit Symmetric products of a
bouquet of k circles and skeleta of torus $\scriptstyle\T^k$
\hfill{25}}

\vskip-1mm\item{\sevenrm 4D.} {\sevenit Algebraic properties of
skeleta of a torus \hfill{26}}

\vskip-1mm\item{\sevenrm 4E.} {\sevenit The Main Theorem
\hfill{27}}

\ms\centerline{\eightrm CHAPTER 5}

\centerline{\eightbf Locally connected generalized manifolds in
product of curves}

\item{\sevenrm 5A.} {\sevenit Definitions and basic properties of certain
generalized manifolds \hfill{29}}

\vskip-1mm\item{\sevenrm 5B.} On locally connected quasi
manifolds. From embeddings into products of curves to embed-\hfill
\vskip-1mm dings into products of graphs\hfill{33}}

\vskip-1mm\item{\sevenrm 5C.} {\sevenit  Ramified manifolds in
products of graphs\hfill{37}}

\vskip-1mm\item{\sevenrm 5D.} {\sevenit Product structure of
generalized manifolds lying in products of graphs\hfill{40}}

\vskip-1mm\item{\sevenrm 5E.} {\sevenit Contractible 2-dimensional
polyhedra in products of two graphs\hfill{43}}

\vskip-1mm\item{}{\sevenit Problems to Chapter 5\hfill{46}}

\ms\centerline{\eightrm CHAPTER 6}

\centerline{\eightbf Embedding surfaces into product of two
curves}

\item{\sevenrm 6A.} {\sevenit Another proof of the Borsuk theorem
\hfill{47}}

\vskip-1mm\item{\sevenrm 6B.} {\sevenit On surfaces lying in
products of two graphs\hfill{47}}

\vskip-1mm\item{\sevenrm 6C.} {\sevenit Surfaces in products of
$\scriptstyle\theta_n$-curves\hfill{49}}

\vskip-1mm\item{\sevenrm 6D.} {\sevenit On Cauty's results about
embeddability of non-orientable surfaces into product of
graphs\hfill{53}}

\vskip-1mm\item{\sevenrm 6E.} {\sevenit Retracting products onto
surfaces\hfill{60}}

\vskip-1mm\item{\sevenrm 6F.} {\sevenit Embedding bordered
surfaces in the "three-page book"\hfill{62}}

\vskip-1mm\item{\sevenrm 6G.} {\sevenit Embedding surfaces in the
symmetric products of curves\hfill{63}}

\vskip-1mm\item{\sevenrm 6H.} {\sevenit Embedding surface-like
continua in products of two curves\hfill{65}}

\vskip-1mm\item{} {\sevenit Problems to Chapter 6\hfill{66}}

\ms\item{}{\eightbf Appendix.} {\sevenit Some results on the
tensor product of groups\hfill{66}}

\ms\item{}{\eightbf References} {\sevenit\hfill{70}}

\vfill \break

\document
%%%%%%%%%% Introduction %%%%%%%%%%%%
\head 1. Introduction
\endhead

Throughout this paper we use the following standard terminology
and notation. All {\it spaces} discussed in this paper are
metrizable and all {\it mappings} (also called {\it maps}) are
continuous. By a {\it compactum} we mean a compact (metric) space,
by a {\it continuum} we mean a connected compactum, and by a {\it
curve} we mean a 1-dimensional continuum. Sometimes we write
$X\approx Y$ to indicate that $X$ is homeomorphic to $Y$.

By $\B^n$, $n\geq1$, we denote the closed unit $n$-ball in the
Euclidean $n$-space $\R^n$; and $\B^0=\R^0$ stands for the
one-point set $\{0\}$. A space homeomorphic to $\B^n$ is called a
(closed) $n$-{\it disc}. By $\Sb^{n-1}$ we denote the unit
$n$-sphere in $\R^n$ - the boundary of $\B^n$; and $\Sb^{-1}$
stands for the empty set. A space homeomorphic to $\Sb^n$ is
called a {\it topological $n$-sphere}; a space homeomorphic to
$\Sb^1$ is called a {\it $($topological$)$ circle} (or a {\it
simple closed curve}). A space homeomorphic to the open unit
$n$-ball $\kn\B^n=\B^n \backslash \Sb^{n-1}$ is called an {\it
open $n$-disc}. As usual, by the $n$-torus $\T^n$, $n\geq1$, we
mean the $n$-fold product $\Sb^1\times\cdots\times \Sb^1$. In
particular, $\T^1=\Sb^1$. By $\T^0$ we mean a one-point space. A
space homeomorphic to $\T^n$ is called a {\it topological
n-torus}. By a {\it dendrite} we mean a non-degenerate locally
connected continuum containing no simple closed curve. A
non-degenerate continuum is said to be a {\it local dendrite} if
every point has a neighborhood which is a dendrite. It is known
that dendrites coincide with 1-dimensional compact absolute
retracts, and local dendrites -- with 1-dimensional compact
absolute neighborhood retracts (cf. \cite{Kur}).

By a {\it complex} we mean a standard $CW$ complex. Terminology
concerning some special complexes introduced in this work will be
presented in appropriate places. Throughout the paper by a {\it
polyhedron} we mean the underlying space $|K|$ of a finite regular
$CW$ complex $K$. Recall that a $CW$ complex $K$ is said to be
{\it regular} if each cell $\sigma\in K$ admits a characteristic
map $\varphi_{\sigma}: {\B}^n \to \sigma$, where $n=\dim \sigma$,
which is a homeomorphism. Any {\it simplicial complex} $K$ is
regular, such $K$ is called a {\it triangulation} of $|K|$.

In 1958 J. Nagata \cite{Na1} discovered the following remarkable
theorem.

\proclaim{Theorem 1.1 (Nagata)} Every $n$-dimensional space,
$n\ge2$, can be embedded in the topological product
$X_1\times\cdots \times X_{n+1}$ of $1$-dimensional spaces.\qed
\endproclaim

\noi (This result has been subsequently discussed and refined by
many authors, see e.g. \cite{Bw}, \cite{I-M}, \cite{L}, \cite{Mi},
\cite{Ol}, \cite{St}, \cite{T}.)

A few years later, in his book devoted to dimension theory, he
also asked the following question: "It is an open problem whether
every $n$-dimensional metric space can be topologically imbedded
in the topological product of $n$ $1$-dimensional metric
spaces?"(\cite{Na2},p.163). This question was answered in the
negative by K.~Borsuk \cite{Bo3} in $1975$. Actually, Borsuk
proved the following interesting result.

\proclaim{Theorem 1.2 (Borsuk)} The $2$-sphere $\Sb^2$ is not
embeddable in any product of two curves. $($Analogous result holds
for all spheres $\Sb^n$, $n\geq 3$.$)$\qed
\endproclaim

In the present paper we shall show that some 2-dimensional
contractible (so, acyclic) polyhedra have this property as well.

The above two results justify the following notion. An
$n$-dimensional space, $n\geq2$, is said to be {\it ordinary} if it
can be embedded in a product of $n$ 1-dimensional spaces; otherwise
we call it {\it exceptional}. Hence $\Sb^2$ and some other spaces
are exceptional.

Thus every class of $n$-dimensional, $n\ge2$, metric spaces splits
in a natural way into two complementary subclasses - the ordinary
and exceptional spaces. Generally speaking, in our paper we make
an attempt at better understanding this splitting. We focus our
attention on some special ordinary compacta. The fundamental
problem of characterization of ordinary $n$-dimensional compacta
(even for $n=2$) seems unattainable at the moment. Using this new
terminology one could have expressed our results in a more concise
form.

Any 1-dimensional compactum can be embedded in the Menger curve
$\mu$. It follows that an $n$-dimensional compactum is ordinary if
and only if it can be embedded in the $n$-fold product $\mu^n$.

Let us recall that a mapping $f:X\to Y$, where $X$ is metric, is
said to be an {\it $\varepsilon$-mapping} if
$\diam{}f^{-1}(y)<\varepsilon$ for each $y\in f(X)$. It is well
known that any $n$-dimensional compactum $X$ admits surjective
$\varepsilon$-mappings onto $n$-dimensional (compact) polyhedra
for each $\varepsilon > 0$. In particular, each curve admits
$\varepsilon$-mappings onto graphs for each $\varepsilon > 0$. (By
a {\it graph} we mean a $1$-dimensional connected polyhedron.) It
follows that if $X$ is $n$-dimensional and ordinary then it admits
$\varepsilon$-mappings into products of $n$ graphs for each
$\varepsilon > 0$ (equivalently, $X$ admits $\varepsilon$-mappings
into $\mu^n$). In Problem 1.4 located at the end of this chapter
we inquire about the reverse implication. If $X$ admits
$\varepsilon$-mappings into $Y$ for each $\varepsilon >0$ then $X$
is said to be {\it quasi-embeddable} into $Y$.

In this context we mention an important and quite large class of
extremely complicated ordinary continua. Namely, R.~Pol \cite{P}
proved a surprising result that every $n$-dimensional, $n\ge2$,
hereditarily indecomposable continuum is ordinary. Moreover, he
showed that such a continuum can be embedded in a product of $n$
hereditarily indecomposable curves.

All {\it manifolds} discussed in this paper are assumed to be
compact and connected (possibly with non-empty boundary), unless
opposite is explicitly stated. The interior of $M$ will be often
denoted by $\kn M$. A manifold $M$ is {\it closed} if its boundary
is empty, $\partial{}M=\emptyset$. A manifold with non-empty
boundary will be called {\it bordered manifold}. A mapping
$f:X\to{M}$, where $M$ is a closed manifold, is said to be {\it
essential} if every mapping $g:X\to{M}$ homotopic to $f$ is
surjective.

The symbol $H^{\ast}(\cdot\,;\,G)$ is used to denote the {\v Cech}
cohomology functor with coefficients in an Abelian group $G$. In
some cases where no confusion is likely to occur we shall write
$f^{\ast}$ instead of $H^{\ast}(f\,;\,G)$, where $f$ is a mapping.
If $H^n(f\,;\,G)\neq 0$ then $f$ is said to be {\it non-trivial
with respect to} $H^n(\cdot\,;\,G)$. A mapping $f:X\to M$, where
$M$ is a closed $n$-manifold, is said to be {\it algebraically
essential} if $H^n(f;G)\neq 0$ for some group $G$ (that is, $f$ is
non-trivial with respect to $H^n(\cdot\,;\,G))$. Any algebraically
essential mapping is essential. (In fact, this readily follows
from the following result of M. Brown \cite{Br, p. 94}:\ss

{\it If $M$ is a bordered manifold then there is a surjective map
$H:(\partial M) \times I\to M$ such that $H(x,0)=x$ for each
$x\in\partial M$, $H|(\partial M)\times [0,1):(\partial M)\times
[0,1)\to M$ is an embedding, $H^{-1}(H((\partial M)\times
\{1\}))=(\partial M)\times \{1\}$, and $\dim H((\partial M)\times
\{1\})\le \dim M-1$}.\ss

It follows that the set $Y=H((\partial M)\times \{1\}))$ has
dimension $< \dim M$ and it is a strong deformation retract of
$M$, because $M$ is homeomorphic to the mapping cylinder of the
mapping $\varphi:\partial M\to Y$ given by
$\varphi(x)=H(x,1)$.)\ss

The cohomology functor $H^{\ast}(\cdot\,,\,\Z)$ with integer
coefficients $\Z$ will be abbreviated to $H^{\ast}(\cdot)$. Thus
the groups $H^{\ast}(X\,,\,\Z)$ and the homomorphisms
$H^{\ast}(f\,,\,\Z)$ will be written briefly $H^{\ast}(X)$ and
$H^{\ast}(f)$, respectively. By the Hopf Classification Theorem,
cf. \cite{Sp, p. 431}, for any $n$-dimensional space $X$, the
group $H^n (X)$ is in one-to-one correspondence with the set of
homotopy classes of maps $X\to \Sb^n$. Non-zero elements
correspond to homotopy classes of essential maps. As usual,
$H_*(X)$ denotes the singular homology functor with integer
coefficients.

Let $g_1,\cdots,g_k$ be elements of an Abelian group $G$. They are
said to be {\it linearly independent} (over $\Z$) if the equality
$n_1g_1+\cdots+n_kg_k=0$, $n_i\in\Z$, implies $n_1=\cdots=n_k=0$.
By the {\it rank} of $G$, denoted $\rank G$, we mean the maximal
number of linearly independent elements in $G$ (over $\Z$). We
write $G\cong H$ to denote that a group G is isomorphic to a group
H. \ss

The original proof of the Borsuk theorem was not elementary. Two
simpler proofs have been subsequently supplied by J.~van Mill and
R.~Pol \cite{M-P}. Here we supply several other proofs (and one,
perhaps the simplest possible, in article 6A). Actually, we shall
see that Borsuk's theorem readily follows from each of several
results proved in our paper.

Generalizing the Borsuk theorem W.~Kuperberg \cite{Ku} showed in
$1978$ that no closed $n$-manifold, $n \geq 2$, with finite
fundamental group, can be embedded in any product of an
$(n-1)$-dimensional compactum and a curve. (Here we obtain this
result as a corollary to Theorem 2B.1, see Corollary 2B.5.) In
particular, the projective plane is not embeddable in any product
of two curves. He also noted that

\proclaim{Theorem 1.3 (Kuperberg)} Any orientable closed surface
different from the $2$-sphere can be embedded in a product of two
graphs. \qed\endproclaim

\noi Using this result one can show that certain class of pretty
complicated 2-dimensional continua consists of ordinary continua.
Such continua have been constructed by Dranishnikov \cite{Dr$_2$}
(called by him {\it fractal Riemann surfaces}). (Any such a
continuum contains no open non-void subset embeddable in the plane
and admits mappings as small as we please onto orientable surfaces
of arbitrary high genus). A common feature of the continua studied
by Dranishnikov and Pol is their infinite rank of the first {\v
C}ech cohomology.

Kuperberg also asked if there is a non-orientable closed surface
embeddable in a product of two curves. The question was answered
in $1984$ by R.~Cauty \cite{C1}, who proved the following
remarkable result.

\proclaim{Theorem 1.4 (Cauty)} Any non-orientable closed connected
surface $M$ can be embedded in a product of two graphs if and only
if genus of $M$ is $\geq 6$. \qed \endproclaim

\noi (The latter condition is equivalent to either
$\rank{}H^1(M)\ge5$, or Euler characteristic of $M$ is $\le -4$.)
We take the opportunity to recall some relations between three
classic numerical invariants - $g(M)$, $\chi(M)$, $r(M)$ -
associated with any closed surface $M$. By $g(M)$ - the {\it
genus} of $M$ - we mean maximal number of mutually disjoint simple
closed curves on $M$ whose union does not separate $M$. As usual,
$\chi(M)$ stands for the Euler characteristic of $M$. And we put
$r(M)=\rank H_1(M)$. From the universal coefficient theorem for
cohomology we have the following short exact sequence:
$$
0 \to {\Ext}(H_0(M),\Z)\to H^1(M)\to {\Hom} (H_1(M),\Z)\to 0.
$$
It follows that $r(M)=\rank H^1(M)$ as well (see \cite{Sp, p. 244,
5.5.4}). If $M$ is orientable then $g(M)$ is equal to the so
called {\it handle number} of $M$ (that is, the number of handles
attached to $\Sb^2$ in a canonical presentation of $M$); for $M$
non-orientable it is the so called {\it crosscup number} of $M$
(that is, the number of the M\"{o}bius strips attached to $\Sb^2$
in a canonical presentation) of $M$. Moreover, we have the
following equalities:

\ms ($\bullet$) $g(M)=1-\frac{1}{2}\chi(M)= \frac{1}{2}r(M)$ {\it
for $M$ orientable}, and \ss

($\sim\bullet$) $g(M)=2-\chi(M)= 1+r(M)$ {\it for $M$
non-orientable}.\bs

The present paper has been inspired by the above results and by
the following theorem (which also implies the Borsuk theorem) due
to Dydak and Koyama \cite{D-K} and proved in 2000. This theorem
readily follows from our Corollary 3A.4.

\proclaim{Theorem 1.5 (Dydak, Koyama)} Let $X$ be a compact subset
of the product of $n$ curves, $n>1$, and let $G$ be an Abelian
group such $H^n(X;G) \neq 0$. Then $H^1(X;G)\neq 0$. \qed
\endproclaim

\centerline {*\hskip5mm  *\hskip5mm  *}

\ms The remaining part of this introduction is basically devoted
to a brief summary of the main results of our paper. Some extra
remarks are also included.

In Chapters 2 and 3 we first study properties of compacta
admitting algebraically non-trivial maps into products of
compacta. Then we investigate the cohomology ring of certain
compacta embeddable into products of curves. Some consequences
concerning the {\it category} of the compacta are derived. The
following theorem is a simplified version of some basic results.
Actually, its conclusion also holds under weaker stipulation. It
is enough to require that there is a mapping $f:X\to
Y_1\times\cdots\times Y_n$ which is non-trivial with respect to
$H^n(\cdot)$.

\ms \noi {\bf Theorem 1.6} (cf. Corollary 2A.5, Corollary 2B.3,
and Corollary 3A.1){\bf .} {\it Let $X$ be a compact subset of the
product $Y_1\times\cdots\times Y_n$ of $n$ curves with $H^n(X)
\neq 0$. Then there is an algebraically essential map $X\to \T^n$.
Consequently, there exist elements $a_1,\cdots a_n \in H^1(X)$
whose cup product $a_1 \ssm \cdots \ssm a_n \in H^n(X)$ is
non-zero. Such elements are linearly independent, hence $\rank
H^1(X)\ge n$. Moreover, $\cat X>n$.}\qed\ms

As an application we see that {\it the Klein bottle $\K$ can not
be embedded in any product of two curves} (because $H^1(\K)\cong
\Z)$. This conclusion can be also derived from the above result of
Cauty. But it does not follow from the Dydak-Koyama theorem
(because $H^1(\K;G)=0$ if and only if $G=0$, as $G$ is a direct
summand of $H^1(\K;G)$).

An analogous theorem holds for symmetric products, see Theorem
4G.1. It follows that no $\Sb^n$, $n\geq 2$, can be embedded in
the $n$th symmetric product of a curve (see Corollary 4G.3). This
is an analogue of the Borsuk Theorem 1.2. And it answers in the
negative the following question posed by Illanes and Nadler
\cite{I-N, Question 83-14}: {\it Is the $2$-sphere embeddable in
the second symmetric product of a curve?} The theorem also implies
that neither the projective plane $\P^2$ nor the Klein bottle $\K$
can be embedded in the symmetric product of a curve (see Corollary
6G.2).

An $n$-dimensional continuum $X$, $n\geq 1$, is said to be a {\it
quasi n-manifold} if for every point $x\in X$ there is an open
neighborhood $V$ of $x$ such that every closed subset of $X$ which
separates $X$ between $x$ and $X\setminus V$ and has dimension
$\leq n-1$ admits an essential map into $\Sb^{n-1}$ (see Section
5A for details). This and some other classes of $n$-dimensional
continua have been defined in Chapter 5. Each class comprises all
$n$-manifolds, and we have the following results.

\ms\noi{\bf Theorem 1.7} (cf. Theorem 5B.1){\bf.} {\it Let $X$ be
a locally connected quasi $n$-manifold, $n\ge2$, with $H^1(X)$ of
finite rank. If $X$ embeds in a product of $n$ curves then there
exists an embedding $g =(g_1, \cdots ,g_n):X \to P_1 \times \cdots
\times P_n$ such that} \roster
\item "(1)" {\it each $g_i$ is a monotone surjection},
\item "(2)" {\it each $P_i$ is a graph with no endpoint}.\qed
\endroster

\ms\noi{\bf Corollary 1.8} (see Corollary 5B.3){\bf.} {\it If a
closed $n$-manifold is embeddable in a product of $n$ curves then
it is also embeddable in a product of $n$ graphs.} \qed

\ms\noindent It follows that {\it if a closed surface can be
embedded in a product of two curves then it can be also embedded
in a product of two graphs}. This answers a question posed by
Cauty. A harder variant of this question where the word "surface"
is replaced by "$2$-dimensional polyhedron" is still open
\cite{C1}.

We also prove a strong version of the following

\ms\noi{\bf Theorem 1.9} (cf. Theorem 5D.5){\bf.} {\it Let $X$ be
a locally connected quasi $n$-manifold lying in a product of $n$
curves. Then $\rank H^1(X)\ge n$.} \qed\ms

An edge of a $2$-dimensional regular $CW$ complex $K$ is said to
be {\it free} if it is incident with exactly one $2$-cell of $K$.
Such an edge is also said to be free in $|K|$.

\proclaim{Corollary 1.10} No contractible $2$-dimensional
polyhedron $|K|$ with no free edge can be embedded in a product of
two curves. \qed
\endproclaim

\noi There are two well known polyhedra satisfying the hypotheses
of this corollary: the Borsuk example \cite{Bo1} (which occurred
in \cite{Z} under the name {\it "dunce hat"}), and the {\it "Bing
house"}, cf. \cite{R-S}. Hence neither can be embedded in a
product of two curves.

We also prove some results on $2$-dimensional polyhedra. Here we
quote three of them.

\ms\noi{\bf Theorem 1.11} (cf. Theorem 5E.1){\bf .} {\it Let $X$
be a $2$-dimensional connected polyhedron. If $X$ can be embedded
in a product of two curves and $\rank{}H_1(X)\le2$, then $X$
collapses to either a point, or a graph, or a torus. In
particular, $X$ is collapsible if $\rank{}H_1(X)=0$.} \qed

\ms\noi{\bf Theorem 1.12} (cf. Theorem 5E.3) {\bf.} {\it Let $X$
be a $2$-dimensional polyhedron. If $X$ is collapsible then $X$
can be embedded in a product of two trees.} \qed

\ms \noi (By a {\it tree} we mean a connected graph containing no
circle.)

\ms\noi{\bf Theorem 1.13} (see Corollary 5E.7){\bf.} {\it The cone
over an $n$-dimensional polyhedron can be embedded in a product of
$n+1$ copies of an $m$-od.} \qed

\ms \noi (By an $m$-{\it od} we mean the cone over an $m$-element
set.)

In Chapter 6 we present a detailed discussion of the problem of
embeddability of surfaces into products of two curves. The results
of Kuperberg (Theorem 1.3) and Cauty (Theorem 1.4) will be given
new proofs. In section 6C the Kuperberg theorem is improved by
showing that each orientable surface different from $\Sb^2$ can be
embedded in $\Theta_n\times \Theta_n$, where $n=g(M)+1$ and
$\Theta_n$ denotes the canonical $\theta_n$-curve (see 6B for
definition). In section 6D we prove the following refinement of
the Cauty result.

\ms \noi{\bf Theorem 1.14} (cf. Theorem~6D.1) {\bf.} {\it Let $M$
be a closed $2$-manifold in the product $Y_1\times Y_2$ of two
curves. Then}

\item\item{\rm(i)} {\it $\rank H_1(M)\le3$ implies $M=P_1\times P_2$,
where each $P_i$ is a circle in $Y_i$};

{\rm(ii)} {\it $\rank H_1(M)=4$ implies $M\subset P_1\times P_2$,
where each $P_i$ is a $\theta$-curve in $Y_i$}. \qed\ms

In section 6E we show that the torus is the only surface which,
when embedded in a product of two curves, is a retract of that
product (see Theorem 6E.1). In section 6F we show that any
bordered surface can be embedded in "the three-page book" (Theorem
6F.1). Section 6G contains some results on non-embeddability in
the second symmetric product of a curve, which have been mentioned
above.

One easily sees that if a space embeds in the product of two
curves $X$ and $Y$ then it also embeds in the second symmetric
product of a curve $Z$ (because $X\times Y$ embeds into the second
symmetric product of $Z=X \vee Y$ - the {\it one-point union} (or,
the {\it bouquet}) of $X$ and $Y$, here $X,Y$ are pointed spaces).
We shall show that the reverse implications is not true. The former
statement implies that any orientable surface different from
$\Sb^2$ and any non-orientable surface of genus $\ge6$ can be
embedded in the second symmetric product of a graph.

We complete this chapter with some remarks on
$\varepsilon$-mappings. Let us begin with a general result.\ms

\proclaim {Theorem 1.15} Let X and Y be n-dimensional compacta,
$n\geq 1$. If $H^n(X)\neq 0$ then there is an $\varepsilon >0$
such that every $\varepsilon$-mapping $f:X\to Y$ is non-trivial
with respect to $H^n(\cdot)$.\endproclaim

\demo{Proof} By the hypothesis $H^n(X)\neq 0$ there is an
essential map $u:X\to \Sb^n$. By a result of Eilenberg \cite{E}
there is an $\varepsilon
>0$ such that for every surjective $\varepsilon$-mapping $v:X\to
Z$ there is a mapping $w:Z \to \Sb^n$ such that dist $(u,w \circ
v)<2$. Consequently, $u\simeq w\circ v$. Now, consider any
$\varepsilon$-mapping $f:X\to Y$, and let $v_0:X\to f(X)$ be
determined by $f$. Then by the Eilenberg theorem there is a
mapping $w_0:f(X)\to \Sb^n$ such that $u\simeq w_0\circ v_0$.
Since $\dim Y\leq n$ and $f(X)$ is a closed subset of $Y$, there
is an extension $g:Y\to \Sb^n$ of $w_0$. It follows that $u\simeq
g\circ f$. Hence $H^n(f)\neq 0$ because $H^n(u)\neq 0$. This
completes the proof.\qed
\enddemo

\proclaim{Corollary 1.16} {\rm(a)} Let X be a compactum
quasi-embeddable in the product $Y_1\times\cdots\times Y_k$ of
finite dimensional compacta. If $H^n(X)\neq 0$, where
$n=n_1+\cdots +n_k$ and $n_i=\dim Y_i \geq 1$, then X admits a
mapping into $\Sb^{n_1}\times\cdots\times\Sb^{n_k}$ which is
non-trivial with respect to $H^n(\cdot)$.

{\rm(b)} Let X be a compactum quasi-embeddable in the symmetric
product $SP^n(Y)$, where Y is a curve. If $H^n(X)\neq 0$ then X
admits a mapping into $\T^n$ which is non-trivial with respect to
$H^n(\cdot)$.
\endproclaim

\demo{Proof} (a) By Theorem 1.15 there is a mapping $f:X\to
Y_1\times\cdots\times Y_k$ which is non-trivial with respect to
$H^n(\cdot)$. Hence the conclusion follows from Corollary 2A.3
(which is independent of Theorem 1.15).

(b) As in case (a), there is a mapping $f:X\to SP^n(Y)$ which is
non-trivial with respect to $H^n(\cdot)$. The conclusion follows
from Theorem 4G.1. \qed
\enddemo

\proclaim {Corollary 1.17} No continuum $X$ with $H^1(X)=0$ and
$H^n(X)\neq 0$, where $n\geq 2$, is quasi-embeddable in either a
product of n curves or the nth symmetric product of a
curve.\qed\endproclaim

\ms\noi {\bf Note.} This corollary improves the Borsuk Theorem
1.2. Its particular case, for $X=\Sb^n$ and the Cartesian product
of $n$ curves, was observed in \cite{M-P}.

\bs\centerline {\bf Problems to Chapter 1}

\bs The following problem is of fundamental importance for the
theory which has been developed in this paper and is certainly
very hard.

\proclaim{Problem 1.1} Characterize n-dimensional ordinary
compacta for each $n\ge 2$.
\endproclaim

Any $n$-dimensional compactum can be regarded as a subset of the
cube $I^{2n+1}$. Hence in an attempt of solving the above problem
it may be helpful to start with studying ordinary $n$-dimensional
compacta lying in lower dimensional cubes. Clearly, all
$n$-dimensional compacta lying in $I^n$ are ordinary by
definition. The following is the first, very special and already
highly nontrivial case of the problem.

\proclaim{Problem 1.2} Characterize $2$-dimensional ordinary
compacta lying in $I^3$. {\rm(In general:} Characterize
n-dimensional ordinary compacta lying in $I^{n+1}$.$)$
\endproclaim

\proclaim{Problem 1.3} Can $\mu^2$ be embedded in $I^4$? {\rm(In
general:} Can $\mu^n$, $n\ge 2$, be embedded in $I^{2n}$?$)$
\endproclaim

Our next problem is of great interest, we believe it has
affirmative solution.\ms

\proclaim{Problem 1.4} Suppose $X$ is a compactum which admits
$\varepsilon$-mappings into $\mu^n$ for each $\varepsilon >0$
$($i.e. X is quasi-embeddable in $\mu^n$$)$. Can X be embedded in
$\mu^n$?
\endproclaim

In connection with Theorems 1.6 and 1.9 our next question is also
of great interest.

\proclaim{Problem 1.5} Let $X$ be a locally connected quasi
$n$-manifold lying in a product of $n$ curves. Does X admit an
essential map $X\to \Sb^n$?
\endproclaim

For a compactum $X$ and a natural $n\ge2$ by $F_n(X)$ we denote
the hyperspace of all non-void subset of $X$ composed of at most
$n$ points, with the Hausdorff metric (or the Vietoris topology).
We shall call it {\it the hyperspace of $X$ of at most $n$
points}. (In \cite{Bo2} it is called {\it the $n$}th {\it potency
of $X$}.) One easily sees that $F_2(\Sb^1)$ is homeomorphic to the
M\"{o}bius strip. An interesting result of Bott (see \cite{Bo},
cf. \cite{Bo2}) asserts that $F_3(\Sb^1)$ is homeomorphic to
$\Sb^3$. In 1947 Wu gave complete description of $H^*(F_n(\Sb^1))$
for each $n\geq 1$ (see \cite{Wu}). In particular, he showed that
$H^{2n}(F_{2n}(\Sb^1))=0$. This implies that $\Sb^{2n}$ is not
embeddable in $F_{2n}(\Sb^1)$ because dim $F_{2n}(\Sb^1)= 2n$. In
a recent paper by Chinen and Koyama \cite{Ch-K} the latter result
has been extended to all odd $n>3$. Thus, no $\Sb^n$, $n\geq 4$,
can be embedded in $F_n(\Sb^1)$. But the following harder problem
is still open.

\proclaim{Problem 1.6} Can the $n$-sphere $\Sb^n$, $n\geq4$, be
embedded in $F_n(X)$, where $X$ is a curve?
\endproclaim

 \head 2. Embeddability into products of curves and the
first cohomology
\endhead

This chapter splits into sections 2A and 2B. In section 2A we
consider mappings from and into $n$-dimensional products of
compacta $Y_1,\cdots,Y_k$, where $n=\dim Y_1+\cdots + \dim Y_k$,
which induce non-trivial homomorphisms of cohomology in dimension
$n$. Theorem 2A.1 states one of principal observations of this
paper. It says that there is no proper subgroup of
$H^n(Y_1\times\cdots\times Y_k;G)$ which contains all images $\im
H^n(g)$, where $g$ is a mapping from $Y_1\times\cdots\times Y_k$
to the product of spheres with corresponding dimensions. This
implies Corollary 2A.2 on the existence of algebraically essential
mappings. Next we obtain an important Corollary 2A.3 which
provides a relation between algebraically essential mappings into
$Y_1\times\cdots\times Y_k$ and analogous mappings into
corresponding products of spheres. In section 2B we study
analogous problems in which the products are more restricted (they
involve 1-dimensional factors). In that case we get a strong
additional information on the first cohomology of the domain space
of a mapping, see Theorem 2B.1. This enables us to give an
alternative argument for the Kuperberg Theorem 1.3.\ms

\centerline{2A. {\it Mapping products of compacta into
corresponding products of spheres }}\ms

To present our results we first recall some standard notation and
auxiliary results.

For the $n$-sphere $\Sb^n$, $n\ge1$, by $\gamma_n$ we denote a
generator of $H^n(\Sb^n)$, and by $\gamma_{n,{\star}}$ -- the
generator of $H^n(\Sb^n,{\star})$ corresponding to $\gamma_n$
under the inclusion mapping $j:\Sb^n\to(\Sb^n,{\star})$, that is
$H^n(j)(\gamma_{n,{\star}})=\gamma_n$.

To simplify notation we often use the letter $\nu$ to denote the
homomorphisms $H^{m}(X)\otimes{G}\to{H}^{m}(X;G)$ from the
universal coefficient formula (for various $X$ and $G$); and the
letter $\mu$ to denote the homomorphism $H^{n_1}(X_1) \otimes
\cdots \otimes H^{n_k}(X_k) \to H^{n}(X_1\times\cdots\times{X_k})$
from the K\"unneth formula \cite{Sp, p. 237}, where
$n=n_1+\cdots+n_k$. The composition
$$\lambda=\nu\circ(\mu\otimes1_G):
H^{n_1}(X_1)\otimes\cdots\otimes{}H^{n_k}(X_k)\otimes{G} \to
H^{n}(X_1\times\cdots\times{X_k};G),$$ where
$\nu:H^{n}(X_1\times\cdots\times{X_k})\otimes{G}\to
H^{n}(X_1\times\cdots\times{X_k};G)$ in this case, will be called
{\it splitting homomorphism} for
$H^{n}(X_1\times\cdots\times{X_k};G)$. This homomorphism is
natural as both $\mu$ and $\nu$ are known to be natural.

In case where $X_1,\cdots,X_k$ are finite dimensional compacta and
$n_i=\dim{X_i}$ for each $i$, $\lambda$ is an isomorphism. In
fact, the K\"unneth formula for the product
$X_1\times\cdots\times{X}_k$ and integer coefficients is the
following short exact sequence
$$0\to [H^{\ast}(X_1)\otimes\cdots\otimes{H^{\ast}(X_k)}]^n @>{\mu}>>
H^n(X_1\times\cdots\times{X_k})\to
[H^{\ast}(X_1)*\cdots*H^{\ast}(X_k)]^{n+1}\to0.$$ By our
assumption $H^l(X_i)=0$ for $l>\dim{X_i}$. So the formula takes on
the form
$$0\to{}H^{n_1}(X_1)\otimes\cdots\otimes{H^{n_k}(X_k)} @>{\mu}>>
H^n(X_1\times\cdots\times{X_k})\to0\to0 .$$ Thus $\mu$ is an
isomorphism. From the universal coefficient formula we have the
following short exact sequence
$$0\to{}H^{n}(X_1\times\cdots\times{X_k})\otimes{G}@>{\nu}>>
H^{n}(X_1\times\cdots\times{X_k};G) \to
H^{n+1}(X_1\times\cdots\times{X_k})*{G}\ .$$ By our assumption
$\dim(X_1\times\cdots\times{X_k})\le{n}$, so
$H^{n+1}(X_1\times\cdots\times{X_k})=0$. It follows that $\nu$ is
an isomorphism. Consequently $\lambda$ is an isomorphism.

\proclaim{Theorem 2A.1} Let $Y_1,\cdots,Y_k$ be finite dimensional
compacta, let $G$ be an Abelian group, and let $H$ be a proper
subgroup of $H^{n}(Y_1\times\cdots\times{Y}_k;G)$, where $n=n_1 +
\cdots + n_k$ and $n_i=\dim Y_i \ge1$. Then there exist mappings
$\varphi_1:Y_1\to\Sb^{n_1}$, $\cdots$,
$\varphi_k:Y_k\to{}\Sb^{n_k}$ such that the image
$\im{H^n}(\varphi_1\times\cdots\times{}\varphi_k;G)\nsubseteq{H}$.

\endproclaim

\demo{First proof} Under our assumptions the splitting
homomorphism $\lambda$ is an epimorphism (in fact -- an
isomorphism), $\lambda^{-1}(H)$ is a proper subgroup of
$H^{n_1}(Y_1)\otimes\cdots\otimes{H_{n_k}}({Y}_k)\otimes{G}$. As
elements of the form $a_1\otimes\cdots\otimes{}a_k\otimes{}g$,
where $a_1\in{}H^{n_1}(Y_1)$, $\cdots$, $a_k\in{}H^{n_k}(Y_k)$ and
$g\in{}G$, generate
$H^{n_1}(Y_1)\otimes\cdots\otimes{H^{n_k}}({Y}_k)\otimes{G}$, one
of them -- say $a_1\otimes\cdots\otimes{}a_k\otimes{}g$ -- is not
in $\lambda^{-1}(H)$. Hence

\itemitem{(i)}
$\lambda(a_1\otimes\cdots\otimes{}a_k\otimes{}g))\notin{H}$.

\noindent Since $\dim{Y_i}=n_i$, by the Hopf-Whitney
Classification Theorem, cf. \cite{Sp, p. 431},  there exist
mappings $\varphi_1:Y_1\to{}\Sb^{n_1}$, $\cdots$,
$\varphi_k:Y_k\to{}\Sb^{n_k}$ such that
\itemitem{(ii)} $a_1=\varphi_1^{\ast}(\gamma_{n_1})$, $\cdots$,
$a_k=\varphi_k^{\ast}(\gamma_{n_k})$.

\noindent By naturality of the splitting isomorphism, the
following diagram commutes
$$
\CD H^n(Y_1\times\cdots\times{Y_k};G)
 @<{H^n(\varphi_1\times\cdots\times\varphi_k;G)}<<
H^n(\Sb^{n_1}\times\cdots\times{\Sb^{n_k}});G) \\
   @A{\lambda}AA
     @A{\lambda'}AA \\
H^{n_1}(Y_1)\otimes\cdots\otimes{}H^{n_k}(Y_k)\otimes{G}
@<{\varphi_1^{\ast}\otimes\cdots\otimes\varphi_k^{\ast}\otimes{\tenbf1}_G}<<
H^{n_1}(\Sb^{n_1})\otimes\cdots\otimes{}H^{n_k}(\Sb^{n_k})
\otimes{G} \ .
\endCD
$$
Hence, by (ii), we infer that
$$\lambda(a_1\otimes\cdots\otimes{}a_k\otimes{g})=
\lambda(\varphi_1^{\ast}(\gamma_{n_1})\otimes\cdots\otimes
\varphi_k^{\ast}(\gamma_{n_k})\otimes{g})=$$
$$\lambda((\varphi_1^{\ast}\otimes\cdots\otimes\varphi_k^{\ast}\otimes{\tenbf1}_G)
(\gamma_{n_1}\otimes\cdots\otimes\gamma_{n_k}\otimes{}g))=
 H^n(\varphi_1\times\cdots\times\varphi_k;G)
(\lambda'(\gamma_{n_1}\otimes\cdots\otimes\gamma_{n_k}\otimes{}g))
\ .
$$
Thus, by (i), this completes the proof. \qed
\enddemo

We are going to supply yet another argument for Theorem 2A.1 using
(instead of the Hopf-Whitney theorem applied in the above proof)
the following standard fact from cohomology theory:

{\it For any $n$-dimensional polyhedron $P$ $($with a fixed
triangulation$)$ the group $H^n(P)$ is generated by elements
$\sigma^*=H^n(v_{\sigma})(\gamma_{n})$, where $\sigma$ runs over
all $n$-simplices of $P$, and $v_{\sigma}:P\to\Sb^n$ is a mapping
transforming the interior of $\sigma$ homeomorphically onto
$\Sb^n\setminus\{{\star}\}$ and carrying each other point to
${\star}$.}

To clarify the claim let us make the following observations. We
begin letting $j:P\to(P,P^{(n-1)})$ and
$j':\Sb^n\to(\Sb^n,{\star})$ denote the inclusions. By
$\tilde{v}_{\sigma}:(P,P^{(n-1)})\to(\Sb^n,{\star})$ we denote the
mapping determined by $v_{\sigma}$, i.e.
$\tilde{v_{\sigma}}(x)=v_{\sigma}(x)$ for each $x\in{P}$. Since
$j'\circ{v}_{\sigma}=\tilde{v}_{\sigma}\circ{j}$, we have the
equalities:
$$\sigma^*=H^n(v_{\sigma})(\gamma_{n})=
H^n(v_{\sigma})(H^n(j')(\gamma_{n,{\star}}))=
H^n(j'\circ{v}_{\sigma})(\gamma_{n,{\star}})=$$
$$H^n(\tilde{v}_{\sigma}\circ{j})(\gamma_{n,{\star}})=
H^n(j)(H^n(\tilde{v}_{\sigma})(\gamma_{n,{\star}})) \ .$$ It is
known that the group $H^n(P,P^{(n-1)})$ is freely generated by the
elements $H^n(\tilde{v}_{\sigma})(\gamma_{n,{\star}})$. Since the
inclusion $j$ induces an epimorphism
$H^n(j):H^n(P,P^{(n-1)})\to{H^n}(P)$, the claim follows.

 \demo{Second proof} By the Freudenthal theorem each $Y_i$ is
the limit of an inverse sequence
$Y_{i,1}\leftarrow{Y_{i,2}}\leftarrow\cdots$ of $n_i$-dimensional
polyhedra. (Each $Y_{i,m}$ is taken with a fixed triangulation.)
Let $p_{i,m}:Y_i\rightarrow{Y_{i,m}}$ denote the projections. Then
$Y_1\times\cdots\times{Y_k}$ is the limit of the sequence
$$Y_{1,1}\times\cdots\times{Y_{k,1}}\leftarrow
{Y_{1,2}\times\cdots\times{Y_{k,2}}}\leftarrow\cdots$$ and
$p_{1,m}\times\cdots\times{p_{k,m}}:Y_{1}\times\cdots\times{Y_{k}}
\rightarrow{Y_{1,m}\times\cdots\times{Y_{k,m}}}$ are the
projections. By the continuity of the {\v C}ech cohomology there
exist an index $m$ such that
$$H'=(H^n(p_{1,m}\times\cdots\times{p_{k,m}};G))^{-1}(H)$$
is a proper subgroup of
$H^n(Y_{1,m}\times\cdots\times{Y_{k,m}};G)$. It follows that
$(\lambda)^{-1}(H')$ is a proper subgroup of
$H^{n_1}(Y_{1,m})\otimes\cdots\otimes{H^{n_k}(Y_{k,m})}\otimes{G}$,
where $\lambda$ is the splitting isomorphism for
$H^n(Y_{1,m}\times\cdots\times{Y_{k,m}};G)$. Referring to the
claim preceding the proof we infer that
$H^{n_1}(Y_{1,m})\otimes\cdots\otimes{H^{n_k}(Y_{k,m})}\otimes{G}$
is generated by elements of the form
$\sigma_1^*\otimes\cdots\otimes\sigma_k^*\otimes{g}$, where
$\sigma_i$ is an $n_i$-simplex of $Y_{i,m}$ and $g\in{G}$. It
follows that there exists an element
$\sigma_1^*\otimes\cdots\otimes\sigma_k^*\otimes{g}
\notin(\lambda)^{-1}(H')$,where

\itemitem{(ii$'$)}
$\sigma_i^*=H^{n_i}(v_i)(\gamma_{n_i})$

\noindent and $v_i:Y_{i,m}\to\Sb^{n_i}$ is the mapping
corresponding to $\sigma_i$. Thus

\itemitem{(i$'$)}
$\lambda(\sigma_1^*\otimes\cdots\otimes\sigma_k^*\otimes{g})
\notin{H'}.$

On the other hand, by an argument similar to that used in the
first proof, we get the equality
$$\lambda(\sigma_1^*\otimes\cdots\otimes\sigma_k^*\otimes{g})=
H^n(v_1\times\cdots\times{v_k};G)
(\lambda'(\gamma_{n_1}\otimes\cdots\otimes\gamma_{n_k}\otimes{}g))
\ .$$ Setting $\varphi_i=v_i\circ{p_{i,m}}:Y_i\to\Sb^{n_i}$ for
$i=1,\cdots,k$, by (i$'$), we easily infer the conclusion.
 \qed
\enddemo\ms

Theorem 2A.1 readily implies the following corollary on existence
of algebraically essential mappings.

\proclaim{Corollary 2A.2} Let $Y_1,\cdots,Y_k$ be finite
dimensional compacta such that
$H^{n}(Y_1\times\cdots\times{Y}_k;G)\neq0$, where $n=n_1 + \cdots
+ n_k$ and $n_i = \dim Y_i\ge1$. Then there exist mappings
$\varphi_1:Y_1\to\Sb^{n_1}$, $\cdots$,
$\varphi_k:Y_k\to{}\Sb^{n_k}$ such that
$H^n(\varphi_1\times\cdots\times{}\varphi_k;G)\neq0$. In
particular,
$\varphi_1\times\cdots\times{}\varphi_k:Y_1\times\cdots\times{Y}_k
\to \Sb^{n_1}\times\cdots \times\Sb^{n_k}$ is algebraically
essential.
\endproclaim

\demo{Proof} By our hypothesis the trivial group $H=0$ is a proper
subgroup of $H^{n}(Y_1\times\cdots\times{Y}_k;G)$. Hence the
conclusion follows from Theorem 2A.1.\qed\enddemo

\ms Our next corollary points out an important relation between
algebraically non-trivial mappings into products of finite
dimensional compacta and algebraically essential mappings into
products of spheres.

\proclaim{Corollary 2A.3} Let $f:X\to Y_1\times\cdots\times Y_k$
be a mapping of a compactum $X$ to the product of finite
dimensional compacta $Y_1, \cdots , Y_k$. If f is non-trivial with
respect to $H^n(\cdot\,;\,G)$, where $n=n_1 + \cdots + n_k$ and
$n_i=\dim Y_i \ge1$, then there exist mappings
$\varphi_1:Y_1\to\Sb^{n_1}$, $\cdots$,
$\varphi_k:Y_k\to{}\Sb^{n_k}$ such that
$$(\varphi_1\times\cdots\times{}\varphi_k)\circ{f}:X \to \Sb^{n_1}\times\cdots
\times\Sb^{n_k}$$ is non-trivial with respect to
$H^n(\cdot\,;\,G)$ as well.
\endproclaim

\demo{Proof} It follows from our hypothesis that the group $H=\ker
f^{\ast}$ is a proper subgroup of $H^{n}(Y_1 \times \cdots \times
Y_k;G)$. Hence by Theorem 2A.1 there exist mappings
$\varphi_1:Y_1\to\Sb^{n_1}$, $\cdots$,
$\varphi_k:Y_k\to{}\Sb^{n_k}$ such that
$\im{H^n}(\varphi_1\times\cdots\times{}\varphi_k;G)\nsubseteq{H}$.
One easily sees that these mappings satisfy the conclusion of our
corollary. \qed\enddemo\ms

Next corollary immediately follows from the above result (applied
to an inclusion) taking into account the following simple lemma
which has been observed in \cite{D-K}.

\proclaim{Lemma~2A.4} Let $X$ be a closed subset of an
$n$-dimensional compactum $Y$, and let $i:X\hookrightarrow{}Y$
denote the inclusion mapping. Then $H^n(i\,;\,G):H^n(Y\,;\,G)\to
H^n(X\,;\,G)$ is an epimorphism. In particular, $H^n(i\,;\,G)\ne0$
if $H^n(X\,;\,G)\ne0$.
\endproclaim

\demo{Proof} Let us consider the following portion of the
cohomology exact sequence of the pair $(Y,X)$:
$$H^n(Y;G) @>{i^{\ast}}>> H^n(X;G) \to
H^{n+1}(Y,X;G) \ .$$ Note that $H^{n+1}(Y,X;G)=0$ as $\dim{}Y=n$.
It follows that $H^n(i\,;\,G)$ is an epimorphism.
 \qed
\enddemo

\proclaim{Corollary 2A.5} Let $X$ be a compactum lying in the
product $Y_1\times\cdots\times Y_k$ of finite dimensional
compacta. If $H^n(X\,;\,G)\neq 0$, where $n=n_1 + \cdots + n_k$
and $n_i=\dim Y_i \ge1$, then there exist mappings
$\varphi_1:Y_1\to\Sb^{n_1}$, $\cdots$,
$\varphi_k:Y_k\to{}\Sb^{n_k}$ such that
$H^n((\varphi_1\times\cdots\times{}\varphi_k)|X;G)\neq 0$. In
particular, $X$ admits algebraically essential mappings into
$\Sb^{n_1}\times\cdots \times\Sb^{n_k}$.\qed
\endproclaim

\ms \centerline{2B. {\it Mappings into products involving
$1$-dimensional factors }}\bs

In case of mappings into products of curves Corollary 2A.3,
combined with some other facts, gives an important information on
the first cohomology of the domain space.

\proclaim{Theorem 2B.1} Let
$f:X\to{}Y_1\times\cdots\times{}Y_k\times{}Y$ be a map of a
compactum $X$ into the product of $1$-dimensional compacta
$Y_1,\cdots,Y_k$, and an $l$-dimensional compactum $Y$. If
$H^{k+l}(f)\neq 0$ then $\rank H^1(X)\ge{k}$.
\endproclaim

The actual proof will be given after Lemmas 2B.2 and 2B.3 below.
In the next lemma we use the following well-known fact (which
follows from a result of M. Brown \cite{Br, p. 94}, see Chapter
1):\ss

{\it If $X$ is a compactum and $M$ is a closed manifold, then any
mapping $f:X \to M$ which is not surjective is homotopic to a
mapping $g$ such that $\dim{g(X)}\le\dim{M}-1$.}\ss

\proclaim{Lemma 2B.2} Let $X$, $Y$ be compacta, let $M$ be a
closed manifold, and let $n=\dim{M}+\dim{Y}$. If
$(f,g):X\to{M}\times{Y}$ is a mapping such that $H^n((f,g);G)$ is
non-trivial, for an Abelian group $G$, then $f$ is essential.
\endproclaim

\demo{Proof} Otherwise, by the result preceding this lemma,
$(f,g)$ is homotopic to a map $h:X\to{M}\times{Y}$ such that
$\dim{h(X)}\le{n-1}$. Then $H^n((f,g);G)=H^n(h;G)=0$, a
contradiction. \qed
\enddemo

Let $A$ be a subset of a space $X$ and let $i:A\hookrightarrow{X}$
be the inclusion. For $a\in{}H^n(X)$, we write briefly $a|A$ for
$i^{\ast}(a)$. Note that for a mapping $f:X\to{Y}$ and
$b\in{}H^n(Y)$ we have $f^{\ast}(b)|A=(f|A)^{\ast}(b)$.

\proclaim{Lemma 2B.3} Let $X$ be a compactum, and let
$f_i:X\to\Sb^1$, $i=1\cdots,k$, be mappings such that
$(f_1,\cdots,f_k):X\to\Sb^1\times\cdots\times\Sb^1$ is essential.
Then elements
$$f_1^{\ast}(\gamma_1),\cdots,f_k^{\ast}(\gamma_1)\in H^1(X)$$
are linearly independent.
\endproclaim

\demo{Proof} Let $\e=(1,0)$ be the unit point of $\Sb^1$. Put
$$X_i=\{x\in{X}:f_j(x)=\e\ {\tenrm for}\ j\ne{i}\}\ .$$
We have a
natural homomorphism
$$\eta:H^1(X)\to  H^1(X_1)\oplus\cdots\oplus H^1(X_k)$$
defined by the formula $\eta(a)=(a|X_1,\cdots,a|X_k)$ for
$a\in{}H^1(X)$. Since $f_i(X_j)\subset\{\e\}$ for $i\ne{j}$, and
$f_i^{\ast}(\gamma_1)|X_j=(f_i|X_j)^{\ast}(\gamma_1)$, we infer
that

\itemitem{(i)} $\eta(f^{\ast}_i(\gamma_1))=(0,\cdots,0,(f_i|X_i)^{\ast}(\gamma_1),0,
\cdots,0)$,

\noi where $(f_i|X_i)^{\ast}(\gamma_1)$ stands in the $i$th
position.

The conclusion of our lemma will follow once we show the images
$\eta(f^{\ast}_1(\gamma_1))$, $\cdots$,
$\eta(f^{\ast}_k(\gamma_1))$ are linearly independent. To this
end, by (i), it is enough to show that each
$(f_i|X_i)^{\ast}(\gamma_1)$ is non-zero (here we refer to
$H^1(Y)$ being torsion free for any space $Y$). But
$(f_i|X_i)^{\ast}(\gamma_1)\ne0$ if and only if
$f_i|X_i:X_i\to\Sb^1$ is essential (this follows from the
Bruschlinsky theorem). Hence we have to prove that $f_i|X_i$ is
essential. Suppose it is not true.

Then there is a homotopy $F:X_i\times{I}\to\Sb^1$ such that
$F_0=f_i|X_i$ and $\e\notin{}F_1(X_i)$. As $f_i|X_i$ extends to
$f_i:X\to\Sb^1$ there is an extension $F':X\times{I}\to\Sb^1$ of
$F$ (by the homotopy extension theorem). Now we define a homotopy
$G:X\times{I}\to\Sb^1\times\cdots\times\Sb^1$ by the formula
$$G(x,t)=(f_1(x),\cdots,f_{i-1}(x),F'(x,t),f_{i+1}(x),\cdots,f_{k}(x))\ .$$
Obviously, $G$ is continuous and $G_0=f$. To complete the proof it
suffices to show that

\itemitem{(ii)} $(\e,\cdots,\e)\notin G_1(X)$ .

\noi (In fact, this contradicts the essentiality of $f$.)

So fix a point $x\in X$. We have to show that
$G_1(x)\ne(\e,\cdots,\e)$. This clearly holds if $f_j(x)\ne\e$ for
some $j\ne i$. Then consider the case $f_j(x)=\e$ for all $j\ne
i$. In such a case $x\in X_i$. Therefore
$$G_1(x)=(\e,\cdots,\e,F'(x,1),\e,\cdots,\e)\ne(\e,\cdots\,e)\ ,$$
because $F'(x,1)=F(x,1)\ne\e$. This proves (ii), and completes the
proof of our lemma. \qed
\enddemo

\centerline{\bf Proof of Theorem~2B.1}

\bs

By Corollary~2A.3 there is a mapping
$(\psi_1,\cdots,\psi_k,\psi):X\to\Sb^1\times\cdots\times\Sb^1\times\Sb^{l}$
such that $H^{k+l}((\psi_1,\cdots,\psi_k,\psi))\ne0$. By
Lemma~2B.2 we infer that the mapping
$(\psi_1,\cdots,\psi_k):X\to\Sb^1\times\cdots\times\Sb^1$ is
essential. Applying Lemma~2B.3 we conclude that
$\psi_1^{\ast}(\gamma_1),\cdots,\psi_k^{\ast}(\gamma_1)$ are
linearly independent in $H^1(X)$. This ends the proof. \qed\bs

\proclaim{Corollary 2B.4} Let $M$ be a closed $(k+l)$-manifold
lying in the product $Y_1\times\cdots\times{}Y_k\times{}Y$  of
$1$-dimensional compacta $Y_1,\cdots,Y_k$, and an $l$-dimensional
compactum $Y$. Then $\rank{}H^1(M)\ge{}k$.
\endproclaim

\demo{Proof} We are going to apply Theorem~2B.1. Notice that
$\dim(Y_1\times\cdots\times{}Y_k\times{}Y)=k+l$. Thus, by
Lemma~2A.3, $H^{k+l}(i)$ is not trivial, where
$i:M\to{}Y_1\times\cdots\times{}Y_k\times{}Y$ denotes the
inclusion. Applying Theorem~2B.1, we get the conclusion. \qed
\enddemo

Now we can give an alternative proof of Kuperberg's theorem \cite
{Kup}.

\proclaim{Corollary 2B.5 (Kuperberg)} Let $M$ be a closed
$(l+1)$-manifold with finite fundamental group $\pi_1(M)$. Then
$M$ is not embeddable in any product $Y_1\times{}Y$, where $Y_1$
is a $1$-dimensional compactum  and $Y$ is an $l$-dimensional
compactum.
\endproclaim

\demo{Proof} Suppose $M$ embeds in $Y_1\times{}Y$, where $Y_1$ and
$Y$ are as above. By Corollary~2B.4, it follows that $H^1(M)$ is
non-trivial. On the other hand, it follows from our assumption
that $H_1(M)$ is finite (being abelisation of the finite group
$\pi_1(M)$, cf. \cite{Sp, pp. 398, 391}). By the theorem on
universal coefficients for cohomology we have the following exact
sequence:
$$0\to\Ext(H_0(M),\Z)\to{}H^1(M)\to\Hom(H_1(M),\Z)\to0 \ .$$
As $H_1(M)$ is finite, the last term is trivial. Since
$H_0(M)\approx\Z$ is free, $\Ext(H_0(M),\Z)$ is trivial as well
(cf. \cite{Sp, 5.5.1, p. 241}). Hence $H^1(M)=0$, a contradiction.
\qed
\enddemo

\head 3. Embeddability into products and the cohomology ring
\endhead

This chapter bas been divided in sections 3A and 3B. In section 3A
we show that algebraically non-trivial mappings into products are
strictly related to certain non-zero cup products in the domain
space. The main result of the entire chapter is Theorem~3A.1.
Corollary 3A.4 generalizes the Dydak-Koyama Theorem 1.5. In
section 3B we apply those results to a study of the categories of
spaces.\ms

\centerline{3A. {\it Algebraically non-trivial mappings into
products of compacta and}} \centerline{\it non-zero cup products
of spherical elements}\ms

In this section we show that mappings into products
$Y_1\times\cdots\times Y_k$ are non-trivial with respect to
$H^n(\cdot\,;\,\bigotimes_{i=1}^{k} G_i)$, where $n=\dim
Y_1+\cdots+\dim Y _k$, if and only if they induce certain non-zero
cup products in the domain space. We need the following notion. An
element $\alpha\in H^m(Y;G)$ is called {\it spherical} if
$\alpha=\nu(\varphi^{\ast}(\gamma_m)\otimes{g})$ for a mapping
$\varphi:Y\to \Sb^m$ and an element $g\in{G}$.

\proclaim{ Theorem 3A.1} Let $f=(f_1,\cdots,f_k): X \to Y_1 \times
\cdots \times Y_k$ be a mapping of a compactum $X$ in the product
of finite dimensional compacta. Then $
H^{n}(f\,;\,\bigotimes_{i=1}^{k}{G_i})\neq 0$, where
$n=n_1+\cdots+n_k$ and $n_i=\dim Y_i\ge1$, if and only if there
exist spherical elements $\alpha_1\in{H^{n_1}}(Y_1;G_1)$,
$\cdots$, $\alpha_k\in{H^{n_k}}(Y_k;G_k)$ such that the cup
product $f_1^{\ast}(\alpha_1)\ssm\cdots\ssm{f_k^{\ast}(\alpha_k)}
\in{}H^n(X;\bigotimes_{i=1}^{k}{G_i})$ is not zero.

In particular, if $X\subset{Y_1}\times\cdots\times{Y_k}$ and
$H^{n}(X;\bigotimes_{i=1}^{k}{G_i})\ne0$  then there exist
elements $a_1\in{H^{n_1}}(X;G_1)$, $\cdots$ ,
$a_k\in{H^{n_k}}(X;G_k)$ such that the cup product
$a_1\ssm\cdots\ssm{a_k}\in{}H^n(X;\bigotimes_{i=1}^{k}{G_i})$ is
not zero. $($So, $H^{n_i}(X;G_i)\ne0$ for each $i$.$)$
\endproclaim

\noi{\bf Remark.} A simpler variant of this theorem is Corollary
3A.6, where one assumes that each $G_i=G$ and that $H^n(f;G)\neq
0$ for some $G$. \ms

We are going to show that Theorem 3A.1 follows from Corollary 2A.4
and Lemma 3A.3 below. A short proof will be given after the proof
of Lemma~3A.3. In the proof of that lemma we use Lemma~3A.2 which
will be proved first. The general K\"unneth theorem for the \v
Cech cohomology of finite products of compact spaces imposes some
limitations on the coefficient groups. In the following lemma we
point out that those limitations are immaterial in case of
additional dimensional restrictions.

\proclaim{Lemma 3A.2} Let $X_1, \cdots , X_k$ be nonempty compacta
with $\dim X_i = n_i$ and let $n=n_1 + \cdots + n_k$. Then the
homomorphism
$$\mu': H^{n_1}(X_1;G_1) \otimes \cdots \otimes H^{n_k}(X_k;G_k)
\to H^{n}(X_1\times\cdots\times{X_k};\bigotimes_{i=1}^k{G_i})$$ is
an isomorphism.
\endproclaim

\demo{Proof} The introductory remarks preceding Theorem 2A.1 show
that
$$\mu:H^{n_1}(X_1)\otimes\cdots\otimes{H^{n_k}(X_k)}\to
H^n(X_1\times\cdots\times{X_k})$$ is an isomorphism.

From the universal coefficient formula we have the exact sequence
$$0\to{}H^{n_i}(X_i)\otimes{G_i}@>{\nu_i}>>H^{n_i}(X_i;G_i)\to
H^{n_i+1}(X_i)*{G_i}\ .$$ As $H^{n_i+1}(X_i)=0$, each $\nu_i$ is
an isomorphism. Hence
$$
(H^{n_1}(X_1)\otimes{G_1})\otimes\cdots\otimes(H^{n_k}(X_k)\otimes{G_k})
@>{\nu_1\otimes\cdots\otimes\nu_k}>>
H^{n_1}(X_1;G_1)\otimes\cdots\otimes{H^{n_k}}(X_k;G_k)
$$
is an isomorphism as well. Let $\varphi$ be a composition of the
canonical isomorphism
$$(H^{n_1}(X_1)\otimes\cdots\otimes{H^{n_k}(X_k)})
\otimes\bigotimes_{i=1}^k{G_i} \to
(H^{n_1}(X_1)\otimes{G_1})\otimes\cdots\otimes(H^{n_k}(X_k)\otimes{G_k})
$$
and $\nu_1\otimes\cdots\otimes\nu_k$. Let
$$
\nu:H^{n}(X_1\times\cdots\times{X_k})\otimes(G_1\otimes\cdots\otimes{G_k})
\to H^{n}(X_1\times\cdots\times{X_k};G_1\otimes\cdots\otimes{G_k})
$$
be the homomorphism from the universal coefficient formula. Since
$\dim(X_1\times\cdots\times{X_k})=n$, we have
$H^{n+1}(X_1\times\cdots\times{X_k})=0$. It follows, as before,
that $\nu$ is an isomorphism.

Let us consider the following diagram:
$$
\CD (H^{n_1}(X_1)\otimes\cdots\otimes{H^{n_k}(X_k)})
\otimes\bigotimes_{i=1}^k{G_i} @>{\mu\otimes{\tenbf 1}}>>
  H^n(X_1\times\cdots\times X_k)\otimes\bigotimes_{i=1}^k{G_i}\\
   @V{\varphi}VV
     @V{\nu}VV \\
  H^{n_1}(X_1;G_1)\otimes\cdots\otimes{H^{n_k}(X_k;G_k)} @>{\mu'}>>
       H^n(X_1\times\cdots\times{X_k};\bigotimes_{i=1}^k{G_i}),
\endCD
$$
where ${\tenbf 1}
:\bigotimes_{i=1}^k{G_i}\to\bigotimes_{i=1}^k{G_i}$ is the
identity isomorphism. By the above considerations we infer that
all the upper horizontal and the vertical homomorphisms are
isomorphisms. To complete the proof it is enough to observe that
the diagram commutes. \qed
\enddemo

\proclaim{Lemma 3A.3} Let $X$ be a compactum, let $G_1,\cdots,G_k$
be a sequence of Abelian groups, and let $n_1,\cdots, n_k$,
$n_i\geq1$, be natural numbers. Then any mapping
$f=(f_1,\cdots,f_k):X\to{\Sb^{n_1}}\times\cdots\times{\Sb^{n_k}}$
is non-trivial with respect to
$H^{n}(\cdot\,;\,\bigotimes_{i=1}^{k}{G_i})$, where
$n=n_1+\cdots+n_k$, if and only if there exist elements
$\alpha_1\in{H^{n_1}(\Sb^{n_1};G_1)}$, $\cdots$,
$\alpha_k\in{H^{n_k}}(\Sb^{n_k};G_k)$ such that the cup product
$$f_1^{\ast}(\alpha_1)\ssm\cdots \ssm f_k^{\ast}(\alpha_k)\in
H^n(X;\bigotimes_{i=1}^{k}{G_i})$$ is  non-zero. In the "only if"
part we can take $\alpha_i=\nu(\gamma_{n_i}\otimes{g}_i)$ for some
$g_i\in G_i$, and consequently,
$f_i^{\ast}(\alpha_i)=\nu(f_i^{\ast}(\gamma_{n_i})\otimes g_i)$;
moreover, if $G_i=\Z$, then we can take $\alpha_i=\gamma_i$ $($and
$g_i=1$$)$.
\endproclaim

\demo{Proof} We start with the "only if" implication. Since
$\dim(\Sb^{n_1}\times\cdots\times{\Sb^{n_k}})=n$, the homomorphism
$$\mu:H^{n_1}(\Sb^{n_1};G_1)\otimes\cdots\otimes{H^{n_k}}(\Sb^{n_k};G_k)\to
H^n(\Sb^{n_1}\times\cdots\times\Sb^{n_k};\bigotimes_{i=1}^{k}{G_i})$$
(from the K\"unneth formula) is an epimorphism, see Lemma~3.2. As
$f^{\ast}(=H^n(f;\bigotimes_{i=1}^{k}{G_i}))$ is not trivial, it
follows that $f^{\ast}\circ\mu$ is neither. Hence there exist
elements $\alpha_1\in{H^{n_1}(\Sb^{n_1};G_1)}$, $\cdots$,
$\alpha_k\in{H^{n_k}}(\Sb^{n_k};G_k)$, such that

\ss
\itemitem{(i)} $(f^{\ast}\circ\mu)(\alpha_1\otimes\cdots\otimes{\alpha_k})\ne0$.
\ss

\noindent(Moreover we can take $\alpha_i=\gamma_i$ if $G_i=\Z$.)

By definition, $\mu(\alpha_1\otimes\cdots\otimes{\alpha_k})=
\alpha_1\times\cdots\times{\alpha_k}$. Thus referring to \cite{Sp,
Corollary~5.6.14, p. 253} we get
$\alpha_1\times\cdots\times{\alpha_k}=
p_1^{\ast}(\alpha_1)\ssm\cdots\ssm{p_k}^{\ast}(\alpha_k)$, where
$p_i:\Sb^{n_1}\times\cdots\times{\Sb^{n_k}}\to{\Sb^{n_i}}$ is the
projection for each $i$. Combining these equalities and \cite{Sp,
Property 5.6.8, p. 251}, we get

\ss
\itemitem{(ii)}
$(f^{\ast}\circ\mu)(\alpha_1\otimes\cdots\otimes{\alpha_k})=
f^{\ast}(p_1^{\ast}(\alpha_1)\ssm\cdots\ssm{p}_k^{\ast}(\alpha_k))=
f^{\ast}(p_1^{\ast}(\alpha_1))\ssm\cdots\ssm{f^{\ast}(p_k^{\ast}(\alpha_k))}=
f_1^{\ast}(\alpha_1)\ssm\cdots\ssm{f_k^{\ast}(\alpha_k)}$. \ss

 This combined
with (i) completes the proof of this implication. Now let us prove
the additional claim.

Consider the following commutative diagram:
$$
\CD H^{n_i}(\Sb^{n_i})\otimes{G_i} @>{f_i^{\ast}\otimes{\tenbf
1}}>>
  H^{n_i}(X)\otimes{G_i}\\
   @V{\nu}VV
     @V{\nu}VV \\
  H^{n_i}(\Sb^{n_i};G_i) @>{f_i^{\ast}}>>
       H^{n_i}(X;G_i) \ .
\endCD
$$
(In the top arrow $f_i^{\ast}=H^{n_i}(f_i)$, while the bottom
arrow is $f_i^{\ast}=H^{n_i}(f_i;G_i)$.) As each element of the
group $H^{n_i}(\Sb^{n_i})\otimes{G_i}$ ($= \, <\gamma_{n_i}>
\otimes\,{G_i}$) can be uniquely written in the form
$\gamma_{n_i}\otimes g$ (for $g\in G_i$) and
$\nu:{H}^{n_i}(\Sb^{n_i})\otimes{G_i}\to H^{n_i}(\Sb^{n_i};G_i)$
is an epimorphism, there exists $g_i\in G_i$ such that
$\alpha_i=\nu(\gamma_{n_i}\otimes{g_i})$. (For $G_i=\Z$, it
follows from the definition of $\nu$ that
$\nu(\gamma_{n_i}\otimes1)=\gamma_{n_i}$.)  By the commutativity
we get
$$f_i^{\ast}(\alpha_i)=\nu(f_i^{\ast}(\gamma_{n_i})\otimes g_i) \ ,$$
which ends the proof of the claim.

The other implication follows from the equality
$f^{\ast}\mu(\alpha_1\otimes\cdots\otimes{\alpha_k})=
f_1^{\ast}(\alpha_1)\ssm\cdots\ssm{f_k^{\ast}(\alpha_k)}$, see
(ii).
 \qed
\enddemo

\centerline{\bf Proof of Theorem 3A.1} \bs

\noi The proof has been divided in two parts.

{\it Part} "if". By our assumption and Corollary 2A.4 there exist
mappings $\varphi_1:Y_1\to\Sb^{n_1}$, $\cdots$,
$\varphi_k:Y_k\to{}\Sb^{n_k}$ such that
$H^n((\varphi_1\times\cdots\times{}\varphi_k)\circ{f}\,;\,G)\neq
0$. Since $(\varphi_1\times\cdots\times{}\varphi_k)\circ{f}=
(\varphi_1 \circ{f_1},\cdots,\varphi_k\circ{f_k})$, by Lemma 3A.3
there exist elements $\beta_1\in{H^{n_1}(\Sb^{n_1};G_1)}$,
$\cdots$, $\beta_k\in{H^{n_k}}(\Sb^{n_k};G_k)$ such that the cup
product
$$(\varphi_1 \circ{f_1})^{\ast}(\beta_1)\ssm\cdots \ssm (\varphi_k \circ{f_k})^{\ast}(\beta_k) \in
H^n(X;\bigotimes_{i=1}^{k}{G_i})$$ is  non-zero. Then each
$\alpha_i=(\varphi_i)^{\ast}(\beta_i)$ is a spherical element of
$H^{n_i}(Y_i;G_i)$, and we have
$$({f_1})^{\ast}(\alpha_1)\ssm\cdots \ssm ({f_k})^{\ast}(\alpha_k)=
(\varphi_1 \circ{f_1})^{\ast}(\beta_1)\ssm\cdots \ssm (\varphi_k
\circ{f_k})^{\ast}(\beta_k),$$ which proves this part.

{\it Part} "only if". By our assumption there exist elements
$\alpha_1\in{H^{n_1}}(Y_1;G_1)$, $\cdots$,
$\alpha_k\in{H^{n_k}}(Y_k;G_k)$ such that the cup product
$$({f_1})^{\ast}(\alpha_1)\ssm\cdots \ssm ({f_k})^{\ast}(\alpha_k) \in
H^n(X;\bigotimes_{i=1}^{k}{G_i})$$ is  non-zero. Let $p_i:
Y_1\times\cdots\times{Y_k}\to{Y_i}$ denote the projection. Then we
have
$$({f_1})^{\ast}(\alpha_1)\ssm\cdots \ssm ({f_k})^{\ast}(\alpha_k)=
({p_1}\circ{f})^{\ast}(\alpha_1)\ssm\cdots \ssm
({p_k}\circ{f})^{\ast}(\alpha_k)=$$
$$f^{\ast}(({p_1})^{\ast}(\alpha_1)\ssm\cdots
\ssm ({p_k})^{\ast}(\alpha_k)),$$ which completes the proof. \qed

\bs

Theorem~3A.1 implies the following generalization of the
Dydak-Koyama theorem.

\proclaim{Corollary 3A.4} Let $f: X \to Y_1 \times \cdots \times
Y_k$ be a mapping of a compactum $X$ into the product of finite
dimensional compacta, and let $n=n_1+\cdots+n_k$, where $n_i=\dim
Y_i\ge1$. If $H^{n}(f\,;\,G)\neq 0$, then $H^{n_i}(X;G)\ne0$ for
each $i$.

In particular, if $X\subset{Y_1\times\cdots\times{Y_k}}$ and
$H^n(X;G)\ne0$, then $H^{n_i}(X;G)\ne0$ for each $i$.
\endproclaim

\demo{Proof} Without loss of generality we can assume $i=1$. Thus
we have to show that $H^{n_1}(X;G)\ne0$. To this end apply
Theorem~3A.1 with $G_1=G$ and $G_j=\Z$ for $j\neq 1$. Thus we
infer that there exist elements $a_1\in{H^{n_1}}(X;G)$,
$a_2\in{H^{n_2}}(X)$, $\cdots$ , $a_k\in{H^{n_k}}(X)$ such that
$a_1\ssm\cdots\ssm{a_k}\ne0$. Consequently, $a_1\ne0$, which
proves the first statement of our corollary.

To prove the second statement, we  repeat the argument used in the
proof of Corollary~2A.4. \qed
\enddemo\ms

The final result of this section - Corollary 3A.6 - represents a
simplified version of Theorem 3A.1 in case where all groups
$G_1,\cdots,G_k$ are equal and somehow restricted. In the proof we
refer to Lemma 3A.5 which makes the simplification possible. The
lemma shows that the hypothesis $H^{n}(f\,;\,\bigotimes^kG)\neq 0$
from Theorem 3A.1 can be replaced by $H^{n}(f;\,G)\neq 0$ (where
$\bigotimes^kG$ denotes $G\otimes\cdots\otimes{}G$). Lemma 3A.5 is
a consequence of Theorem A1 from the Appendix.

\proclaim{Lemma 3A.5} Let $f:X\to{}Y$ be a mapping between
compacta and let $\dim{}Y=n$. Suppose $G$ is either a direct sum
of cyclic groups or a non-torsion Abelian group. Then
$H^n(f;G)\neq 0$ implies $H^n(f;\bigotimes{}^kG)\neq 0$.
\endproclaim

\demo{Proof} By the universal coefficient theorem the diagram
$$
\CD H^{n}(Y)\otimes{G'} @>{f^{\ast}\otimes1_{G'}}>>
  H^{n}(X)\otimes{G'}\\
   @V{\nu}VV
     @V{\nu}VV \\
H^{n}(Y;G') @>H^n(f;G')>> H^{n}(X;G')
\endCD
$$
commutes for any group $G'$. Since $H^n(f;G)\neq 0$ and the left
vertical homomorphism is an epimorphism (because $\dim Y\le{n}$),
the homomorphism
$$
f^{\ast}\otimes{1}_{G}: H^{n}(Y)\otimes{G} \to
  H^{n}(X)\otimes{G}
$$
is non-trivial (in fact, take $G'=G$ in the diagram). Therefore,
since $G$ is either a direct sum of cyclic groups or a non-torsion
Abelian group, by Theorem A1 from the Appendix the homomorphism
$$
f^{\ast}\otimes1_{\bigotimes^kG}: H^{n}(Y)\otimes{\bigotimes^kG}
\to
  H^{n}(X)\otimes{\bigotimes^kG}
$$
is non-trivial as well. Consequently, since the right vertical
homomorphism is a monomorphism, (take $G'=\bigotimes^kG$ in the
diagram) the homomorphism $H^n(f\,;\,\bigotimes^kG)$ is
non-trivial (in fact, take $G'=\bigotimes^kG$ in the diagram),
which completes the proof.
\enddemo

\proclaim{Corollary 3A.6} Let $f: X \to Y_1 \times \cdots \times
Y_k$ be a mapping of a compactum $X$ into the product of finite
dimensional compacta, and let $n=n_1+\cdots+n_k$, where $n_i=\dim
Y_i\ge1$. Suppose $G$ is either a direct sum of cyclic groups or a
non-torsion Abelian group. Then $H^{n}(f\,;\,G)\neq 0$ implies
that there exist elements $a_1\in{H^{n_1}}(X;G),\cdots ,
a_k\in{H^{n_k}}(X;G)$ such that the cup product
$a_1\ssm\cdots\ssm{a_k}\in{}H^n(X;\bigotimes{}^{k}{G})$ is not
zero .
\endproclaim

\demo{Proof} This corollary directly follows from Theorem 3A.1 and
Lemma 3A.5. \qed
\enddemo\bs

\centerline{3B. {\it Categories of spaces}}

\bs

In this section we are going to show that for any compactum its
embeddability into products is related to its category.

Let us recall the classic definition of the category (cf. [Sp]): A
space $X$ is said to have {\it category} $\le k$ (written: $\cat
X\le k$) if there exists a closed covering $\{F_1,\cdots,F_k\}$ of
$X$ such that each $F_i$ is contractible in $X$. We shall deal
with a modification of this notion.

Let $G_1,\cdots,G_k$ be Abelian groups and let $n_1,\cdots,n_k$ be
natural numbers. A space $X$ is said to have {\it category}
$(G_1,\cdots,G_k;n_1,\cdots,n_k)$ (written: $X\in$

\noindent $\cat(G_1,\cdots,G_k;n_1,\cdots,n_k)$) if there exists a
closed covering $\{F_1,\cdots,F_k\}$ of $X$ such that each
homomorphism
$$  H^{n_i}(F_i\hookrightarrow{}X;G_i)$$
is trivial. In such a case, it is sometimes said that $X$ {\it has
trivial decomposition relatively}
$(G_1,\cdots,G_k;n_1,\cdots,n_k)$. To shorten notation we write
$(G;n_1,\cdots,n_k)$ for $(G,\cdots,G;n_1,\cdots,n_k)$.

Let us recall some basic properties of the cup product which will
be used in the proof of our next lemma.

Let $f:X\to{}Y$ be a mapping and suppose
$f(A_1)\subset{}B_1,\cdots,f(A_k)\subset{}B_k$ for some subsets
$A_1,\cdots,A_k$ of $X$ and $B_1,\cdots,B_k$ of $Y$. Let
$\underline{f}:
(X,A_1\cup\cdots\cup{}A_k)\to(Y,B_1\cup\cdots\cup{}B_k)$ and
$f_i:(X,A_i)\to(Y,B_i)$, $i=1,\cdots{}k$, be defined by $f$. Then
for any sequence
$b_1\in{}H^{n_1}(Y,B_1;G_1),\cdots,b_k\in{}H^{n_k}(Y,B_k;G_k)$ we
have
\itemitem{(1)}
$b_1\ssm\cdots\ssm{}b_k\in{}
H^n(Y,B_1\cup\cdots\cup{}B_k;\bigotimes^k_{i=1}G_i)$, where
$n=n_1+\cdots+n_k$, and
\itemitem{(2)}
$\underline{f}^{\ast}(b_1\ssm\cdots\ssm{}b_k)=
f_1^{\ast}(b_1)\ssm\cdots\ssm{}f_k^{\ast}(b_k)$.

\proclaim{Lemma 3B.1} Let $X$ be a compactum with category
$(G_1,\cdots,G_k;n_1,\cdots,n_k)$. Then for every sequence
$a_1\in{}H^{n_1}(X;G_1),\cdots,a_k\in{}H^{n_k}(X;G_k)$ the product
$a_1\ssm\cdots\ssm{}a_k\in{}H^n(X;\bigotimes^k_{i=1}G_i)$ is equal
to $0$.
\endproclaim

\demo{Proof} From the hypothesis we infer that there exist a
closed covering $\{F_1,\cdots,F_k\}$ of $X$ such that each
homomorphism $H^{n_i}(F_i\hookrightarrow{}X;G_i)$ is trivial.
Consider the following portion of the cohomology exact sequence of
$(X,F_i)$:
$$ H^{n_i}(X,F_i;G_i)@>{j_i^{\ast}}>>H^{n_i}(X;G_i)@>{}>>H^{n_i}(F_i;G_i) \ $$
where $j_i:X\to(X,F_i)$ is defined by $\id_X$. Hence the right
homomorphism is trivial. So, as $a_i\in{}H^{n_i}(X;G_i)$ there is
an element $b_i\in H^{n_i}(X,F_i;G_i)$ such that
$j_i^{\ast}(b_i)=a_i$. Setting $Y=X$, $f=\id_X$, $A_i=\emptyset$
and $B_i=F_i$, we have
$\underline{\id}_X:X\to(X,F_1\cup\cdots\cup{}F_k)$ and
$(\id_X)_i=j_i$. It follows  from (2) that
$$
\underline{\id}^{\ast}_X(b_1\ssm\cdots\ssm{}b_k)=
j_1^{\ast}(b_1)\ssm\cdots\ssm{}j_k^{\ast}(b_k)=
a_1\ssm\cdots\ssm{}a_k \ .
$$
Since $X=F_1\cup\cdots\cup{}F_k$, we have
$H^n(X,F_1\cup\cdots\cup{}F_k;\bigotimes^k_{i=1}G_i)=0$. Thus
$b_1\ssm\cdots\ssm{}b_k=0$, hence $a_1\ssm\cdots\ssm{}a_k=0$,
which completes the proof. \qed
\enddemo

\proclaim{Theorem 3B.2} Let $X$ be an $n$-dimensional compactum
and let $Y_1,\cdots,Y_k$ be compacta such that $n=n_1+\cdots+n_k$,
where $n_i=\dim{}Y_i\ge1$. If there is a mapping
$f:X\to{}Y_1\times\cdots\times{}Y_k$ such that
$H^n(f\,;\,\bigotimes_{i=1}^k{}G_i)\neq 0$ then $X$ is not of
category $(G_1,\cdots,G_k; n_1,\cdots,n_k)$.
\endproclaim

\demo{Proof} Suppose, on the contrary, that
$X\in\cat(G_1,\cdots,G_k;n_1,\cdots,n_k)$. Then by Theorem 3A.1
there exist elements
$a_1\in{}H^{n_1}(X;G_1),\cdots,a_k\in{}H^{n_k}(X;G_k)$ such that
$a_1\ssm\cdots\ssm{}a_k\in{}H^n(X;\bigotimes_{i=1}^k{}G_i)$ is not
$0$. This contradicts Lemma 3B.1, and completes the proof. \qed
\enddemo

\proclaim{Corollary 3B.3} Let $X\subset Y_1\times\cdots\times Y_k$
be a compactum lying in the product of compacta, and let
$G_1,\cdots,G_k$ be Abelian groups. If
$H^n(X;\bigotimes_{i=1}^k{}G_i)\ne0$, where $n=n_1+\cdots+n_k$ and
$n_i=\dim{}Y_i\ge1$, then $X$ is not of category
$(G_1,\cdots,G_k;n_1,\cdots,n_k)$.
\endproclaim

\demo{Proof} This corollary follows from Theorem 3B.2 because
$H^n(X\hookrightarrow{}Y_1\times\cdots\times{}
Y_k;\bigotimes_{i=1}^k{}G_i)\ne0$ (cf. the proof of the second
assertion of Corollary 3A.4). \qed
\enddemo

\proclaim{Corollary 3B.4} Let $X$ be a compactum embeddable in the
product $Y_1\times \cdots \times Y_n$ of $n$ curves. If
$H^n(X)\ne0$ then $\cat X
> n$.\qed
\endproclaim\bs

 \head 4. Embedding compacta into symmetric products of
curves
\endhead

\bs In this chapter we establish a basic result on symmetric
products of curves, see Theorem 4E.1. Its principal assertion is
analogous to Corollary 2A.3 and reads: any compactum which admits
algebraically non-trivial mapping into the $n$th symmetric product
of a curve admits an algebraically non-trivial mapping into the
$n$-torus $\T^n$. In order to prove this result we need several
auxiliary observations which will be presented in a series of
sections. The most important result applied in the proof is the
result of B.W. Ong \cite{On} on the symmetric products of a
bouquet of circles, we recall it in section 4D. The remaining
assertions of Theorem 4E.1 follow then from Lemma 2B.3 and Theorem
3A.1. (The results of this chapter formulated for compacta hold
for compact Hausdorff spaces as well.) \bs

\centerline{4A. {\it Symmetric products as functors}}\ms

As usual, $S_n$, $n\ge1$, denotes the $n$th {\it symmetric group},
i.e. the permutation group of the set $\{1,\cdots,n\}$. There is a
standard (left) action of $S_n$ on $X^n$ (permuting coordinates of
points) for any non-void set $X$. Points of the diagonal $\Delta^n
(X)=\{(x,\cdots,x):x\in X\}$ are the only fixed points of this
action. Hence any orbit $S_n\cdot(x,\cdots,x)$ is a one-point set
$\{(x,\cdots,x)\}$. \ss

First, recall the definition of the symmetric product. Given a
space $X\neq\emptyset$ and $n\ge 1$ the $n$th {\it symmetric
product} of $X$, denoted $SP^n(X)$, is the space of all {\it
orbits} $S_n\cdot x$, where $x\in X^n$, of the standard
(continuous) action of $S_n$ on $X^n$. We identify $SP^1(X)$ with
$X$ under the assignment $\{x\}\to x$. Let $q:X^n \to SP^n(X)$
denote the quotient mapping taking each point to the orbit of that
point, i.e. $q(x)=S_n\cdot x$. In other words, for any $x,y\in
X^n$ we have\ss

    (1) $q(x)=q(y)$ {\it if and only if } $y=\alpha\cdot x$ {\it for
    some }$\alpha \in S_n$.

\ss \noi It follows that \ss

    (2) {\it q is open and closed}.

\ss Notice that $SP^n$ is a {\it functor} from the category of
topological spaces to itself. In fact, if $f:X\to Y$ then
$f^n:X^n\to Y^n$ is an equivariant mapping in the following
sense:\ms

    (3) $f^n(\alpha \cdot x)=\alpha \cdot f^n(x)$
    {\it for each} $x\in X^n$ and $\alpha \in S_n$. \ms

\noi Consequently, $f^n$ takes orbits in $X^n$ onto orbits in
$Y^n$. Therefore, there is a uniquely determined mapping
$SP^n(f):SP^n(X)\to SP^n(Y)$ making the following diagram
commutative:
$$
\CD X^n @ >{f^n}>> Y^n \\
 @ V{q}VV @V{q}VV  \\
SP^n(X) @ >{SP^n(f)}>> SP^n(Y) \ .
\endCD
$$

Notice that $SP^n$ sends injections to injections and surjections
to surjections. In particular, $SP^n(f):SP^n(X)\to SP^n(Y)$ is an
embedding if $X$ is compact and $f:X\to Y$ is an injection.

If $Y_1,\cdots,Y_n$ are pointed metrizable spaces then the
one-point union $Y=\bigvee_{i=1}^{n}Y_i$ is metrizable and each
$Y_i$ can be regarded as a closed subset of $Y$. Then
$Y_1\times\cdots\times Y_n\subset Y^n$ is closed. Hence, by (2),
$q:Y^n\to SP^n(Y)$ embeds $Y_1\times\cdots\times Y_n$ in
$SP^n(Y)$. If all $Y_i$ are $1$-dimensional then so is $Y$. This
combined with the Nagata Theorem 1.1 implies that {\it any
$n$-dimensional metrizable space can be embedded in the symmetric
product} $SP^{n+1}(Y)$, {\it where $Y$ is metric $1$-dimensional}.
For compacta we have more: {\it any $n$-dimensional compactum can
be embedded in the symmetric product $SP^{n+1}(\mu)$}, where $\mu$
stands for the Menger curve. This holds because $SP^{n}(\mu)$
contains a copy of the product $\mu^{n}$ (as $\mu$ is locally
homogeneous).

Moreover, $SP^n$ is a {\it homotopy functor}, i.e. if $f\simeq g$
then $SP^n(f)\simeq SP^n(g)$. Hence $SP^n$ preserves the homotopy
type of spaces. This property is important for our discussion.
Another important property of that functor is {\it continuity}: if
the projections $\{p_i:X \to X_i\}$ represent the limit of the
inverse sequence $\{X_1 @<{p_{1,2}}<< X_2 @<{p_{2,3}}<< \cdots \}$
then the projections $\{SP^n(p_i):SP^n(X) \to SP^n(X_i)\}$
represent the limit of the inverse sequence $\{SP^n(X_1)
@<{SP^n(p_{1,2})}<< SP^n(X_2) @<{SP^n(p_{2,3})}<< \cdots \}$. \ms

\centerline{4B. {\it Symmetric products of a pointed space}}\ms

For $1\leq k \leq l$ there is a canonical monomorphism
$s_{k,l}:S_k \to S_l$ defined as follows: $s_{k,l}(\beta)$, for
$\beta \in S_k$, is a permutation on $\{1,\cdots,l\}$ which acts
as $\beta$ on $\{1,\cdots,k\}$ and keeps each element $> k$ fixed.

Now consider a pointed space $X$ with a base point $\ast$. Then
each space $X^k$,$k\geq1$, is regarded as pointed with
$(\ast,\cdots,\ast)$ as the base point. Put $X^0=\{\ast\}$ and
regard it as pointed as well. Let $u_{k,l}:X^k\to X^l$ be
canonical embedding given by:
$(x_1,\cdots,x_k)\to{}(x_1,\cdots,x_k,\ast,\cdots,\ast)$. This
embedding preserves the base points and is equivariant in the
following sense:\ms

    (4) $u_{k,l}(\beta \cdot x)=s_{k,l}(\beta)\cdot u_{k,l}(x)$ {\it for each} $x\in
    X^k$.\ms

\noi Thus orbits in $X^k$ go under $u_{k,l}$ into orbits in $X^l$.
Moreover, we have the following property:\ms

    (5){\it If } $u_{k,l}(y)=\alpha\cdot u_{k,l}(x)$ {\it for} $x,y\in X^k$ {\it
    and}
    $\alpha\in S_l$, {\it then} $y=\beta\cdot x$ {\it for some}
    $\beta\in S_k$.\ms

    \noi It follows that \ms

    (6) $(S_l \cdot u_{k,l}(x))\cap u_{k,l}(X^k)= u(S_k \cdot x)$ {\it for each}
    $x\in X^k$.\ms

\noi Thus $u_{k,l}$ induces a mapping $v_{k,l}:SP^k(X)\to
SP^l(X)$, $S_k\cdot x\to S_l\cdot u_{k,l}(x)$, which makes the
following diagram commutative:
$$
\CD X^k @ >{u_{k,l}}>> X^l \\
 @ V{q}VV @V{q}VV  \\
SP^k(X) @ >{v_{k,l}}>> SP^l(X) \ .
\endCD
$$
One can define analogous commutative diagram for $k=0$ as well,
setting $SP^0(X)=\{\ast\}$, $u_{0,l}(\ast)=(\ast,\cdots,\ast)$ and
$v_{0,l}(\ast)=\{(\ast,\cdots,\ast)\}$. Notice that
$v_{k,k}=id_{SP^k(X)}$ and $v_{l,m}\circ v_{k,l}=v_{l,m}$ for
$k\leq l \leq m$. By (5), $v_{k,l}$ is injective, hence $v_{k,l}$
is an embedding if $X$ is compact - we call it {\it canonical
embedding}. Therefore, the image $SP^{k,l}(X)=v_{k,l}(SP^k(X))$ is
a copy of $SP^k(X)$ in $SP^l(X)$ in that case. In this way one
obtains a {\it canonical direct sequence of embeddings} of
symmetric products of $X$: $SP^0(X) @>{v_{0,1}}>> SP^1(X)
@>{v_{1,2}}>>SP^2(X)@>{v_{2,3}}>> \cdots $. All the mappings
belong to the category of pointed spaces.\ms

A pointed space $X$ with a base point \, $\ast$ \, is said to be a
{\it bouquet} of pointed subspaces $(X_j, \ast)$, for $j\in J$, if
each $X_j$ is a closed subset of $X$, $X_j\neq \{\ast\}$,
$X=\bigcup_{j\in J}X_j$, and $X_j\cap X_{j^\prime}=\{\ast\}$ for
$j\neq j^\prime$. The spaces $X_j$ are called {\it leaves}.

\ms\centerline{4C. {\it Symmetric products of a bouquet of k
circles and skeleta of torus $\T^k$ }}\ms

\noi If $K_1,\cdots, K_n$ are $CW$ complexes then we put
$$K_1\kw\cdots\kw{}K_n=\{\sigma_1\times\cdots\times\sigma_n:\;
\sigma_1\in{}K_1,\cdots,\sigma_n\in{}K_n \}. $$ Notice that
$K_1\kw\cdots\kw{}K_n$ is a $CW$ complex, it is called {\it
product cell complex} of the complexes. Given characteristic maps
of cells of $K_i$'s the cells of $K_1\kw\cdots\kw{}K_n$ have
natural characteristic maps obtained by taking products of the
corresponding characteristic maps of factors.\ss

The circle $\Sb^1$ will be regarded as a $CW$ complex with cell
structure $S=\{\sigma^0, \sigma^1\}$ composed of a $0$-cell
$\sigma^0=\{1\}$ and a 1-cell $\sigma^1=\Sb^1$. Consequently, for
any $k\geq1$, the $k$-torus $\T^k=\Sb^1\times\cdots\times \Sb^1$,
$k$ times, will be regarded as a $CW$ complex with the cell
structure $S \kw \cdots \kw S$, $k$ times. The cell structure will
be called {\it canonical} $CW$ structure on $\T^k$. The space $|(S
\kw \cdots \kw S)^{(n)}|$ of the $n$th skeleton of $S \kw \cdots
\kw S$, where $n\geq0$, carries the $CW$ complex structure $(S \kw
\cdots \kw S)^{(n)}$; with that structure it will be denoted by
$(\T^k)^{(n)}$ and called the $n$th {\it skeleton of the torus}
$\T^k$. For each $J\subset\{1,\cdots,k\}$ by $\sigma_J$ we denote
the cell $\sigma_1\times\cdots\times{}\sigma_k\in (S \kw \cdots
\kw S)$, where $\sigma_j=\sigma^1$ if $j\in{}J$, and
$\sigma_j=\sigma^0$ otherwise. It is a cell of dimension $|J|$,
$\dim\sigma_J=|J|$. Topologically, it is a torus of the same
dimension. The cell $\sigma_{\emptyset}$ is the unique 0-cell of
$S \kw \cdots \kw S$, hence lies in every cell. One easily sees
that \ms

\item\item{(1)} $\sigma_J\cap\sigma_{J^\prime}=\sigma_{J\cap J^\prime}.$

\ms\noi Each cell $\sigma_J$ is a retract $\T^k$ by a retraction
$\T^k \to \sigma_J$ which takes each point $(x_1,\cdots,x_k)$ to
$(y_1,\cdots,y_k)$, where $y_i=x_i$ if $i\in J$, and $y_i=1$
otherwise. This retraction will be called {\it canonical}.

Any two $m$-cells of the canonical $CW$ structure on $\T^k$ can be
joined by a chain of $m$-cells so that any two consecutive
$m$-cells meet in an $(m-1)$-cell. If $n\leq k$ then any cell of
the $n$th skeleton lies in an $n$-cell. Thus, we see that, using
terminology of 5A, the $n$th skeleton (for $k\geq 2$) is close to
be a simple ramified $n$-manifold complex (it is not such a
complex because it is not regular).

Now we recall the main result concerning the bouquets of circles
(see \cite{On, Theorem 3.1}).

\proclaim{Theorem 4C.1 (Ong)} Let $\bigvee_{i=1}^k S_i$, $k\geq1$,
be a bouquet of circles. Then, for any $n\geq1$, there is a
homotopy equivalence $\psi:SP^n(\bigvee_{i=1}^k S_i)\to
(\T^k)^{(n)}$.\qed
\endproclaim

\noi{\bf Note.} Moreover, one can construct $\psi$ to be a
surjective mapping with fibers being absolute retracts.\bs

\centerline{4D. {\it Algebraic properties of skeleta of a torus}}

\ms \noi Assume $n\leq k$. Then put $\J_n = \{J: J\subset
\{1,\cdots,k\}, |J|=n\}$. For each $J\in \J_n$ let
$$r_J:(\T^k)^{(n)} \to \sigma_J$$
denote the retraction induced by the canonical retraction $\T^{k}
\to \sigma_J$ (under our assumption
$\sigma_J\subset(\T^k)^{(n)}$).

\proclaim{Lemma 4D.1} For $n\leq k$ the homomorphism
$$h:\bigoplus_{J\in\J_n}H^n(\sigma_J)\to H^n((\T^k)^{(n)}),
\;\;\;\;\;\; (a_J)_{J\in\J_n}\to \sum_{J\in\J_n} H^n(r_J)(a_J),$$
is an isomorphism.
\endproclaim

\demo{Proof} We may assume that $n\geq1$. (For $n=0$ the
homomorphism $h$ is the identity isomorphism on $H^0(\sigma_0))$.
Let $\rho:(\T^k)^{(n)}\to(\T^k)^{(n)}/(\T^k)^{(n-1)}$ denote the
quotient mapping. For each $J\in\J_n$ the image $\rho(\sigma_J)$
is an $n$-sphere which we denote by $S_J$. Since
$(\T^k)^{(n)}=\bigcup_{J\in\J_n}\sigma_J$ and
$\sigma_J\cap{}\sigma_{J'}\subset(\T^k)^{(n-1)}$ for
$J,J'\in\J_n$, $J'\ne{}J$, we have
$(\T^k)^{(n)}/(\T^k)^{(n-1)}=\bigvee_{J\in\J_n}S_J$ (a bouquet of
$n$-spheres).

For each $J\in\J_n$, let $r'_J:\bigvee_{J\in\J_n}S_J\to S_J$ be
the retraction which takes $S_{J'}$, $J'\ne J$, to the base point.
Then the following diagram
$$
\CD (\T^k)^{(n)} @>{r_J}>> \sigma_J \\
  @V{\rho}VV @V{\rho|\sigma_J}VV \\
 \bigvee_{J\in\J_n}S_J @ >{r'_J}>> S_J
\endCD
$$
commutes. Consequently, the diagram
$$
\CD  H^n((\T^k)^{(n)})@ <{h}<<
\bigoplus_{J\in\J_n}H^n(\sigma_J) \\
 @A{H^n(\rho)}AA
@A{\bigoplus{}H^n(\rho|\sigma_J)}AA \\
H^n(\bigvee_{J\in\J_n}S_J) @ <{h'}<< \bigoplus_{J\in\J_n}H^n(S_J))
\endCD
$$
commutes, where $h'$ is a homomorphism analogous to $h$.

One easily sees that $H^n(\rho)$ is an epimorphism. (Indeed, this
follows from examination of a portion of the homomorphism
$CS(\bigvee_{J\in\J_n}S_J, \ast)\to CS ((\T^k)^{(n)},
(\T^k)^{(n-1)})$ between the cohomology exact sequences induced by
the relative homeomorphism $\rho:((\T^k)^{(n)}, (\T^k)^{(n-1)})\to
(\bigvee_{J\in\J_n}S_J, \ast)$.) The homomorphism
$\bigoplus{}H^n(\rho|\sigma_J)$ is an isomorphism since each
$H^n(\rho|\sigma_J)$ is an isomorphism. It is well known that $h'$
is an isomorphism. Thus, by commutativity of the above diagram,
$h$ is an epimorphism.

For each $J\in\J_n$, let $i_J:\sigma_J\to(\T^k)^{(n)}$ denote the
inclusion. If $J,J'\in\J_n$ and $J'\ne J$, then
$(r_{J'}\circ{}i_J)(\sigma_J)=\sigma_{J\cap
J'}\subset(\sigma_J)^{(n-1)}$, hence $H^n(r_{J'}\circ{}i_J)=0$.
And $H^n(r_{J}\circ{}i_J)=id_{H^n(\sigma_J)}$, since
$r_{J}\circ{}i_J=id_{\sigma_J}$. It follows that the composition
$$
\bigoplus_{J\in\J_n}H^n(\sigma_j) @ <({H^n(i_J))}<<
H^n((\T^k)^{(n)})@ <{h}<< \bigoplus_{J\in\J_n}H^n(\sigma_J)
$$
is the identity homomorphism. Thus $h$ is a monomorphism.
Consequently, $h$ is an isomorphism. This ends the proof of the
lemma.
 \qed \enddemo\ms

 \centerline{4E. {\it The Main Theorem}}\ms

Now we have at disposal all necessary tools needed for proving the
following main theorem of this chapter.

\proclaim{Theorem 4E.1} Let $f: X \to SP^n(Y)$ be a map from a
compact space $X$ to the $n$th symmetric product of a curve $Y$
such that the induced homomorphism $f^{\ast}: H^n(SP^n(Y)) \to
H^n(X)$ is non-trivial. Then there is a map $g:SP^n(Y)\to \T^n$
such that $H^n(g\circ f):H^n(\T^n)\to H^n(X)$ is non-zero.
Consequently, $\rank H^1(X)\ge n$ and there exist elements
$a_1,\cdots,a_n\in H^1(X)$ such that $a_1\ssm\cdots\ssm a_n\neq0$.
\endproclaim

\noi{\bf Note.} An analogous theorem holds if we take any Abelian
group as the coefficient group for cohomology in place of $\Z$.

\ms\demo{Proof} We can present $Y$ as the limit $Y=\invlim
\{P_1\longleftarrow P_2\longleftarrow\cdots\}$, where $P_i$ are
graphs. Then $SP^n(Y)= \invlim \{SP^n(P_1)\longleftarrow
SP^n(P_2)\longleftarrow\cdots\}$. Since $f^{\ast}: H^n(SP^n(Y))
\to H^n(X)$ is non-trivial, by the continuity of {\v{Cech}
cohomology, there is an index $l\ge 1$ such that
$$(pr_l\circ f)^{\ast}: H^n(SP^n(P_l)) \to H^n(X)$$ is nontrivial,
where $pr_l:SP^n(Y)\to SP^n(P_l)$ is the projection.

Since functor $SP^n$ preserves homotopy type, we may assume that
$P_l$ is a bouquet of circles, $P_l=\bigvee_{i=1}^k S_i$, where
$k\ge 1$. By Theorem 4C.1 there is a homotopy equivalence
$\psi:SP^n(\bigvee_{i=1}^k S_i)\to (\T^k)^{(n)}$. Therefore, the
homomorphism
$$(\psi\circ pr_l\circ f)^{\ast}: H^n((\T^k)^{(n)}) \to H^n(X)$$
is nontrivial. It follows that $k\ge n$. By Theorem 4F.1 there is
$J\in\J_n$ such that the homomorphism
$$(r_J\circ \psi\circ pr_l\circ f)^{\ast}: H^n(\sigma_J) \to H^n(X)$$
is nontrivial. As $\sigma_J$ is homeomorphic to $\T^n$, this
concludes the proof of the first assertion of the theorem. To
obtain the remaining assertions we apply Lemma 2B.3 and Theorem
3A.1. \qed
\enddemo

\ms\proclaim{Corollary 4E.2} Let $X$ be a compactum embeddable in
the nth symmetric product of a curve. If $H^n(X)\neq 0$ then
$\rank H^1(X)\geq n$ and there exist elements $a_1,\cdots,a_n\in
H^1(X)$ such that $a_1\ssm\cdots\ssm a_n\neq0$.
\endproclaim\ms

\noi{\bf Note.} In his recent paper R. Cauty  \cite{C2} formulates
the following related result: {\it under the above assumptions, if
$H^n(X)\neq 0$ then $H^1(X)\neq 0$}. His proof resorts to more
sophisticated techniques.\qed

\ms\demo{Proof} Let $f:X\to SP^n(Y)$ be an embedding, where $Y$ is
a curve. Since ${\dim}\;SP^n(Y)=n$, combining the hypothesis
$H^n(X)\neq 0$ with the Hopf Extension Theorem, we infer that
there is an element $\alpha\in H^n(SP^n(Y))$ such that
$H^n(f)(\alpha)\neq 0$. Then the conclusion follows from the above
theorem.\qed
\enddemo\ms

Next corollary directly follows from the above corollary because
$H^n(M)\neq 0$ for any closed manifold $M$.

\ms\proclaim{Corollary 4E.3} If a closed n-manifold M, $n\geq 2$,
can be embedded in the nth symmetric product of a curve then
$\rank H^1(M)\geq n$. In particular, $\Sb^n$ is not embeddable in
any $n$th symmetric product of a curve. The same holds for any
$\Sb^{n_1}\times\cdots\times\Sb^{n_k}$, where $n_1+\cdots+n_k=n$
and $n_i\geq 2$.\qed
\endproclaim\ms

\noi{\bf Note.} The second assertion is an analog of the Borsuk
Theorem. It follows from \cite{C2} as well.

\head 5. Locally connected generalized manifolds in products of
curves
\endhead

This chapter splits into sections 5A-5E. In section 5A we
introduce broad classes of continua (each including in particular
all closed manifolds) and call them quasi manifolds, pseudo
manifolds, para manifolds and ramified manifolds, and establish
some basic results. In Theorem 5B.1 we prove that embeddings of
locally connected quasi $n$-manifolds in products of $n$ curves
can be factored through special embeddings in products of $n$
graphs with no endpoint. Then we construct an example of a closed
surface lying in a product of two curves whose image under either
projection is not a graph (Example 5B.2). Theorem 5B.1 has
noteworthy consequences. For example, it follows that no locally
connected and unicoherent quasi $n$-manifold can be embedded in
any product of $n$ curves. In section 5C we present a list of
basic properties of ramified pseudo $n$-manifolds lying in
products of $n$ graphs. To obtain these properties we carefully
study the "fibers" of the projections restricted to the
pseudo-manifold. In section 5D we prove a fundamental Theorem 5D.5
on algebraic structure of locally connected quasi $n$-manifolds
lying in products of $n$ curves. In particular, this implies that
there exist contractible 2-dimensional polyhedra not embeddable in
products of two curves. Thus we reveal acyclic polyhedra which
have the same property as $\Sb^2$ in the Borsuk theorem. In
section 5E we prove that any $2$-dimensional collapsible
polyhedron (in particular, any cone over a graph) can be embedded
in a product of two trees.

\ms \centerline{5A. {\it Definitions and general properties of
certain generalized manifolds}}\ms

Let $K$ be a $CW$ complex. Then the open cells of $K$ (that is,
the interiors $\kn\sigma$ of the cells $\sigma \in K$) form a {\it
partition} of $|K|$, i.e. they are mutually disjoint and cover
$|K|$. It follows from the definition of a $CW$ complex that for
each skeleton $K^{(n)}$ the space $|K^{(n)}|$ is the union of open
cells with dimension $\leq n$. If $\dim\sigma =n$ then
$\partial\sigma\subset |K^{(n-1)}|$ and
$\kn\sigma\cap|K^{(n-1)}|=\emptyset$. A cell $\sigma \in K$ is
said to be {\it proper} if it is a union of open cells of $K$. If
each cell of $K$ is proper then $K$ is said to {\it have proper
cells}.

\proclaim{Proposition 5A.1} A CW complex K has proper cells if and
only if for each two cells $\sigma, \tau \in K$ the condition
$\kn\sigma\cap\tau\neq\emptyset$ implies $\sigma\subset\tau$
$($that is, $\sigma$ is a face of $\tau$$)$. \qed\endproclaim

\proclaim{Corollary 5A.2} If $K_1,\cdots, K_n$ are CW complexes
with proper cells then $K_1\kw\cdots\kw{}K_n $ has proper cells as
well. \qed\endproclaim

Notice that there exist finite $CW$ complexes with proper cells
which are not regular. (The canonical $CW$ structure on the torus
$\T^k$ has this property.) But the converse is true:

\proclaim{Lemma 5A.3} Any regular CW complex K has proper cells.
\endproclaim

 \demo{Proof} We must show that any cell $\tau\in K$
is the union of some open cells. By induction we may assume that
this holds for all cells with dimension $<n+1=\dim \tau$, $n\geq
0$. Note that $\partial\tau$ is a subset of the union of the
$n$-cells $\sigma_1,\cdots,\sigma_k$ of $K$ such that
$\kn\sigma_i\cap\tau\neq\emptyset$. It remains to show that
$\sigma_i\subset\tau$ for each $i$. To this end fix $i$. As
$\partial\tau\cap\kn\sigma_i$ is a non-void closed subset of
$\kn\sigma_i$, to get the conclusion, it is enough to prove that
$\partial\tau\cap\kn\sigma_i$ is open in $\kn\sigma_i$. Since
$\partial\tau$ is an $n$-sphere in $|K^{(n)}|$ and $\kn\sigma_i$
is open in $|K^{(n)}|$, for each point
$x\in\partial\tau\cap\kn\sigma_i$ there is an open $n$-cell
containing $x$ and wholly lying in this intersection. By the
Brouwer Domain Invariance Theorem such a ball is a neighborhood of
$x$ in $\kn\sigma_i$. Thus $\partial\tau\cap\kn\sigma_i$ is open
in $\kn\sigma_i$, which completes the proof.\qed
\enddemo\ms

  From the classic Borsuk Separation Theorem relating closed sets separating
${\Sb}^n$ to their essential mappings into ${\Sb}^{n-1}$ (cf.
\cite{E-S, p. 302}) we infer the following fact.

\proclaim{Lemma 5A.4} For any $n$-manifold $M$ and a point $x_0 \in
M$ there is an open neighborhood $V$ of $x_0$ in $M$ such that every
closed subset $F$ of $M$ separating $M$ between $x_0$ and $M
\setminus V$ admits an essential map into $\Sb^{n-1}$. $($In fact,
this holds for every neighborhood $V$ which is an open n-disc.$)$
\qed
\endproclaim

The revealed property can be used to define, for each natural $n\ge
1$, a new class of $n$-dimensional continua (comprising in
particular all connected $n$-manifolds) playing an important role in
our investigations. Namely, an $n$-dimensional continuum $X$ is said
to be a {\it quasi} $n$-{\it manifold at a point} $x \in X$ if there
is an open neighborhood $V$ of $x$ in $X$ such that every closed
subset $F$ of $X$ with $\dim F \leq {n-1}$ separating $X$ between
$x$ and $X \setminus V$ admits an essential map into $\Sb^{n-1}$.
(Recall that a closed set $F\subset X$ is said to {\it separate X
between subsets} $A$ and $B$ if there exist disjoint open sets $U$
and $V$ such that $X\setminus F=U\cup V$, $A\subset U$ and $B\subset
V$.) Notice that $X$ is a quasi 1-manifold at $x\in X$ if and only
if $x$ is not an endpoint of $X$. (A point of a space is said to be
its {\it endpoint} if that point admits arbitrarily small open
neighborhoods whose boundaries are one-point sets.) If $X$ is a
quasi $n$-manifold at $x$ then every neighborhood of $x$ is
$n$-dimensional. If $V$ is as in the first definition then any other
open neighborhood $W$ of $x$ contained in $V$ has the separation
property as well. Also note that if a closed set $F\subset V$
separates $\overline{V}$ between $x$ and $\partial V$ then $F$
separates $X$ between $x$ and $X\setminus V$. If $X$ is a quasi
$n$-manifold at every point of $X$ then it is called a {\it quasi}
$n$-{\it manifold}. Notice that an $n$-dimensional continuum which
is a union of quasi $n$-manifolds is a quasi $n$-manifold as well.

An $n$-dimensional continuum $X$ is said to be a {\it para
n-manifold} if each point of $X$ belongs to an open $n$-disc lying
in $X$ (not necessarily open in $X$). In other words, $X$ is a union
of open $n$-discs.

If $X$ is an $n$-dimensional continuum and a point $x\in X$ is an
element of an open $n$-disc lying in $X$ then $X$ is a quasi
$n$-manifold at $x$. (This follows from the observation that if
$E$ is an open disc in $X$ then $E$ is open in $\overline{E}$.)
Hence any para $n$-manifold is a quasi $n$-manifold. This simple
criterion can be used to detect many interesting quasi
$n$-manifolds which are not $n$-manifolds. For instance, both the
$\lq \lq$Bing house" and the $\lq\lq$dunce hat" are para
2-manifolds, so they are also quasi 2-manifolds. Those examples
are widely known primarily for being 2-dimensional contractible
and not collapsible polyhedra.\ms

And it is convenient to introduce other generalizations of
$n$-manifolds. First, for an $n$-dimensional continuum $X$ define
the following subsets: \ms $P(X)=\{x \in X: \text {\it x is an
element of an open n-cell lying in X and open in X} \}$; \ms $R(X)=
\{x \in X: \text{\it x is an element of an open n-cell lying in X}
\}$.\ms

\noi Thus $P(X)=X$ if and only if $X$ is an $n$-manifold, and
$R(X)=X$ if and only if $X$ is a para $n$-manifold.

Now we define the generalizations. An $n$-dimensional continuum $X$
is said to be a {\it pseudo $n$-manifold} ({\it ramified
$n$-manifold}, respectively) if $P(X)$ ($R(X)$, respectively) is
dense in $X$ and $\dim [X\setminus P(X)]\leq n-2$ ($\dim [X\setminus
R(X)]\leq n-2$, respectively). If, in addition, $P(X)$ ($R(X)$,
respectively) is connected then $X$ is said to be a {\it simple
pseudo n-manifold} ({\it simple ramified n-manifold}, respectively).

\ms The set of $n$-cells of a regular $CW$ complex $K$, $n\ge1$, is
said to be {\it chain connected} if every two $n$-cells of $K$ can
be joined by a finite sequence of $n$-cells of $K$ such that every
two successive cells meet along an $(n-1)$-cell. If each
$(n-1)$-cell of $K$ is a face of an $n$-cell, the above property
holds if and only if the space $|K^{(n)}|\setminus |K^{(n-2)}|$ is
connected.

Let $K$ be an $n$-dimensional finite regular $CW$ complex. Then
$K$ is said to be a {\it para n-manifold} (resp.,{\it quasi
$n$-manifold, pseudo $n$-manifold, ramified $n$-manifold}) {\it
complex} if the polyhedron $|K|$ is a para $n$-manifold (resp.,
quasi $n$-manifold, pseudo $n$-manifold, ramified $n$-manifold).
To name the complexes from the classes defined above for which the
set of $n$-cells is chain connected, we add the adjective {\it
simple} to the basic names, e.g. a {\it simple para $n$-manifold
complex}, etc.

Notice that any ramified $n$-manifold complex $L$ lying in a
simple pseudo $n$-manifold complex $K$, coincides with $K$. In
particular, this holds if $|K|$ is an $n$-manifold (as $K$ is a
simple pseudo $n$-manifold complex in this case).

(In the literature, the term {\it pseudo $n$-manifold} is often
used in a more restrictive sense to mean the space $|K|$, where
$K$ is a simple pseudo $n$-manifold complex.)

\proclaim{5A.5 Proposition} Let K be an n-dimensional finite regular
CW complex. Then\ss

{\rm(i)} $K$ is a pseudo $n$-manifold $($simple pseudo $n$-manifold,
{\it respectively}$)$ complex if and only if each cell of K is a
face of an n-cell, each $(n-1)$-cell of K is incident with exactly
two n-cells $($and the set of $n$-cells of $K$ is chain connected,
respectively$)$;\ss

{\rm(ii)} $K$ is a ramified $n$-manifold $($simple ramified
$n$-manifold, respectively$)$ complex if and only if each cell of K
is a face of an n-cell, each $(n-1)$-cell of K is incident with at
least two n-cells $($and the set of $n$-cells of $K$ is chain
connected, respectively$)$.\qed\endproclaim

\ms Our main observation in this subsection is Theorem~5A.8 below
which describes a basic property of quasi $n$-manifolds and ramified
$n$-manifolds lying in $n$-dimensional polyhedra. In the proof Lemma
5A.7 is needed. In the proof of Lemma 5A.7 we need in turn
Lemma~5A.6 below; it is relatively simple but not trivial. Lemma
5A.6 is certainly known to many topologists, and can be proved using
various arguments. We supply possibly the shortest one.

\proclaim {Lemma 5A.6} No proper closed subset of $\Sb^n$ admits
an essential map into $\Sb^n$.\endproclaim

\demo{Proof} We may assume that $n\ge1$. Consider a proper closed
subset $F$ of $\Sb^n$. Since $\Sb^n$ is a compact subset of
$\R^{n+1}$, the set $F$ can be regarded as a compact subset of
$\R^{n+1}$. Then note that $\R^{n+1}\setminus F$ is connected. So,
by the Borsuk Separation Theorem \cite{E-S, p. 302}, $F$ admits no
essential map into $\Sb^n$.\qed
\enddemo

\proclaim{Lemma 5A.7}{\rm(a)} Let $X$ be a quasi $n$-manifold at a
point $x$. If $U$ is a neighborhood of $x$ in $X$ and $h:U\to
\R^n$ is an embedding, then $h(U)$ is a neighborhood of $h(x)$ in
$\R^n$.

{\rm(b)} Let $X$ be a ramified n-manifold and let U be a non-void
open subset of X. If $h:U\to \R^n$ is an embedding such that
$h(U)$ is closed, then $h(U)=\R^n$.
\endproclaim

\demo{Proof of $(${\rm a}$)$} Suppose, to the contrary, that
$h(x)\in\partial h(U)$. Let $V$ be an open neighborhood of $x$ in
$X$ with $\overline{V}\subset U$ such that any closed subset $F$ of
$X$ separating $X$ between $x$ and $X\setminus V$ admits an
essential map into $\Sb^{n-1}$. Since $h(U\setminus V)$ is a closed
subset of $h(U)$ not containing $h(x)$ there is an open ball
$B(h(x),\varepsilon_0)$ in $\R^n$ such that
$B(h(x),\varepsilon_0)\cap h(U\setminus V)=\emptyset$. Also, there
is a sphere $S=\partial B(h(x),\varepsilon)$, $0 < \varepsilon
<\varepsilon_0$, such that $S\nsubseteq h(U)$. Since $S$ separates
$\R^n$ between $x$ and $h(U\setminus V)$ , the intersection $E=S\cap
h(U)$ separates $h(U)$ between $h(x)$ and $h(U\setminus V)$. Note
that $E$ is a proper subset of $S$. As $E=S\cap h(V)=S\cap
h(\overline V)$, the set $E$ is a compact subset of $h(V)$ which
separates $h(\overline V)$ between $h(x)$ and $h(\partial V)$. It
follows that $F=h^{-1}(E)$ is a compact subset of $V$ which
separates $\overline V$ between $x$ and $\partial V$. Hence
$F(\subset V)$ is a closed subset of $X$ which separates $X$ between
$x$ and $X\setminus V$. By our choice of $V$, there is an essential
map $F\to \Sb^{n-1}$. But $F(\approx E)$ is homeomorphic to a proper
closed subset of $\Sb^{n-1}$, which contradicts Lemma 5A.6.\enddemo

\demo{Proof of $(${\rm b}$)$} It is enough to show that $h(U)$ is
dense in $\R^n$. Let $V$ = int $h(U)$. First we show that $V$ is
dense in $h(U)$. In fact, since $R(X)$ is dense $X$ and $U$ is open,
the image $h(R(X)\cap U)$ is dense in $h(U)$. On the other hand,
$h(R(X)\cap U)$ is open in $\R^n$, by the Brouwer Invariance of
Domain Theorem. Hence $h(R(X)\cap U)\subset V$, so $V$ is dense in
$h(U)$. Therefore, $\partial V=h(U)\setminus V$. Consequently,
$\partial V \subset h(U\setminus R(X))$, hence $\dim
\partial V \leq n-2$, by the definition of a ramified
$n$-manifold. Therefore, $\R^n \setminus \partial V$ is connected
(cf. \cite {E, Theorem 1.8.13, p. 77}). It follows that $V$ is
dense in $\R^n$, for otherwise $\partial V$ separates $\R^n$.
Hence $h(U)=\R^n$, which completes the proof. \qed\enddemo

An easy modification of the above argument proves the following

\proclaim{Corollary 5A.8} If $f:X\to Y$ is an embedding of a
ramified n-manifold into a simple pseudo n-manifold, then
$h(X)=Y$.\qed\endproclaim

\proclaim{Theorem 5A.9} Let $X$ be either a quasi $n$-manifold or a
ramified $n$-manifold. If $f:X\to |K|$ is an embedding, where $|K|$
is an $n$-dimensional polyhedron, then $h(X)=|L|$, where $L$ is a
ramified $n$-manifold subcomplex of $K$. Moreover, if $X$ is a
pseudo $n$-manifold, then $L$ is a pseudo n-manifold complex.
\endproclaim

\demo{Proof} First we prove that $f(X)$ is a union of $n$-cells of
$K$. To this end, it is enough to show that each point $y\in f(X)$
is an element of an $n$-cell of $K$ which lies in $f(X)$. Notice
that there is an open neighborhood $V$ of $y$ in $f(X)$ such that,
for each cell $\sigma \in K$, the condition $V\cap
\sigma\neq\emptyset$ implies $y \in \sigma$. Since ${\dim}V=n$,
there is an $n$-cell $\sigma_0\in K$ such that
$V\cap\kn\sigma_0\neq\emptyset$. Hence $y\in \sigma_0$. It remains
to show that $\sigma_0\subset f(X)$. As $f(X)$ is closed, it is
enough to show that $\kn\sigma_0\subset f(X)$. Notice that $f(X)
\cap \kn \sigma_0$ is a non-void closed subset of $\kn\sigma_0$.
On the other hand, this set is open in $f(X)$. It follows that
$U=f^{-1}(\kn\sigma_0)$ is non-void and open in $X$, and
$f(U)=f(X)\cap\kn\sigma_0$. Since $\kn\sigma_0 \approx \R^n$, by
Lemma 5A.7 we infer that $f(X)\cap\kn\sigma_0$ is also open in
$\kn\sigma_0$, in both cases under discussion. Consequently, $f(X)
\cap \kn\sigma_0=\kn\sigma_0$. Hence $\kn\sigma_0 \subset f(X)$,
as desired. The remaining assertions simply follow from the Borsuk
Separation Theorem, cf. \cite {E-S, p. 302}.\qed\enddemo

By the above theorem, for any polyhedron $|K|$ which is either a
quasi $n$-manifold or a ramified $n$-manifold, $K$ is a ramified
$n$-manifold complex. Moreover, the converse holds for $n=1,2$. In
general, the converse fails: the suspension of the "dunce hat" (or
the "Bing house") can be presented as a polyhedron $|K|$ such that
$K$ is a ramified $3$-manifold complex but $|K|$ is no quasi
$3$-manifold.

Let $K$ be a ramified $n$-manifold complex. By its {\it
combinatorial component} we mean any maximal simple ramified
$n$-manifold subcomplex of $K$. One easily sees that $K$ is the
union of combinatorial components, and any two different
combinatorial components meet in a subcomplex of dimension
$\leq{n-2}$.

We have the following diagram of inclusions in the class of
connected polyhedra:
$$
\CD \{para \ n{-}manifolds\}\ & \supset & \ \ \{special\
para \ n{-}manifolds\} \\
    \cap & {} & \cap\\
\{quasi \ n{-}manifolds\} \ & \supset & \ \ \{special\
quasi \ n{-}manifolds\} \\
    \cap & {} & \cap\\
\{ramified \ n{-}manifolds\}\ & \supset & \ \ \{special\
ramified \ n{-}manifolds\} \\
\cup & {} & \cup \\
\{pseudo \ n{-}manifolds\} & \supset & \ \ \{special\ pseudo\
n{-}manifolds\}.
\endCD
$$

\ms\centerline{5B. {\it On locally connected quasi manifolds.}}
\centerline{\it From embeddings into products of curves to
embeddings into products of graphs}\ms

 Here we prove a useful theorem on factorization of embeddings of quasi
manifolds into products of curves through embeddings into product
of graphs.

\proclaim{Theorem 5B.1} Let $X$ be a locally connected quasi
$n$-manifold such that $H^1(X)$ has finite rank. If $f = (f_1,
\cdots ,f_n):X \to Y_1 \times \cdots \times Y_n$ is an embedding
of $X$ in the product of $n$ curves, then there exist mappings $g
= (g_1, \cdots ,g_n):X \to P_1 \times \cdots \times P_n$ and $h=
h_1 \times \cdots \times h_n:P_1 \times \cdots \times P_n \to Y_1
\times \cdots \times Y_n$ such that $f_i = h_i \circ g_i$ for each
$i=1,\cdots ,n$ $($hence $f=h \circ g$$)$, where $g_i:X \to P_i$
is a monotone surjection, $P_i$ is a graph with no endpoint
$($that is, $P_i$ is a quasi $1$-manifold$)$, and $h_i:P_i\to Y_i$
is $0$-dimensional.

In particular, if $X$ is embeddable in a product of $n$ curves,
then there exists an embedding $(g_1, \cdots ,g_n):X \to P_1
\times \cdots \times P_n$, where each $g_i:X \to P_i$ is a
monotone surjection, $P_i$ is a graph with no endpoint, and $
rank$ $H^1(P_i)\le rank$ $H^1(X)$.
\endproclaim

\ms\noindent{\bf Note.} It follows that {\it if} $f_i:X \to Y_i$
{\it is monotone, then} $f_i(X)$ {\it is a graph.} In fact, in
this case $f_i(X)=h_i(P_i)$ and $h_i:P_i\to Y_i$ is an embedding.
If $f_i$ is not monotone then $f_i(X)$ need not be a graph, see
Example 5B.3. The proof given below shows that {\it if $X$ is a
non-degenerate connected polyhedron $($or any non-degenerate
locally connected continuum whose $H^1(X)$ has finite rank$)$ then
$f_i(X)$ is a local dendrite.} \qed\bs

\demo{Proof} By the Whyburn factorization theorem, there is a
factorization $f_i=h_i \circ g_i$ ,
$$
X @>{g_i}>> P_i @>{h_i}>> Y_i,
$$
where $g_i$ is a monotone surjection, and $h_i$ is
$0$-dimensional. Since $Y_i$ is 1-dimensional and $h_i$ is
0-dimensional we infer that $\dim P_i \leq1$ (by a theorem of
Hurewicz). Clearly, $g = (g_1, \cdots ,g_n):X \to P_1 \times
\cdots \times P_n$ is an embedding. Since $\dim X=n$, it follows
that $\dim P_i>0$ for each $i$. Therefore, $P_i$ is a locally
connected curve, as $g_i$ is a surjection. Since $g_i$ is a
monotone surjection and $H^1(X)$ has finite rank, $P_i$ is
actually a local dendrite (see \cite{Kr, Lemma 3.1}). Hence each
point of $P_i$ has a closed neighborhood which is a dendrite.
First we show that\ms

 (*) {\it $P_i$ has no endpoint}.\ms

\noi For suppose $P_i$ has an endpoint $z_0$. Since $g_i(X) =
P_i$, there is a point $x_0 \in X$ such that $g_i(x_0) = z_0$.
Since $X$ is a quasi $n$-manifold at $x_0$ there is an open
neighborhood $V$ of $x_0$ in $X$ such that

$(1)$ {\it any closed $(n-1)$-dimensional subset of} $X$ {\it
separating} $X$ {\it between} $x_0$ {\it and} $X \setminus V$ {\it
admits an essential map to} $\Sb^{n-1}$.

Now we shall show that there is an open neighborhood $U$ of
$g(x_0)$ in $P_1 \times \cdots\times P_n$ such that
\roster
\item "(2)" $\overline{U} \cap g(X \setminus V)=\emptyset$,
\item "(3)" $\dim \partial U = n-1$,
\item "(4)" $\partial U$ {\it is contractible}.
\endroster

\noi To construct $U$ we assume, without loss of generality, that
$i=1$. Then $g(x_0)=(y_1,y_2,\cdots,y_n)$, where $y_1=z_0$. Note
that $g(V)$ is a neighborhood of $g(x_0)$ in $g(X)$, hence any
small enough $U$ satisfies (2). Since $z_0$ is an endpoint of
$P_1$, and each $P_j$ is a local dendrite, there exist sets
$U_1,\cdots,U_n$, as small as we please, such that each $U_j$ is
an open and connected neighborhood of $y_j$ in $P_j$ with
$\partial U_j$ finite, each $\overline{U_j}$ is a dendrite, and
$\partial U_1$ is a one-point set. Then the set
$U=U_1\times\cdots\times U_n$ has the desired properties. In fact,
as $U_j$'s are small, $U$ satisfies (2). Then note that
$$\partial U=\bigcup_{j=1}^{n} (
\overline{U_1}\times\cdots \times \overline{U_{j-1}}
\times
\partial U_j\times \overline{U_{j+1}}\times\cdots\times
\overline{U_n}).$$ \noi Hence (3) follows. To prove (4), note that
$(\partial U_1)\times\overline{U_2}\times\cdots\times
\overline{U_n}$ is a strong deformation retract of $\partial U$
(because $\partial U_1$, as a one-point set, is a strong
deformation retract of $\overline{U_1}$). Hence (4) follows from
the fact that $(\partial
U_1)\times\overline{U_2}\times\cdots\times \overline{U_n}$ is
contractible.

Now consider the set $F =\partial_{g(X)} (U \cap g(X))$. Observe
that it is a closed subset of $g(X)$ such that \roster
\item "(5)" $F \subset \partial U$,
\item "(6)" $F$ {\it separates} $g(X)$ {\it between} $g(x_0)$ {\it and} $g(X \setminus V)$.
\endroster
It follows that \roster
\item "(7)" $g^{-1}(F)$ {\it is closed}, $(n-1)$-{\it dimensional, and separates X
between} $x_0$ {\it and} $X\setminus V$.\endroster

Now we are ready to complete the proof of (*). Note that by (1)
and (7) there is an essential map $\varphi: F \to \Sb^{n-1}$. By
(3) and (5) there is a continuous extension $\varphi^*:
\partial U \to \Sb^{n-1}$ of $\varphi$. However, by (4),
$\varphi^*$ is null-homotopic, hence so is $\varphi$, a
contradiction. This proves (*).

Next we show that\ms

(**) {\it $P_i$ is a graph}.

\ms To prove (**) recall that $P_i$ is a local dendrite. Since
$P_i$ has no endpoint, it contains a circle. (Otherwise it is a
dendrite, hence contains an endpoint.) It follows that the union
of all simple closed curves in $P_i$ is a (not necessarily
connected) graph. Enlarging this set by the union of a finite
collection of arcs (e.g., adding arcs in $P_i$ irreducibly
connecting different components of the union), we get a connected
graph $Q_i (\subset P_i)$ such that for each component $C$ of $P_i
\setminus Q_i$ we have

$(8)$ $\overline{C}$ {\it is a dendrite and} $\partial C$ {\it
consists of a single point}.

\noi To prove (**) it suffices to show that $Q_i = P_i$.

Suppose, on the contrary, that $P_i \setminus Q_i \neq \emptyset$.
Then consider a component $C$ of $P_i \setminus Q_i$. It is an
open set in $P_i$. By (1), $\overline{C}$ is a dendrite and
$\partial C$ is a one-point set. There is a point $z_0 \in C$
which is an endpoint of the dendrite $\overline{C}$. Since $C$ is
open, $z_0$ is an endpoint of $P_i$ as well. This contradicts (*)
and ends the proof of (**).

As $g_i$ is a monotone surjection, the induced homomorphism
$H^1(g_i):H^1(P_i)\to H^1(X)$ is a monomorphism by the
Vietoris-Begle Theorem for $n=1$ (see e.g. \cite{Sp, 6.9, Theorem
15}). Therefore, rank $H^1(P_i) \le$ rank $H^1(X)$. This completes
the proof. \qed
\enddemo

It is well known that for any closed $n$-manifold $M$ the group
$H^1(M)$ has finite rank. Consequently, Theorem 5B.1 implies the
following

\proclaim{Corollary 5B.2} If a closed $n$-manifold is embeddable
in a product of $n$ curves, then it is also embeddable in a
product of $n$ graphs.\qed
\endproclaim

\noi {\bf{Note.}} This corollary answers a question posed (for
surfaces) by R. Cauty. A harder variant of this question (for
$n$-dimensional polyhedra) is still open \cite{C1}.\qed\ms

\proclaim{Example 5B.3} There exist a curve $X$ which is not a
graph, and a closed orientable surface $M$ in the product $X
\times X$ such that both projections $pr_i:X \times X \to X$ map
$M$ onto $X$. Moreover, $M$ is invariant under the canonical
involution on $X \times X$ which interchanges the coordinates.
\endproclaim

\demo{Proof} First we construct a closed orientable surface $N$ in
the product $Y_1 \times Y_2$ of two curves such that \roster
\item "(1)" $Y_1$ {\it is not a graph},
\item "(2)" {\it the projections} $q_i:Y_1\times Y_2 \to Y_i$
{\it map} $N$ {\it onto} $Y_i$.
\endroster
Define $Y_1$ to be the union $Y_1=\alpha_0 \cup \alpha_1
\cup\beta_1 \cup \beta_2$ of four arcs with common endpoints $a,b$
such that $\alpha_0 \cup \alpha_1 \cup \beta_1$ ia a
$\theta$-curve, $\beta_2\cap(\kn\alpha_0 \cup
\kn\alpha_1)=\emptyset$, and $\beta_1\cap \kn\beta_2$ is a compact
set with infinitely many components. Then $Y_1$ satisfies (1).
Define $Y_2$ to be a graph given by the formula: $Y_2=T_0 \cup
T_1$, where $T_0$, $T_1$ are two oriented circles whose
intersection $T_0 \cap T_1=L_0 \cup L_1$, where $L_0,L_1$ are
disjoint oriented arcs coherently oriented with each $T_i$. In
such a case, each $T_i$ can be presented as the union of four arcs
with disjoint interiors, $T_i=A_i \cup B_i \cup L_0 \cup L_1$,
such that $S_1=A_0 \cup A_1$ and $S_2=B_0 \cup B_1$ are disjoint
circles. We define the surface $N$ in $Y_1 \times Y_2$ by the
formula
$$
N=\alpha_0 \times T_0 \cup \alpha_1 \times T_1 \cup \beta_1\times
S_1 \cup \beta_2\times S_2.
$$

\noindent One easily verifies that $N$ is an orientable surface
satisfying $(2)$.

To construct the promised example we proceed as follows. Choose a
homeomorphism $h:\alpha_0 \to A_0$. Then define $X$ to be the
quotient space $X=(Y_1\sqcup Y_2)/x\sim h(x)$ for each $x\in
\alpha _0$). By $(1)$ we infer that $X$ is a curve but not a
graph. Let $X_i = h_i(Y_i)$, where $h_i:Y_i \to X$ are canonical
embeddings. Clearly, $X = X_1 \cup X_2$ and $X_1 \cap X_2 = A$,
where $A = h_1(\alpha_0) = h_2(A_0)$ is an arc. Let $t:Y_1 \times
Y_2 \to Y_2 \times Y_1$ denote the map given by $t(y,z) = (z,y)$.
Then we define $M (\subset X \times X)$ as follows:
$$
M = [(h_1 \times h_2)(N) \cup (h_2 \times h_1)(t(N))] \setminus
(\kn A \times \kn A).
$$
One easily verifies that $M$ is invariant under canonical
involution on $X \times X$. Since $\alpha_0 \times A_0 \subset N
\subset Y_1 \times Y_2$, $A_0 \times \alpha_0 \subset t(N) \subset
Y_2 \times Y_1$. Hence $A \times A = (h_1 \times h_2)(\alpha_0
\times A_0) \subset (h_1 \times h_2)(N) \subset (h_1 \times
h_2)(Y_1 \times Y_2) = X_1 \times X_2$. Likewise, $A \times A
\subset (h_2 \times h_1)(t(N)) \subset X_2 \times X_1$. Moreover,
$(X_1 \times X_2) \cap (X_2 \times  X_1) = A \times A$. Hence we
have
$$
(h_1 \times h_2)(N) \cap (h_2 \times h_1)(t(N)) = A \times A.
$$
Thus $M$ is the connected sum of orientable surfaces $(h_1 \times
h_2)(N)$ and $(h_2 \times h_1)(t(N))$, hence it is an orientable
surface. Applying $(2)$, we easily see that both projections
$pr_i:X \times X \to X$ map $M$ onto $X$. \qed
\enddemo

\ms In connection with the above proof let us notice the following
fact.

\ms \noi{\bf Note.} Let $M$ be any compactum lying in the product
$P_1\times P_2$ of two graphs. If $A$ is an arc in $P_1$ with
$p_1^{-1}(A)=A\times(S_1\cup\dots\cup S_k)$, $k\ge2$, where
$S_1,\dots,S_k$ are disjoint circles in $P_2$, then $M$ can be
embedded in the product $P_1'\times P_2$ in such a way that $P_1'$
is a curve and the image of $M$ under the projection $P_1'\times
P_2\to P_1'$ is not a graph. (In fact, $P_1'$ can be obtained from
$P_1$ by adding an arc $B$ with the same endpoints as that of $A$
in such a way that $A\cup B$ is not a graph.) \qed\ms

The surface $M$ constructed in Example 5B.3 meets the diagonal of
$X\times X$. Below we present another example of a surface in the
product $P\times P$, where $P$ is a graph, which is disjoint from
the diagonal of the product and is invariant under the canonical
involution on $P\times P$.

\proclaim{Example 5B.4} There exist a graph $P$ and a closed
orientable surface $M$ in $P \times P$ such that: $M$ is disjoint
with the diagonal of $P\times P$, both projections $pr_i:P \times
P \to P$ map $M$ onto $P$, and $M$ is invariant under the
canonical involution on $P \times P$.
\endproclaim

\demo{Proof} Fix any number $n\ge 4$. The graph $P$ is defined to
be a subset of $\Sb^1\times I$ given by
$$P=(\Sb^1\times \{0,1\}) \cup \{z_0,\cdots,z_{n-1}\}\times I,$$
where $z_j=\exp(2\pi i \frac{j}{n})$, $j=0, \cdots, n-1$. Then
define arcs $A_j\times\{0\}$, $A_j\times\{1\}$, $I_j$ and circles
$S_j$ in $P$ as follows:
$$A_j=\{\exp(2\pi it):t\in [\frac{j}{n},\frac{j+1}{n}]\}, \;\;
I_j=\{z_j\}\times I,\;\; S_j=I_j \cup (A_j \times \{0,1\}) \cup
I_{j+1}.$$ \noi (All indices in this construction are reduced
modulo $n$.) Finally, define tori $T_j$ to be the subsets of
$P\times P$ given by
$$T_j=S_j\times S_{j+2}.$$ \noi Notice that the intersection
$$D_j=T_i \cap T_{j+1}=I_{j+1}\times I_{j+3}$$
is a disc. Now we are ready to define the surface $M$, put
$$M=(T_0 \cup \cdots \cup T_{n-1})\setminus (\kn D_0\cup\cdots\cup\kn D_{n-1}).$$
One easily verifies that $M$ has all the desired properties.\qed
\enddemo\ms

\centerline{5C. {\it Ramified manifolds in products of graphs}}\ms

Here we establish some properties of ramified $n$-manifolds lying
in products of $n$ graphs. These properties will find essential
applications in the next section and in Chapter 6.\ms

\centerline{\it{Throughout this section we consider fixed graphs
$P_1=|K_1|$, $\dots$,$P_n =|K_n|$, $n\ge2$,}} \centerline
{\it{where each $K_i$ is a regular $1$-dimensional $CW$
complex.}}\ms \noi (So, each 1-cell of $K_i$, i.e. an {\it edge},
is an arc; its endpoints are called {\it vertices}.) Let us recall
that by $K_1\kw\cdots\kw{}K_n$ we denote the cell structure on
$P_1\times\cdots\times{}P_n$ defined by
$$K_1\kw\cdots\kw{}K_n=\{\sigma_1\times\cdots\times\sigma_n:
\sigma\in{}K_1,\cdots,\sigma_n\in{}K_n \}. $$
\centerline{\it{Also, we consider a fixed ramified $n$-manifold
$M=|K(M)|$ lying in}}
\centerline{\it{$P_1\times\cdots\times{}P_n$, where $K(M)$ is a
subcomplex of $K_1\kw\cdots\kw{}K_n$.}}\ms \noi By Theorem 5A.9 we
have

\proclaim{Property (a)} $K(M)$ is a ramified $n$-manifold complex;
in particular, $M$ is the union of $n$-cells of $K(M)$.
\qed\endproclaim

We adopt the following notation. For a {\it non-void subset} $J$
of the index set $\{1,\cdots,n\}$ let: \ss

- $P_J=\prod _{j\in J}P_j$,

- $K_J$ denote the cell structure on $P_J$ induced by
$\{K_j:j\in{}J\}$,

- $p_J$ denote the restriction to $M$ of the projection
$pr_J:P_1\times\cdots\times{}P_n\to{}P_J$ (in particular,
$p_{\{1,\cdots,n\}}:M\to{}P_1\times\cdots\times{}P_n$ is the
inclusion mapping,

- $n_J=|J|$,

- $J^c=\{1,\cdots,n\}\setminus{}J$ (therefore, $n_{J^c}=n-n_J$).

\noindent Notice that $P_{\{j\}}=P_j$, $K_{\{j\}}=K_j$ and
$pr_{\{j\}}=pr_j$. We abbreviate $p_{\{j\}}$ to $p_j$.

\ms For any cell $\sigma=\sigma_1\times\cdots\times\sigma_n\in{}K$,
$\sigma_j\in{}K_j$, the restriction $p_J|\sigma$ is the projection
onto $\sigma_J=\prod_{j\in{}J}\sigma_j$. In this sense we say that
$p_J$ {\it preserves} the cell structures $K$ and $K_J$. It follows
from Property (a) that

\proclaim{Property (b)} $p_J(M)=|K'_J|$, where $K'_J$ is a
subcomplex of $K_J$. Moreover, $|K'_J|$ is a ramified
$n_J$-manifold; if $M$ is a $($simple$)$ pseudo $n$-manifold, then
$p_J(M)$ is a $($simple$)$ pseudo $n_J$-manifold. {\rm (In the
sequel we abbreviate $K'_{\{j\}}$ to $K'_j$.)}\qed\endproclaim

\ms From this point on to the end of this section we assume that
$J$ {\it is a proper non-void subset of} $\{1,\cdots,n\}$.

\ms For every cell $\tau\in{}K'_{J^c}$, we define $P_{J}(\tau)$ to
be the union of all $n_{J}$-cells $\sigma\in{K}_{J}$ such that
$\sigma\times\tau\subset{}M$. From Property (a) we infer

\proclaim{Property (c)}
$M=\bigcup\{P_{J}(\tau)\times\tau:\tau\in{}K'_{J^c}\ is\ an\
n_{J^c}{-}cell\}$. \qed\endproclaim

\proclaim{Property (d)} If $\tau$ is a face of a cell
$\tau'\in{}K'_{J^c}$ then $P_{J}(\tau)\supset{}P_{J}(\tau')$.
\qed\endproclaim

\ms Moreover, for any point $y\in{p}_{J^c}(M)$ let
$P_{J}(y)=\{x\in{}P_{J}:(x,y)\in{}M\}$. Thus,
$P_{J}(y)\times\{y\}=p_{J^c}^{-1}(y)$, and
$P_{J}(y)\subset{}p_{J}(M)$.

\proclaim{Property (e)} For any cell $\tau\in K'_{J^c}$ and any
point $y\in\kn\tau$ we have $P_{J}(y)=P_{J}(\tau)=\bigcup
\{P_{J}(\tau'):\tau'\in{}K'_{J^c}\ is\ an\ n_{J^c}{-}cell\ with\
face\ \tau\}.$
\endproclaim

\demo{Proof} Note that
$P_{J}(y)\supset{}P_{J}(\tau)\supset{}P_{J}(\tau')$ for each cell
$\tau'\in{}K'_{J^c}$ with face $\tau$. So, it remains to justify
the inclusion
$$P_{J}(y)\subset\bigcup\{P_{J}(\tau'):\tau'\in{}K'_{J^c}\ is\
an\ n_{J^c}{-}cell\ with\ face\ \tau \}.$$ To this end, consider a
point $x\in{}P_{J}(y)$. Then $(x,y)\in{}M$. By Property (a),
$(x,y)$ belongs to an $n$-cell $\sigma\times\tau'\subset{}M$,
where $\sigma$ is an $n_{J}$-cell in $K_{J}$ and $\tau'$ is an
$n_{J^c}$-cell in $K'_{J^c}$. As $y\in\kn\tau$, $\tau$ is a face
of $\tau'$. It follows that
$x\in\bigcup\{P_{J}(\tau'):\tau'\in{}K'_{J^c}\ is\ an\
n_{J^c}{-}cell\ with\ face\ \tau \}$, which ends the proof. \qed
\enddemo

\proclaim{Property (f)}  For any cell $\tau\in{}K'_{J^c}$ the set
$P_{J}(\tau)$ is a finite disjoint union of ramified
$n_{J}$-manifolds in $P_{J}=|K_{J}|$. Moreover, if $M$ is a pseudo
$n$-manifold and $\tau$ is an $n_{J^c}$-cell then $P_{J}(\tau)$ is
a finite union of disjoint pseudo $n_{J}$-manifolds.
\endproclaim

\demo{Proof} If $\tau$ is an $n_{J^c}$-cell then both assertions
follow from the fact that each $(n{-}1)$-cell
$\sigma\times\tau\in{}K(M)$, where $\sigma\in{}K_{J}$ is an
$(n_{J}{-}1)$-cell, is a face of at least two (exactly two if $M$
is a pseudo $n$-manifold) $n$-cells $\sigma_1\times\tau,\
\sigma_2\times \tau\in{}K(M)$. Consequently, the first assertion
for arbitrary $\tau\in{}K'_{J^c}$ follows from Property (e). \qed
\enddemo

\proclaim{Property (g)} $p_{J}(M)=
\bigcup\{P_{J}(\tau):\tau\in{}K'_{J^c}\ is\ a\ k{-}cell\}$ for
each $k=0,\cdots,n_{J^c}$.
\endproclaim

\demo{Proof} For $k=n_{J^c}$ this follows from Property (c).
Applying Property (e) we obtain the general case. \qed\enddemo

\proclaim{Property (h)} If $P_{J}(w)$ is an $n_{J}$-manifold for
each vertex $w\in{}K'_{J^c}$, then $p_{J}(M)=P_{J}(w_0)$ for any
vertex $w_0$ of $K'_{J^c}$.
\endproclaim

\demo{Proof} By Property (g) (with $k=0$), it suffices to prove
that

\ss ($\ast$) $P_{J}(w_1)=P_{J}(w_2)$ for any two vertices
$w_1,w_2\in{}K'_{J^c}$. \ss

\noindent To this end, consider a 1-cell $\tau\in{}K'_{J^c}$ with
vertices $w$ and $w'$. Then, by Properties (f) and (d),
$P_{J}(\tau)$ is a finite union of ramified $n_{J}$-manifolds
contained in both $n_{J}$-manifolds $P_{J}(w)$ and $P_{J}(w')$. It
follows that $P_{J}(w)=P_{J}(\tau)=P_{J}(w')$, as no proper
ramified $n_{J}$-manifold  is contained in an $n_{J}$-manifold.
Thus the condition ($\ast$) is a consequence of the connectivity
of the complex $K'_{J^c}$. \qed
\enddemo

\proclaim{Property (i)} If $p_{J}(M)$ is a simple pseudo
$n_{J}$-manifold then $M=p_{J}(M)\times{}p_{J^c}(M)$.
\endproclaim

\demo{Proof} For every $\tau\in{}K'_{J^c}$ the set $P_{J}(\tau)$ is
a ramified $n_{J}$-manifold contained in the simple pseudo
$n_{J}$-manifold $p_{J}(M)$, so $P_{J}(\tau)=p_{J}(M)$. Hence the
assertion follows from Property(c).  \qed
\enddemo

\proclaim{Property (j)} Let $j\in \{1,\cdots, n\}$. Then $p_j(M)$ is
a circle if and only if $P_j(v)$ is a circle for each vertex $v\in
K_{\{j\}^c}'$.
\endproclaim

\demo{Proof}This follows from Property (h) combined with Property
(i).\qed\enddemo\ms

Let us define yet another symbol:
$$J(M)=\{j\in\{1,\cdots,n\}:p_j(M)~ \text{\it {is a circle}}\}.$$
If $J(M)$ is not empty then $M$ is said to {\it have projections
onto a circle}.

\proclaim{Property (k)} Suppose $J(M)$ is a proper non-void subset
of $\{1,\cdots,n\}$ then $M=p_{J(M)}(M)\times p_{J(M)^c}(M)$,
where $p_{J(M)}(M)=\prod_{j\in J(M)}p_j(M)$ is a
$|J(M)|$-dimensional torus, and $p_{J(M)^c}(M)$ $($$\subset
\prod_{j\in J(M)^c}P_j$$)$ is a ramified $|J(M)^c|$-dimensional
manifold with no projection onto a circle. If
$J(M)=\{1,\cdots,n\}$ then $M=p_1(M)\times \cdots \times p_n(M)$
is an $n$-dimensional torus.
\endproclaim

\demo{Proof} First notice that $p_{J(M)}(M)$ is a
$|J(M)|$-dimensional torus $\prod_{j\in J(M)}p_j(M)$. That follows
from the fact that $p_{J(M)}(M)$ is a ramified
$|J(M)|$-dimensional manifold (see Property (b)) lying in the
torus $\prod_{j\in J(M)}p_j(M)$ of the same dimension. Next notice
that $M=p_{J(M)}(M)\times p_{J(M)^c}(M)$ by Property (i).
Analogous arguments can be used to prove the remaining
assertions.\qed
\enddemo

\ms\centerline{5D. {\it Product structure of generalized manifolds
lying in products of graphs}}\ms

To prove the main result of this section we need the following

\proclaim{Lemma 5D.1} Let $p_i:X \to P_i$ and $q_i:S_i \to X$,
$i=1,\cdots,n$, be any mappings such that
$$H_1(p_i\circ q_j):H_{1}(S_j)\to H_1 (P_i)$$ is a monomorphism
for each $i=j$, and the zero homomorphism for $i<j$. Then the
homomorphism
$$\varphi:H_1(S_1)\oplus\cdots\oplus{}H_1(S_n)\to H_1(X),$$ defined by
$\varphi(x_1,\cdots,x_n)=H_1(q_1)(x_1)+\cdots+H_1(q_n)(x_n)$, is a
monomorphism. Consequently, $\rank{}H_1(X)\ge
\rank{}H_1(S_1)+\cdots+\rank{}H_1(S_n)$.
\endproclaim

\demo{Proof} To prove the lemma consider an element
$\alpha=(\alpha_1,\cdots,\alpha_n) \in H_1(S_1)\oplus \cdots\oplus
H_1(S_n)$ such that $\varphi(\alpha)=0$. It is enough to show that
$\alpha=0$. We prove this by finite induction. So, assume that for
some $k$, $1\le{k}\le{n}$, $\alpha_i=0$ if $i\le{}k-1$. (For $k=1$
this condition is empty.) It remains to show that $\alpha_k=0$. To
this end note that
$$H_1(p_k)(\varphi(\alpha))=
H_1(p_k)(H_1(q_k)(\alpha_k)+\cdots+H_1(q_n)(\alpha_n))=
H_1(p_k\circ{}q_k)(\alpha_k).$$ Thus
$H_1(p_k\circ{}q_k)(\alpha_k)=0$, consequently, $\alpha_k=0$,
which ends the proof. \qed \enddemo

In the following lemma we keep the notation and assumptions of
section 5C.

\proclaim{Lemma 5D.2} Let $M$ be a ramified $n$-manifold lying in
the product $P_1\times\cdots\times{}P_n$ of graphs. Then
$\rank{}H_1(M)\ge n$. If $\rank{}H_1(M)=n+k$ and $k<n$, then
$|J(M)|\geq n-k$. Therefore, for $0<k<n$, we have
$M=p_{J(M)}(M)\times p_{J(M)^c}(M)$, where
$p_{J(M)}(M)=\prod_{j\in J(M)}p_j(M)$ is a $|J(M)|$-dimensional
torus and $p_{J(M)^c}(M)$ is a ramified $|J(M)^c|$-dimensional
manifold with no projection onto a circle. For $k=0$ the set $M$
coincides with the $n$-torus $p_1(M)\times\cdots\times{}p_n(M)$.
\endproclaim

\demo{Proof} Let $v_j$, $j=1,\cdots,n$, denote a vertex of
$K'_{\{j\}^c}$. By Property (f), $P_j(v_j)$ is a finite union of
ramified 1-manifolds. To continue the proof we apply Lemma 5D.1 as
follows.

Let $q_j:P_j(v_j)\to{}M$ be the map such that $p_j\circ{}q_j$ is
the inclusion $P_j(v_j)\hookrightarrow{}P_j$ and
$(p_{\{j\}^c}\circ{}q_j)(P_j(v_j))=\{v_j\}$. Since $p_j\circ{}q_j$
is an inclusion into a graph, $H_1(p_j\circ{}q_j)$ is a
monomorphism. If $i\ne{}j$ then $H_1(p_i\circ{}q_j)$ is the zero
homomorphism as the image of $p_i\circ{}q_j$ is a point. Thus, by
Lemma 5D.1, we obtain

\ss($\ast$)
$\rank{}H_1(M)\geq\rank{}H_1(P_1(v_1))+\cdots+H_1(P_n(v_n))$.

\ss \noindent Notice that $\rank{}H_1(P_j(v_j))\geq 1$. Hence
$\rank{}H_1(M)\geq n$, which proves the first assertion.

Now we prove the second one. To this end, pick the vertices $v_j$
so that $\rank{}H_1(P_j(v_j))$ $\geq \rank{}H_1(P_j(w_j))$ for
each vertex $w_j\in{}K'_{\{j\}^c}$. Let $$J_0= \{j\in
\{1,\cdots,n\}:\rank H_1(P_j(v_j))=1\}.$$ Since
$\rank{}H_1(M)=n+k$, by ($\ast$) we infer that
$\rank{}H_1(P_j(v_j))\geq2$ for at most $k$ indices $j$. Hence
$J_0$ consists of at least $n-k$ indices. If $j\in J_0$ then
$\rank{}H_1(P_j(w_j))=1$ for each vertex $w_j\in{}K'_{\{j\}^c}$.
Hence $P_j(w_j)$ is a circle. By Property~(h) we infer that
$p_j(M)$ is a circle for each $j\in{}J_0$. It follows that $J_0
\subset J(M)$ (actually, $J_0=J(M)$). Hence $|J(M)|\geq n-k$. The
remaining assertions follow directly from Property (k).
\qed\enddemo

\proclaim{Lemma 5D.3}Let $X$ be a compactum and let $A$ be a
closed subset of $X$. If $\dim X\leq m$ and $H^m(X)=0$ then
$H^n(A)=0$ for each $n\geq m$.
\endproclaim

\demo{Proof} Since $H^n(X)=0$ and $H^{n+1}(X, A)=0$, the
conclusion follows from the cohomology exact sequence of the pair
$(X, A)$. \qed\enddemo

\proclaim{Lemma 5D.4} Let $X_i$, $i=1,2$, be non-degenerate
continua such that each point of $X_i$ admits a closed
neighborhood with trivial $n_i$-dimensional cohomology, where
$n_i=\dim X_i$. If $X_1\times{}X_2$ is a quasi
$(n_1+n_2)$-manifold then each $X_i$ is a quasi $n_i$-manifold.

\endproclaim\ms

\noi{\bf Note}. Every polyhedron $P$ satisfies the condition from
this lemma: each point of $P$ admits a closed neighborhood which
is contractible.\qed\ms

\demo{Proof} Let $x_1\in X_1$ and $x_2\in X_2$ be arbitrary
points. We must construct open neighborhoods $V_i$ of $x_i$ in
$X_i$ satisfying the definition of a quasi $n_i$-manifold. Let
$n=n_1+n_2$. Since $X_1\times{}X_2$ is a quasi $n$-manifold, there
is an open neighborhood $V$ of the point $(x_1,x_2)$ in
$X_1\times{}X_2$ such that for every closed $(n{-}1)$-dimensional
set $F$ separating $X_1\times{}X_2$ between $(x_1,x_2)$ and
$(X_1\times{}X_2)\setminus{}V$ we have $H^{n-1}(F)\ne0$. Since the
same condition holds for any open neighborhood of $(x_1,x_2)$
contained in $V$, by Lemma 5D.3 and our hypothesis about $X_i$, we
may assume that $V=V_1\times{}V_2$, where each $V_i$ is an open
neighborhood of $x_i$ in $X_i$ such that
$H^{n_i}({\overline{V}_i})=0$, where ${\overline{V}_i}$ stands for
the closure of $V_i$ in $X_i$. We shall show that these $V_i$'s
are the desired neighborhoods.

By Lemma 5D.3 again, it follows that

\ms ($\ast$) {\it for any closed subset $A$ of $X_i$ contained in
${\overline{V}_i}$ and any $k\ge{}n_i$ we have $H^{k}(A)=0$}.\ms

\noi Now consider a closed $(n_i{-}1)$-dimensional subset $F_i$ of
$X_i$ separating $X_i$ between $x_i$ and $X_i\setminus{}V_i$. Then
$X_i\setminus F_i=U_i \cup W_i$, where $U_i, W_i$ are disjoint
open sets in $X_i$ such that $x_i\in U_i$ and $X_i\setminus
V_i\subset W_i$. Then $\overline{U}_i\subset{}V_i$, and the
boundary $\partial{}U_i$ ($\subset F_i$) separates $X_i$ between
$x_i$ and $X_i\setminus{}V_i$. It follows that
$\partial(U_1\times{}U_2)$ separates $X_1\times{}X_2$ between
$(x_1, x_2)$ and $(X_1\times{}X_2)\setminus{}V$. Thus
$H^{n-1}(\partial(U_1\times{}U_2))\ne0$.

Then consider the following portion of the Mayer-Vietoris
cohomology exact sequence of the couple $\{\partial{}U_1\times
\overline{U}_2,\overline{U}_1\times\partial{}U_2\}$:
$$H^{n-2}(\partial{}U_1\times\partial{}U_2)
@>{\delta^{\ast}}>> H^{n-1}(\partial(U_1\times{}U_2)) \to
H^{n-1}(\partial{}U_1\times \overline{U}_2)\oplus{}
H^{n-1}(\overline{U}_1\times\partial{}U_2)\ .$$ (The sequence
takes this form because $\partial{}U_1\times\partial{}U_2=
(\partial{}U_1\times
\overline{U}_2)\cap{}(\overline{U}_1\times\partial{}U_2)$ and
$\partial(U_1\times{}U_2)= (\partial{}U_1\times \overline{
U}_2)\cup{}(\overline{U}_1\times\partial{}U_2)$.) By the
K{\"u}nneth formula and ($\ast$), $H^{n-1}(\partial U_1 \times
\overline{U}_2)=0$ and
$H^{n-1}(\overline{U}_1\times\partial{}U_2)=0$. It follows that
$\delta^{\ast}$ is an epimorphism. Thus
$H^{n-2}(\partial{}U_1\times\partial{}U_2)\ne0$ since
$H^{n-1}(\partial(U_1\times{}U_2))\ne0$. Again, by the K{\"u}nneth
formula and ($\ast$), $H^{n-2}(\partial U_{1} \times
\partial{}U_2)= H^{n_1-1}(\partial{}U_1)\otimes
H^{n_2-1}(\partial{}U_2)$. It follows that both
$H^{n_1-1}(\partial{}U_1)$ and $H^{n_2-1}(\partial{}U_2)$ are not
trivial.

Note that $H^{n_i}(F_i,\partial{}U_i)=0$ since $\dim{}F_i=n_i-1$.
Thus, from the cohomology exact sequence of the pair
$(F_i,\partial U_i)$ and ($\ast$), it follows that the
homomorphism $H^{n_i-1}(F_i)\to{}H^{n_i-1}(\partial{}U_i)$ induced
by the inclusion $\partial{}U_i\hookrightarrow{}F_i$ is an
epimorphism. Consequently, each $H^{n_i-1}(F_i)$ is not trivial,
which concludes the proof of our lemma. \qed \enddemo

\proclaim{Theorem 5D.5} Let $X$ be a locally connected quasi
$n$-manifold, $n\ge1$, lying in the product
$Y_1\times\cdots\times{}Y_n$ of $n$ curves. Then $\rank{}H_1(X)\ge
n$. Moreover, if $\rank{}H_1(X)=n+k$ and $k<n$ then the set
$$J(X)=\{j\in\{1,\cdots,n\}: pr_j(X)\ is\ a\ circle \},$$ where
$pr_j: Y_1\times\cdots\times Y_n \to Y_j$ stand for the
projections, contains at least $n-k$ elements and $X=(\prod_{j\in
J(X)}pr_j(X))\times{}X'$, where $X'$ is a quasi
$|J(X)^c|$-manifold in $Y_{J^c}=\prod_{i\in{J^c}}Y_i$. $($For
$k=0$ we have $X=pr_1(X)\times\cdots\times pr_n(X)$, hence $X$ is
homeomorphic to the torus $\T^n$ in this case.$)$
\endproclaim

\demo{Proof} Let
$f=(f_1,\cdots,f_n):X\to{}Y_1\times\cdots\times{}Y_n$ be the
inclusion of $X$ in the product of $n$ curves, i.e. each $f_i:X
\to{}Y_i$ is the restriction of the projection
$pr_i:Y_1\times\cdots\times{}Y_n\to{}Y_i$. Clearly, we can assume
that $\rank{}H_1(X)$ is finite. By Theorem 5B.1, there exist
mappings $g =(g_1,\cdots,g_n):X\to{}P_1\times\cdots\times{}P_n$
and $h=h_1\times\cdots\times{}h_n:P_1\times\cdots\times{}P_n\to
Y_1\times \cdots\times Y_n$ such that $f_i=h_i\circ{}g_i$ for each
$i=1,\cdots,n$ (so $f=h\circ g$), where $g_i:X\to P_i$ is a
monotone surjection, and $P_i$ is a graph.

By Theorem 5A.8, $M=g(X)$, is a ramified $n$-manifold lying in the
product $P_1 \times \cdots \times P_n$. Notice that
$\rank{}H_1(M)=\rank{}H_1(X)$.

By Lemma~5D.2, $\rank{}H_1(X)\ge n$. Moreover, if
$\rank{}H_1(X)=n+k$, where $k < n$, then $J(X)$ contains at least
$n-k$ indices. And $M=p_{J(X)}(M)\times{}p_{J(X)^c}(M)$, where
$pr_{J(X)}(M)=\prod_{j\in{J(X)}}pr_j(M)$. (If $k=0$ then
$M=pr_1(M)\times\cdots\times{}pr_n(M)$.) We have $pr_j(M)=P_j$ for
each $j\in{}J(X)$ as each $g_j$ is surjective.

Since $f=h\circ{}g$ is an inclusion it follows that $h|M$ is an
embedding of $M=\prod_{j\in{J(X)}}P_j\times{}pr_{J(X)^c}(M)$ in
the product $Y_1\times\cdots\times{}Y_n$. It follows that
$h_j:P_j\to{}Y_j$ is an embedding for each $j\in{}J(X)$ and also
$h_{J^c}|p_{J^c}(M)$ is an embedding of $p_{J^c}(M)$ into
$Y_{J(X)^c}$, where $h_{J(X)^c}=\prod_{i\in{J(X)^c}}h_i$. Observe
that $h_j(P_j)=f_j(X)$ for each $j\in{J(X)}$ and
$h_{J(X)^c}(p_{J(X)^c}(M))=f_{J(X)^c}(M)$, where
$f_{J(X)^c}=\prod_{i\in{J(X)^c}}f_i$. It follows that
$X=f_1(X)\times\cdots\times{}f_n(X)$ if $k=0$ and that
$X=(\prod_{j\in{J(X)}}f_j(X))\times{}f_{J(X)^c}(X)$ if $0<k<n$.
Also note that $f_j(X)\subset{}Y_j$ is a circle for each
$j\in{}J(X)$ and that $f_{J(X)^c}(X)\subset{}Y_{J(X)^c}$.

If $0<k<n$, by Lemma~5D.4 and Property (b), $p_{J(X)^c}(M)$ is a
quasi $(n-|J(X)|)$-manifold. It follows that $f_{J(X)^c}(M)$ is a
quasi $(n-|J(X)|)$-manifold too, as the image of $p_{J(X)^c}(M)$
by the embedding $h_{J(X)^c}|p_{J(X)^c}(M)$. This ends the proof
of the theorem. \qed\enddemo

Applying the K\"{u}nneth formula we infer the following

\proclaim{Corollary 5D.6}No product $\T^{n-k} \times M^k$, where
$2\leq k\leq n$ and $M^k$ is a closed $k$-manifold with
$H^1(M^k)=0$, can be embedded in a product of $n$ curves. \qed
\endproclaim\ms

\noindent{\bf Remark.} This corollary implies the Borsuk theorem
\cite {B3}.\ms

The "Bing house" and the "dunce hat" are contractible quasi
$2$-manifolds. Hence the first cohomology group of both examples
is trivial. Thus the theorem implies that

\proclaim{Corollary 5D.7} Both the "Bing house" and the "dunce
hat" are $2$-dimensional compact contractible polyhedra, and
neither can be embedded in a product of two curves.\qed
\endproclaim\ms

\noindent {\bf Note.} This corollary shows that the number of
factors in the Nagata embedding theorem cannot be lessened to $n$,
even for contractible $n$-dimensional polyhedra.\qed\bs

\proclaim{Corollary 5D.8} Let $X$ be an n-manifold $($in general,
let X be a locally connected quasi n-manifold which is also a
pseudo $n$-manifold$)$, $n\ge1$, lying in the product
$Y_1\times\cdots\times{}Y_n$ of $n$ curves. If $\rank{}H_1(X)\le
n+1$ then $X=S_1\times\cdots\times{}S_n$, where each $S_j$ is a
circle in $Y_j$.
\endproclaim

\demo{Proof} By Theorem 5D.5, there is a set
$J\subset\{1,\cdots,n\}$ composed of $n-1$ elements such that

\ss $(\ast) \ \ X=(\prod_{j\in{}J}S_j)\times{}X'\ $,

\ss \noi where each $S_j$ is a circle in $Y_j$ and $X'$ is a quasi
$1$-manifold in $Y_i$, where $i$ is the element of the set
$\{1,\cdots,n\}\setminus{}J$. Since $X'$ is locally connected and
$\rank{}H_1(X')\le2$ it follows that $X'$ is a graph with no
endpoint. As $X$ is a pseudo $n$-manifold, by ($\ast$) it follows
that $X'$ contains no triod, hence it is a circle. This completes
the proof.
\enddemo

\centerline{5E. {\it Contractible $2$-dimensional polyhedra in
products of two graphs}}\ms

In this section we prove a result which in a particular case gives
a noteworthy property of 2-dimensional polyhedra acyclic in
dimension 1 and embeddable in products of two curves. As neither
the "Bing house" nor the "dunce hat" have this property, we get
another argument for non-embeddability of those examples in
products of two curves.

\proclaim{Theorem 5E.1} Let $|K|$ be a $2$-dimensional connected
polyhedron, where $K$ is a regular $CW$ complex. If $|K|$ can be
embedded in a product of two curves and $\rank{}H_1(|K|)\le2$,
then $K$ collapses to either a point, or a quasi $1$-manifold, or
a torus. In particular, $K$ is collapsible if $\rank{}H_1(|K|)=0$.
\endproclaim

\noindent{\bf Remark.} Also this theorem implies the Borsuk
theorem \cite {B3}.

\ms \demo{Proof} Let
$\kappa=\{K=K_0\searrow{}K_1\searrow\cdots\searrow{}K_n\}$ be a
maximal sequence of subcomplexes of $K$ such that each successive
complex is obtained from the preceding one by an elementary
collapsing. It follows that

\ss($\ast$) $K_n$ {\it is connected and} $H_1(|K_n|)=H_1(|K|)$.

\ss Suppose $|K_n|$ is 1-dimensional. Then, since $\kappa$ is
maximal, $K_n$ has no endpoint. Hence $|K_n|$ is a quasi
1-manifold.

Now suppose $|K_n|$ is 2-dimensional. By $X$ we denote the union
of all 2-cells of $K_n$. Then $X=|K'_n|$, where $K'_n$ is a
subcomplex of $K_n$. Since $\kappa$ is maximal, each 1-cell of
$K'_n$ is a face of at least two different $2$-cells of $K'_n$.
Thus each component of $X$ is a ramified pseudo 2-manifold. Since
$|K_n|\setminus{}X$ is 1-dimensional, by exactness of the homology
sequence of $(|K_n|,X)$, it follows that
$H_1(X\hookrightarrow|K_n|)$ is a monomorphism. Consequently,
$\rank{}H_1(X)\le2$. Thus, by Theorem~5D.5, $X$ is homeomorphic to
a torus. To complete the proof it suffices to show that $X=|K_n|$.
Suppose, to the contrary, that $X$ is a proper subset of $|K_n|$.
Let $C$ be a component of the closure (in $|K_n|$) of
$|K_n|\setminus{}X$. Note that $C$ is a 1-dimensional connected
subpolyhedron of $|K_n|$ intersecting $X$ in a finite (nonzero)
number of points. Since $\kappa$ is maximal each endpoint of $C$
belongs to $X$. Observe that if $C$ has exactly one endpoint then
it contains a circle $S$, and if $C$ has at least two endpoints
then it contains an arc $L$ with end points in $X$. Note that, in
the first case $\rank{}H_1(X\cup{}S)=3$ and in the second case
$\rank{}H_1(X\cup{}L)=3$. It follows that $\rank{}H_1(|K_n|)\ge3$,
contrary to our assumption. Thus in this case $|K_n|$ is
homeomorphic to a torus.

Note that $\rank{}H_1(|K_n|)\ge1$ if $|K_n|$ is a quasi 1-manifold
and $\rank{}H_1(|K_n|)=2$ if $|K_n|$ is homeomorphic to a torus.
It follows that $|K_n|$ is a point if $\rank{}H_1(|K_n|)=0$. This
ends the proof. \qed
\enddemo

Theorem~5E.3 below is a partial inverse of Theorem~5E.1. In the
proof of 5E.3 we need the following

\proclaim{Lemma 5E.2} Let $K_1$ and $K_2$ be regular
$1$-dimensional $CW$ complexes, and let $A$ be an oriented arc in
$|K_1|\times|K_2|$ which is a union of $1$-cells of $K_1\kw{}K_2$.
Then there exist regular $1$-dimensional $CW$ complexes
$K'_1\supset K_1$ and $K'_2\supset K_2$, and a disc $D \subset
|K'_1|\times|K'_2|$, such that

\roster

\item "(i)" each component of $|K'_i|\setminus|K_i|$ is a $1$-cell of $K'_i$
with one endpoint removed,

\ss\item "(ii)" $D$ is a union of $2$-cells of $K'_1\kw{}K'_2$ and
$D\cap(|K_1|\times|K_2|)=A$.

\endroster
\endproclaim

\demo{Proof} Let $p$ denote the initial point of $A$. Without loss
of generality we may assume that $A$ is the union of arcs
$A_1,\cdots A_n$, $n\ge1$, such that: $p\in A_1$, $A_j\subset
pr_1^{-1}(v_j)$ for $j$ odd, where $v_j$ is a vertex of $K_1$,
$A_j\subset pr_2^{-1}(w_j)$ for $j$ even, where $w_j$ is a vertex
of $K_2$, $A_j$ meets $A_{j+1}$ in a common endpoint for each $j$,
and no other arcs intersect. To get $K'_1$ we enlarge $K_1$ by
$1$-cells $v'_jv_j$ for $j$ odd with all $v'_j\notin K_1$
different. Similarly, to get $K'_2$ we enlarge $K_2$ by $1$-cells
$w'_jw_j$ for $j$ even with all $w'_j\notin K_2$ different. Hence
(i) holds. Define $D$ to be the union of the following discs:
\roster
\itemitem{}  $v'_jv_j\times{}pr_2(A_j)$ for $j$ odd,

\ss\itemitem{}  $pr_1(A_j)\times{}w'_jw_j$ for $j$ even,
\ss\itemitem{}  $v'_jv_j\times{}w'_{j+1}w_{j+1}$ for $j$ odd, and
\ss\itemitem{}  $v'_{j+1}v_{j+1}\times{}w'_jw_j$ for $j$ even.
\endroster
\noindent One can easily verify that property (ii) holds as well.
\qed
\enddemo

\proclaim{Theorem 5E.3} Let $K$ be a regular $2$-dimensional $CW$
complex. If $K$ is collapsible then $|K|$ is embeddable in a
product of two trees.
\endproclaim

\demo{Proof} We shall prove a stronger version of this theorem:

\ms(0){\it there exist trees $|K_1|$, $|K_2|$ and an embedding
$h:|K|\to|K_1|\times|K_2|$ such that $h(\sigma)$ is a union of
cells of $K_1\kw{}K_2$ for each cell $\sigma\in{}K$.}\ms

Let$$|K|=|L_n|\searrow\cdots\searrow|L_0|=\{\star\}$$ be a
sequence of elementary collapses of $|K|$ to a point $\star$. The
proof of (0) will be completed ones we show that (0) holds for
each $L_m$, $m=0,\cdots, n$, in place of $K$. This in turn will be
proved by induction on $m$. Evidently, if $m=0$ the assertion (0)
is true. Now assume (0) holds for $m-1\geq 0$. We prove it for
$m$. By our assumption (0) holds for $L_{m-1}$. Hence there exist
an embedding $h':|L_{m-1}|\to|K_1|\times |K_2|$ as in (0). Since
$|L_m|\searrow |L_{m-1}|$ is an elementary collapsing, $|L_m|$ is
a union of $|L_{m-1}|$ and $\tau$, where $\tau$ is either a 1-cell
or a 2-cell of $L_m$. If $\tau$ is a 1-cell then $|L_{m-1}|\cap
\tau$ is a vertex $u_0\in L_{m-1}$. Then $h'(u_0)=(v_1, v_2)$,
where $v_i$ is a vertex of $K_i$. If $\tau$ is a 2-cell then
$|L_{m-1}|\cap\tau$ is an arc $A'$ which is a union of 1-cells of
$L_{m-1}$, so the arc $A=h'(A')$ is a union of 1-cells of
$K_1\kw{}K_2$. One easily sees that in the first case $|L_m|$
embeds in $|K_1|\times |K'_2|$ as in (0), where $K'_2=(K_2,
v_2)\vee(I, 0)$ (the one-point union). In the second case
condition (0) follows from Lemma~5E.2. This ends the proof. \qed
\enddemo

Final results of this section are devoted to embeddability of
cones over polyhedra into products.

\proclaim{Theorem 5E.4} Let $P$ be a $(k+l+1)$-dimensional
polyhedron, where $k,l\ge0$. Then there exist polyhedra $P'$ and
$P''$ with $\dim{}P'=k$ and $\dim{}P''=l$ such that the cone over
$P$ can be embedded in the product of cones over $P'$ and $P''$.
\endproclaim

The above proposition is a consequence of the following two
lemmas.

\proclaim{Lemma 5E.5} Let $P$, $k$ and $l$ be as 5E.4. Then there
exist polyhedra $P'$ and $P''$ with $\dim{}P'=k$ and $\dim{}P''=l$
such that $P$ can be embedded in the join $P'*P''$.
\endproclaim

\demo{Proof} Let $P=|K|$, where $K$ is a simplicial complex. Put
$P'=|K^{(k)}|$, where $K^{(k)}$ is the $k$-skeleton of $K$. Define
$P''$ to be {\it the dual} to $P'$ in $P$, i.e. $P''$ is the union
of all simplices of the barycentric subdivision of $K$ which are
disjoint with $P'$. Then $\dim{P'}=k$ and $\dim{P''}=l$. Observe
that $P$ is pl isomorphic to a subpolyhedron of the join
$P'\ast{P''}$. This follows from the fact that for any simplex
$\sigma\in K$ with $\dim\sigma\ge{k}$ we have
$\sigma=\sigma'\ast\sigma''$, where $\sigma'$ is the $k$-skeleton
of $\sigma$ (with respect to the standard simplicial structure on
$\sigma$) and $\sigma''$ is the dual to $\sigma'$ in $\sigma$.
\qed
\enddemo

\proclaim{Lemma 5E.6} Let $P'$ and $P''$ be polyhedra. Then the
product of cones over $P'$ and $P''$ is homeomorphic to the cone
over the join $P'*P''$.
\endproclaim

\demo{Proof} For a polyhedron $Q$ let $aQ$ denote the cone with
vertex $a$ and base $Q$. According to the definition of link (cf.
\cite{R-S, p. 2}), $Q$ may be considered as a link of the vertex
$a$ in $aQ$. The conclusion of our lemma is the following special
case of a known formula (see \cite{R-S, p. 24, Ex. (3)}) for the
link of a point in the product of two polyhedra: $\lk((a', a''),
a'P'\times a''P'')\cong P'*P''$.\qed
\enddemo

\proclaim{Corollary 5E.7} The cone over any $n$-dimensional
polyhedron can be embedded in the product of $n+1$ copies of an
$m$-od.
\endproclaim

\demo{Proof} We can prove this result by induction on the
dimension $n$ of the polyhedron. We start the induction with
$0$-dimensional polyhedra, in which case the proof is obvious. The
inductive step is proven applying Theorem 5E.4 for $k=0$ and
$l=n-1$.\qed
\enddemo

Note that neither Lemma~5E.5 nor Corollary~5E.7 (and thus
Theorem~5E.4) can be extended to a more general case of continua
in place of polyhedra. For example the Menger curve can not be
embedded in the join of two 0-dimensional compacta. Neither the
cone over the Menger curve is embeddable in the product of two
cones over 0-dimensional compacta.

\bs\centerline {\bf Problems to Chapter 5}

\proclaim {Problem 5A.1} Is it possible to characterize a
polyhedron $|K|$ which is a quasi $n$-manifold in terms of the
complex $K$ itself?
\endproclaim

\proclaim {Problem 5B.1} Let $X$ be a locally connected quasi
$n$-manifold, $n\geq2$, with $H^1(X)$ of finite rank, and let
$(f_1,\cdots,f_n):X \to Y_1\times\cdots\times Y_n$ be an embedding
in a product of curves. Is it possible to approximate mappings
$f_i$ by mappings $f_i':X \to Y_i$ so that $(f_1',\cdots,f_n')$ is
still an embedding and each $f_i'(X)$ is a graph?\endproclaim

\proclaim{Problem 5D.1} Characterize quasi $n$-manifolds
embeddable in products of $n$ graphs.\endproclaim

\proclaim {Problem 5D.2}Let X be a locally connected pseudo
n-manifold lying in a product of n curves. Is X a locally connected
quasi n-manifold?
\endproclaim

\proclaim{Problem 5D.3} Characterize ordinary closed
$3$-manifolds. Must such a manifold be a product of non-degenerate
factors?\endproclaim\ms

\head 6. Embedding surfaces into products of two curves
\endhead

This chapter splits into sections 6A - 6G. In section 6A we
present a short and elementary proof of the Borsuk theorem on
non-embeddability of $\Sb^2$ in products of two curves. In section
6B we establish some general properties of closed surfaces lying
in products of two graphs. In section 6C we show that any
orientable surface of genus $g$ can be monotonically embedded in
the product of two $\theta_{g+1}$-curves (see section 6B for
definitions). The main objective of the entire chapter is
discussed in section 6D devoted to certain results of R.~Cauty
\cite{C} concerning embeddability of non-orientable surfaces in
products of curves. In section 6E we prove that among closed
surfaces lying in a product of two curves the torus is the only
one which is a retract of that product. In section 6F we show that
any surface with non-empty boundary can be embedded in the
"three-page book". In section 6G we consider embeddability in the
second symmetric product of curves.

\bs \centerline{6A. {\it Another proof of the Borsuk theorem}}

\ms Suppose there is an embedding $f:\Sb^2\to Y_1\times Y_2$,
where $Y_1$ and $Y_2$ are curves. Then $f=(f_1,f_2)$, where
$f_i:\Sb^2\to Y_i$, $i=1,2$. By the Whyburn factorization theorem,
there exist a continuous monotone surjection $g_i:\Sb^2\to X_i$
and a 0-dimensional mapping $h_i:X_i\to Y_i$ such that
$f_i=h_i\circ g_i$. As $g_i$ is continuous and monotone, $X_i$ is
a locally connected continuum contractible with respect to $\Sb^1$
(equivalently: $H^1(X_i)=0$) \cite{Kur, \S57, I, Theorem 2, p.
434}. As $h_i$ is 0-dimensional, $X_i$ is 1-dimensional by the
Hurewicz theorem \cite{Kur, \S45, VI, Theorem 1, p. 114}. It
follows that $X_i$ is an absolute retract (a dendrite) (see e.g.
\cite{Kur, \S57, III, Corollary 8, p. 442 and \S53, IV, Theorem
16, p. 344 }). Then we can continue in two different ways:

(a) Since $\dim(Y_1\times Y_2)=2$, there is an extension
$\varphi:Y_1\times Y_2\to\Sb^2$ of the inverse
$f^{-1}:f(\Sb^2)\to\Sb^2$ (see e.g. \cite{Kur, \S53, VI, Theorem
1, p. 354}). Note that $f$ is homotopic to a constant mapping,
because each $f_i$ is (as $X_i$ is contractible). It follows that
$\varphi\circ f:\Sb^2\to\Sb^2$ is homotopic to a constant. But
$\varphi\circ f=1_{\Sb^2}$, a contradiction.

(b) Note that $X_1\times X_2$ is a 2-dimensional absolute retract.
On the other hand, one easily sees that $(g_1,g_2):\Sb^2\to
X_1\times X_2$ is an embedding. Thus $S=(g_1,g_2)(\Sb^2)$ is a
topological 2-sphere in $X_1\times X_2$. It follows that $S$ is a
retract of $X_1\times X_2$ (see e.g. \cite{Kur, \S53, VI, Theorem
1, p. 354}), hence $S$ is an absolute retract as well, a
contradiction. \qed

\bs\centerline{6B. {\it On surfaces lying in products of two
graphs}} \ms

This section is devoted to general properties of closed surfaces
lying in products of two graphs. We apply the notation and
conventions of Section 5C for $n=2$. Below we rephrase the general
Properties (a)-(j) from section 5C for this special dimension.
Therefore, no proof of the obtained properties will be given. (The
restated properties will be designated by adding primes to the
letters designating corresponding properties from Section 5C.)\ms

\centerline{{\it Throughout this section we consider a special
pseudo $2$-manifold $M\subset{}P_1\times{}P_2$,}} \centerline{{\it
where $P_i=|K_i|$ are graphs and $K_i$ are regular $1$-dimensional
$CW$ complexes.}}
\ms \noi From Theorem 5A.4 we infer that
$M=|K(M)|$, where $K(M)$ is a subcomplex of $K_1\kw{}K_2$.
Moreover, $K(M)$ is a special pseudo 2-manifold complex, hence
{\it any ramified $2$-manifold subcomplex of $K(M)$ coincides
with} $K(M)$. Let us recall that $p_i=pr_i|M$.\ms

\proclaim{Property (a$'$)} $M$ is the union of $2$-cells of
$K(M)$. \qed\endproclaim

\proclaim{Property (b$'$)} $p_i(M)=|K_i'|$, {\it where} $K_i'$
{\it is a subcomplex of} $K_i$ {\it with no endpoint}. \qed
\endproclaim

\ms For any cell $\tau\in{}K'_2$ the set $P_1(\tau)$ is the union
of all $1$-cells $\sigma\in{K}_1$ such that
$\sigma\times\tau\subset{}M$. (The set $P_2(\sigma)$ for
$\sigma\in K_1$ has analogous description.) For any point
$y\in{p}_{2}(M)$ the set $P_{1}(y)=\{x\in{}P_{1}:(x,y)\in{}M\}$.
Thus, $P_{1}(y)\times\{y\}=p_{2}^{-1}(y)$, and
$P_{1}(y)\subset{}p_{1}(M)$. (The set $P_2(x)$ for $x\in p_1(M)$
has analogous description.) We list certain properties of the
sets $P_1(\tau)$ and $P_1(y)$; analogous properties are valid for
the other pair of sets.)

\proclaim{Property (c$'$)}
$M=\bigcup\{P_1(\tau)\times\tau:\tau\in{}K'_2\ is\ a \
1{-}cell\}$. \qed\endproclaim

\proclaim{Property (d$'$)} If $w$ is a vertex of a $1$-cell
$\tau\in{}K'_2$ then $P_1(w)\supset{}P_1(\tau)$. \qed\endproclaim

\proclaim{Property (e$'$)} For any cell $\tau\in K'_{2}$ and any
point $y\in\kn\tau$ we have $P_{1}(y)=P_{1}(\tau)$. If
$\tau=\{v\}$ is a vertex then $P_1(\tau)=\bigcup
\{P_1(\tau'):\tau'\in{}K'_{2}\ is\ a\ {1}{-}cell\ with\ vertex\
v\}.$ \qed\endproclaim

\proclaim{Property (f$'$)}  For any $1$-cell $\tau\in{}K'_{2}$ the
set $P_{1}(\tau)$ is a finite union of disjoint circles.
\qed\endproclaim

\proclaim{Property (g$'$)} $p_{1}(M)=
\bigcup\{P_{1}(\tau):\tau\in{}K'_{2}\ is\ a\ k{-}cell\}$ for $k=0,
1$. \qed\endproclaim

\proclaim{Property (h$'$)} If $P_{1}(w)$ is a circle for each
vertex $w\in{}K'_{2}$, then $p_{1}(M)=P_{1}(w_0)$ for any vertex
$w_0$ of $K'_{2}$. \qed\endproclaim

\proclaim{Property (i$'$)} If $p_{1}(M)$ is a circle then
$M=p_{1}(M)\times{}p_{2}(M)$. \qed\endproclaim

\proclaim{Property (j$'$)} $p_1(M)$ is
a circle if and only if $P_1(v)$ is a circle for each vertex $v\in
K_2'$.
\endproclaim

\ms If both $p_i:M\to P_i$ are surjections (i.e. $p_i(M)=P_i$, or
$K_i'=K_i$), then $M$ is said to be {\it surjectively embedded} in
$P_1\times P_2$.\ms

For any two mappings $h_i:P_i\to Y_i$, where $Y_i$ is a curve, we
have the following

\proclaim{Property (l)} {\it Suppose} $M$ {\it is surjectively
embedded in} $P_1 \times P_2$. {\it If} $h_1\times h_2:P_1\times
P_2\to Y_1 \times Y_2$ {\it is injective on} $M$ {\it and every
two} $P_2(\sigma), \sigma \in K_1$, {\it intersect, then} $h_1:P_1
\to Y_1$ {\it is an embedding}.
\endproclaim

\demo{Proof} For suppose $h_1(x)=h_1(x')$ for some $x, x'\in P_1$;
we have to show that $x=x'$. Let $x \in \kn\sigma$ and
$x'\in\kn\sigma'$ for some $\sigma,\sigma'\in K_1$. Then
$P_2(x)=P_2(\sigma)$ and $P_2(x')=P_2(\sigma')$ by Property
(e$'$). By our assumption, there is a point $y \in P_2(\sigma)\cap
P_2(\sigma')=P_2(x)\cap P_2(x')$. Then points $(x,y),(x',y)$ are
different, both belong to $M$ and are sent to $(h_1(x), h_2(y))$
by $h_1 \times h_2$, a contradiction. This proves Property
(l).\qed
\enddemo

\ss Let $\Theta_n$, $n\ge1$, denote the graph in $\Sb^2$ composed
of $n$ meridians $\mu_0^n, ..., \mu_{n-1}^n$ such that each
$\mu_j^n$ passes through $z_j=(\cos2\pi{j\over n}, \sin2\pi{j\over
n},0)$. (So, the meridians are equally and cyclicly spaced on
$\Sb^2$; if no confusion is likely to occur, we abbreviate
$\mu_j^n$ to $\mu_j$.) The poles $p=(0,0,1)$ and $-p=(0,0,-1)$ are
common endpoints of the meridians. Hence $\Theta_n=|K(n)|$, where
$K(n)$ denotes the regular 1-dimensional $CW$-complex with the
meridians as 1-cells (edges) and with the poles as 0-cells
(vertices). Any space homeomorphic to $\Theta_n$ will be called a
$\theta_n$-{\it curve}, its points corresponding to the poles of
$\Theta_n$ will be called {\it poles} of the curve. Note that a
$\theta_3$-curve is commonly known as the $\theta$-{\it curve}.

\proclaim{Property (m)} {\it Suppose $M$ is surjectively embedded
in $P_1 \times P_2$, where $P_2=\Theta_{n+1}$, $n\ge1$. Then, for
any $1$-cell $\sigma\in K_1$, the set $P_2(\sigma)$ is a circle
$($composed of two meridians$)$. Analogous property holds in case}
$P_1=\Theta_{m+1}$. \qed\endproclaim

\ms In case where both $p_1:M\to P_1$ and $p_2:M\to P_2$ are
monotone, $M$ is said to be {\it monotonically
embedded} in $P_1\times P_2$.

\bs\centerline {6C. \it{Surfaces in products of}
$\theta_n$-\it{curves}}

 \noi\bs We keep notation of section 6B. First we prove that
any orientable surface can be monotonically embedded in the
product of two $\theta_n$-curves. Applying some results of the
previous article we show that any two copies of an orientable
surface lying in the product of two $\theta_n$-curves are
equivalent up to a special automorphism of the entire product.

For each natural $m\ge1$ define $M^0_{2m}$ to be the subset of
$\Theta_{m+1} \times\Theta_{m+1}$ given by
$$(\alpha)\;\;\;\;\; M^0_{2m}= \bigcup_{j
=0}^{m}\mu_j^{m+1}\times(\mu_{j}^{m+1}\cup\mu_{j+1}^{m+1}).$$ \noi
Here $\mu_{m+1}^{m+1}=\mu_{0}^{m+1}$. One easily sees that this
space can be also expressed in the "dual" form:
$$(\beta)\;\;\;\;\;M^0_{2m}= \bigcup_{j = 0}^{m}(\mu_{j-1}^{m+1}\cup
\mu_{j}^{m+1})\times \mu_{j}^{m+1}.$$ Here
$\mu_{-1}^{m+1}=\mu_{m}^{m+1}$. Let $p_i:M^0_{2m}\to \Theta_{m+1}$
denote the restriction to $M^0_{2m}$ of the projection
$pr_i:\Theta_{m+1} \times\Theta_{m+1}\to \Theta_{m+1}$. The
mappings $q_1:M^0_{2m}\to \Theta_{m+1}\times\{-p\}$ and
$q_2:M^0_{2m}\to \{-p\}\times\Theta_{m+1}$ defined by
$q_1(z_1,z_2)=(z_1,-p)$ and $q_2(z_1,z_2)=(-p,z_2)$, for all
$(z_1,z_2)\in M^0_{2m}$, are monotone retractions (equivalent to
the corresponding mappings $p_i$) onto $\theta_{m+1}$-curves. And
their diagonal mapping
$$(q_1,q_2):M^0_{2m}\to (\Theta_{m+1}\times\{-p\})\times
(\{-p\}\times\Theta_{m+1})\;\;(\approx \Theta_{m+1}
\times\Theta_{m+1})$$ is an embedding (equivalent to the inclusion
$(p_1,p_2):M^0_{2m}\to \Theta_{m+1} \times\Theta_{m+1})$. The
properties of $q_i^{\prime} s$ follow from corresponding properties of
$p_i^{\prime} s$ (see the proof below).

\proclaim{Lemma 6C.1} The set $M^0_{2m}$ is a closed orientable
surface of genus $m$ monotonically embedded in
$\Theta_{m+1}\times\Theta_{m+1}$.
\endproclaim

\demo{Proof} By definition it follows that both $p_i$ are
surjective. Hence the monotonicity of the embedding is a consequence
of the equalities:
$p_1^{-1}(x)=\{x\}\times(\mu_{j}^{m+1}\cup\mu_{j+1}^{m+1})$ for
any $x\in\kn\mu_{j}^{m+1}$,
$p_2^{-1}(y)=(\mu_{j-1}^{m+1}\cup\mu_{j}^{m+1})\times\{y\}$ for
any $y\in\kn\mu_{j}^{m+1}$; $p_1^{-1}(v)=\{v\}\times\Theta_{m+1}$,
$p_2^{-1}(w)=\Theta_{m+1}\times\{w\}$ for $v,w\in\{p,-p\}$. So, to
complete the proof, it remains to show that:\ms

($\ast_m$) {\it $M^0_{2m}$ is a closed orientable surface of genus
$m$.}

\ms \noi We give two proofs of this assertion.\ss

{\it First proof of} ($\ast_m$). This proof is by induction on
$m$. Clearly, $(\ast_1)$ holds because $M^0_2=\mu_0^2
\times(\mu_0^2 \cup\mu_1^2 )\cup\mu_1^2 \times(\mu_1^2
\cup\mu_0^2)=\Theta_2\times\Theta_2$ is a torus. Then suppose
$(\ast_m)$, $m\ge 1$, has been proved; it remains to show that
$(\ast_{m+1})$ holds too. To this end, observe that
$M^0=\mu_0^{m+2}
\times(\mu_{0}^{m+2}\cup\mu_{1}^{m+2})\cup\mu_{1}^{m+2}\times(\mu_{1}^{m+2}\cup\mu_{2}^{m+2})\cup
\cdots\cup\mu_{m}^{m+2}\times(\mu_{m}^{m+2}\cup\mu_{0}^{m+2})$ is
a topological copy of $M_{2m}^0$ in
$\Theta_{m+2}\times\Theta_{m+2}$. Hence, by the inductive
assumption, $M^0$ is a closed orientable surface of genus $m$. Let
$T$ denote the torus
$(\mu_{m}^{m+2}\cup\mu_{m+1}^{m+2})\times(\mu_{m+1}^{m+2}\cup\mu_{0}^{m+2})$.
Then note that the intersection $M^0\cap T$ is the disc
$\mu_m^{m+2}\times\mu_{0}^{m+2}$, and
$M^0_{2(m+1)}=(M^0\cup{}T)\setminus\kn\mu_{m}^{m+2}\times\kn\mu_{0}^{m+2}$.
It follows that $M^0_{2(m+1)}=M^0\#T$, so $M^0_{2(m+1)}$ is an
orientable surface of genus $m+1$, which proves $(\ast_{m+1})$.\ss

{\it Second proof of }($\ast_m$). This argument appeals to the
following obvious fact: If $N$ is a 2-disc with $m$ holes then the
doubling $M=2N$ is an orientable closed surface of genus $m$. Let
us recall that by definition the {\it doubling} $2N$ is the union
of two copies of $N$ glued along $\partial N$. In other words,
$2N$ is the subset of the product $N\times I$ defined by
$$2N=\partial(N\times I)=(\partial N)\times I \cup N\times\partial I.$$
So, it suffices to show that $M^0_{2m} \approx 2N$. Actually, we
shall prove more: $M^0_{2m}$ can be obtained as the image of a
special embedding $\varphi:2N\to\Theta_{m+1}\times\Theta_{m+1}$.

To this end we take $N$ quite specific. Notice that
$\Sb^2=D_0\cup\cdots\cup D_{m}$, where each $D_j$ stands for the
2-disc in $\Sb^2$ from $\mu_{j-1}^{m+1}$ to $\mu_j^{m+1}$ (that
is, $D_j$ consists of the points which lie on meridians passing
through points $(\cos 2\pi t, \sin 2\pi t,0)$ with $t\in
[\frac{j-1}{m+1}, \frac{j}{m+1}]$). (To make the notation less
complicated we abbreviate $\mu_j^{m+1}$ to $\mu_j$). For each $j$
choose a 2-disc $E_j$ lying in the interior of $D_j$, and put
$C_j=D_j\setminus\kn E_j$. Then each triple $(C_j,\partial
D_j,\partial E_j)$ is homeomorphic to the cylinder triple
$(\Sb^1\times I,\Sb^1\times \{0\}, \Sb^1\times \{1\})$. Now,
define $N$ to be the set
$$N=\Sb^2 \setminus (\kn E_0\cup\cdots\cup\kn E_{m}).$$ \noi
One easily sees that $N=C_0\cup\cdots\cup C_m$ and $C_j\cap
C_{j+1}=\mu_j$ for each $j$.

On the other hand we have
$$2N=\bigcup_{j=0}^{m}(C_j\times \{0,1\})\cup (\partial E_j \times I).$$
\noi Notice that, for each $j$, both sets
$(C_j\times\{0,1\})\cup(\partial E_j\times I)$ and $(\mu_{j-1}\cup
\mu_{j})\times \mu_{j}$ are homeomorphic to the cylinder
$\Sb^1\times I$. Moreover, there exist homeomorphisms
$\varphi_j:(C_j\times\{0,1\})\cup(\partial E_j\times I)\to
(\mu_{j-1}\cup \mu_{j})\times \mu_{j}$ such that
$\varphi_j(z,0)=(z,-p)$ and $\varphi_j(z,1)=(z,p)$ for each
$z\in\mu_{j-1}\cup \mu_{j}$. Combining these homeomorphisms we get
the desired embedding $\varphi$. (See Figure 6C.1 where the doubling
$2N$, with $m=2$ and a special embedding in $\R^3$, is depicted.)
\qed\enddemo\ms

Any sequence $\nu_0,\cdots,\nu_k$ of different meridians of $\Theta_n$ induces
a sequence $$S_0,\cdots,S_k$$ of topological circles in $\Theta_n$,
where $S_j=\nu_j\cup\nu_{j+1}$ for $j=0,\cdots,k-1$ and $S_k=\nu_k\cup\nu_0$. This sequence
will be called {\it a cycle of circles induced by} $\nu_0,\cdots,\nu_k$.
The following assertion immediately follows from this definition and Lemma 6C.1.\bs

\proclaim{Theorem 6C.2} Let $M$ be any subset of the product
$\Theta_n\times\Theta_{n'}$, $n,n'\ge2$, which can be expressed in
the form
$$M= \sigma_0\times S_0\cup\cdots\cup\sigma_m \times S_m,$$ where
$\sigma_0,\cdots,\sigma_m$ are different meridians in $\Theta_n$,
and $S_0,\cdots,S_m$ is a cycle of circles in $\Theta_{n'}$. Then
M is a closed orientable surface of genus $m$.

Moreover, for any embeddings $h:\Theta_{m+1}\to\Theta_n$ and
$h':\Theta_{m+1}\to\Theta_{n'}$ $($preserving the poles$)$ such
that $h(\mu_j^{m+1})=\sigma_j$ and $h'(\mu_j^{m+1})=\tau_j$, for
each $j=1,\cdots,m$, where $\tau_0,\cdots,\tau_m$ induces the
cycle $S_0,\cdots,S_m$, the embedding
$h\times{}h':\Theta_{m+1}\times\Theta_{m+1}\to\Theta_n\times\Theta_{n'}$
takes $M^0_{2m}$ onto M. \qed\endproclaim

\ms\proclaim{Corollary 6C.3} For each $m\ge1$ the product
$\Theta_{m+1}\times \Theta_{m+1}$ contains a closed orientable
surface M of genus m in such a position that it is invariant under
the canonical involution on that product. Moreover, M meets the
diagonal of $\Theta_{m+1}\times \Theta_{m+1}$ along a circle which
does not separate M.
\endproclaim

\demo{Proof} We shall show that the subset $M$ of
$\Theta_{m+1}\times \Theta_{m+1}$ given by the formula (here we
write $\mu_j$ instead of $\mu_j^{m+1}$)
$$M=(\mu_0\times \mu_0)\cup
(\bigcup_{i=0}^{m-1}(\mu_i\times \mu_{i+1}\cup\mu_{i+1}\times
\mu_i))\cup(\mu_m\times \mu_m).$$ \noi has the desired properties.
Obviously, $M$ is symmetric under the canonical involution on
$\Theta_{m+1}\times \Theta_{m+1}$. Moreover, $M$ meets the
diagonal of $\Theta_{m+1}\times \Theta_{m+1}$ along the circle
$S=D_0\cup D_m$, where $D_j$ denotes the diagonal of
$\mu_j\times\mu_j$. One can easily verify that $M\setminus S$ is
connected (the vertices $(-p,p)$ and $(p,-p)$ can be connected by
an arc lying in $\mu_0\times\mu_1 \setminus S$ and any other point
of $M\setminus S$ can be connected to one of these points by an
arc lying in $M\setminus S$). Thus, by Theorem 6C.2, it remains to
show that $M$ can be expressed in the form
$$M= \sigma_0\times (\tau_0\cup\tau_1)\cup\cdots\cup\sigma_m
\times (\tau_m\cup\tau_0),$$ where $\sigma_0,\cdots,\sigma_m$ and
$\tau_0,\cdots,\tau_m$ are sequences of different meridians in
$\Theta_{m+1}$. The constructions of the sequences $(\sigma)$ and
$(\tau)$ depend on parity of $m$.

If {\it m is even} then we put (instead of $\mu_i^{m+1}$'s we
write just their subscripts)

\ms\hskip1.7cm $(\sigma)=(0,1,3,\cdots, m-1,m,m-2,m-4,\cdots,
4,2),$

\ms\hskip1.7cm$(\tau)=(1,0,2,4,\cdots, m-2,m,m-1,m-3,\cdots,
5,3).$

\ms\noi If {\it m is odd} then we put

\ms\hskip1.7cm $(\sigma)=(0,1,3,\cdots,m-2,m,m-1,m-3,\cdots,
4,2),$

\ms\hskip1.7cm $(\tau)=(1,0,2,4,\cdots, m-1,m,m-2,m-4,\cdots,
5,3).$

\ms\noi One can verify that these sequences have the desired
properties.\qed\enddemo

\proclaim{Theorem 6C.4} Let $M$ be a special pseudo $2$-manifold
with $\rank H_1(M)=2m$ lying in the product
$\Theta_n\times\Theta_{n'}$, $n,n'\ge2$. Then M can be expressed
in the form
$$M=\sigma_0\times P_2(\sigma_0)\cup\cdots\cup\sigma_m \times P_2(\sigma_m),$$
where $\sigma_0,\cdots,\sigma_m$ is a sequence of different
meridians in $\Theta_n$, and the sets
$P_2(\sigma_0),\cdots,$ $P_2(\sigma_m)$ form a cycle of circles in
$\Theta_{n'}$. In particular, $M$ is an orientable closed surface
of genus $m$.
\endproclaim

\demo{Proof} Consider a sequence $\sigma_0\cdots\sigma_k$ of
length $k+1$, with $k\ge0$, such that

\ss(1) $\sigma_0\cdots\sigma_k$ {\it are different} 1-{\it cells
of} $K'(n)$,

\ss(2) $P_2(\sigma _j)=\tau_j\cup\tau_{j+1}$ {\it for each}
$j=0,\cdots,k$,

\ss(3) $\tau_0\cdots\tau_k$ {\it are different} 1-{\it cells of}
$K'(n').$

\ss\noi (Such sequences exist: the one-element sequence
$\sigma_0$, for any 1-cell $\sigma_0$ of $K'(n)$, has these
properties. In fact, by Property (f$'$) the set $P_2(\sigma_0)$ is
a circle. Hence $P_2(\sigma_0)=\tau_0\cup\tau_1$ and
$\tau_0\ne\tau_1$.) We investigate the union
$$M[k]=\bigcup_{j=0}^{k}\sigma_j\times(\tau_j\cup\tau_{j+1}).$$
(The 2-cells of $K(M)$ lying in $M[k]$ make the sequence
$$\sigma_0\times\tau_0,\ \sigma_0\times\tau_1,\
\sigma_1\times\tau_1,\ \cdots,\ \sigma_k\times\tau_k,\
\sigma_k\times\tau_{k+1}$$ such that each cell borders upon the
next one along in turn $$\sigma_0\times\{-p,p\},
\{-p,p\}\times\tau_1,\cdots,\sigma_k\times\{-p,p\}.)$$

In subsequent discussion we take for $k$ the maximal number for
which the conditions (1)-(3) hold. Then we have \ss (4)
$\tau_{k+1}=\tau_0$. \ss\noi For suppose $\tau_{k+1}\neq\tau_0$.
Then $\tau_{k+1}\neq\tau_j$ for each $j=0,\cdots,k$. In fact, else
$\tau_{k+1}=\tau_j$ for some $0<j<k$. Thus
$P_1(\tau_j)=\sigma_{j-1}\cup\sigma_j\cup\sigma_k$. So, by (1),
$P_1(\tau_j)$ is not a circle, contrary to Property(m). Next, let
$P_1(\tau_{k+1})=\sigma_k\cup\sigma_{k+1}$. Then
$\sigma_{k+1}\neq\sigma_j$ for each $j$. Else
$\sigma_{k+1}=\sigma_j$ for some $0\leq j<k$. Consequently,
$P_2(\sigma_j)=\tau_j\cup\tau_{j+1}\cup\tau_{k+1}$ is not a
circle, because $\tau_0,\cdots,\tau_{k+1}$ are all different. This
again contradicts Property (m). Combining these properties, one
sees that $\sigma_0,\cdots, \sigma_{k+1}$ is a sequence of length
$k+2$ satisfying (1)-(3), contrary to maximality of $k$.

It follows from (1)-(4) that
$$M[k]=\sigma_0\times P_2(\sigma_0)\cup\cdots\cup\sigma_k \times P_2(\sigma_k)$$
is a subset of $\Theta_n\times\Theta_{n'}$ of the form which has been
discussed in Theorem 6C.2. Consequently,
$M[k]$ is a closed orientable surface of genus $k$ contained in $M$.
This implies $M[k]=M$. In addition, $k=m$ because $2m=\rank H_1(M)=\rank H_1(M[k])=2k$,
which completes the proof. \qed
\enddemo

\ss The following corollary is a special case of the above theorem.

\proclaim{Corollary 6C.5} Let $M$ be a closed surface of genus $m$
lying in the product $\Theta_n\times\Theta_{n'}$, $n,n'\ge2$. Then
M can be written in the form
$$M=\sigma_0\times P_2(\sigma_0)\cup\cdots\cup\sigma_m \times P_2(\sigma_m),$$
where $\sigma_0,\cdots,\sigma_m$ is a sequence of different meridians in $\Theta_n$,
and $P_2(\sigma_0),\cdots,P_2(\sigma_m)$ is a cycle of circles in $\Theta_{n'}$.
In particular, $M$ is orientable.\qed
\endproclaim\ms

Two final assertions directly follow from the above discussion.

\ms\proclaim{Corollary 6C.6}Let $M$ be a closed surface of genus
$m$ lying in the product $\Theta_n\times\Theta_{n'}$, $n,n'\ge2$.
Then both $p_1(M)$ and $p_2(M)$ are $\theta_{m+1}$-curves and $M$
is monotonically embedded in $p_1(M)\times p_2(M)$
.\qed\endproclaim

\ss\proclaim{Corollary 6C.7} Let $M_1$ and $M_2$ be two copies of
a closed surface lying in the product $\Theta_n\times\Theta_{n'}$,
$n,n'\ge2$. Then there are homeomorphisms $h:\Theta_n\to\Theta_n$,
and $h':\Theta_{n'}\to\Theta_{n'}$ such that $M_2=(h\times
h')(M_1)$.\qed
\endproclaim

\ss \bs \centerline{6D. {\it On Cauty's results about
embeddability of non-orientable surfaces}} \centerline {\it into
products of graphs}

\ms In this section we mostly comment on some beautiful results of
Cauty \cite{C1} concerning embeddability as in the title. In
particular, we give slightly different and hopefully simpler
descriptions of the embeddings occurring in the Cauty proof. The
main results are essentially due to Cauty.

\ms\proclaim{Theorem~6D.1} {\rm (cf. \cite{C1})}  Let $M$ be a
special pseudo $2$-manifold in the product $Y_1\times Y_2$ of two
curves. Then

{\rm(i)} $\rank H_1(M)\le3$ implies $M=P_1\times P_2$, where each
$P_i$ is a circle in $Y_i$;

{\rm(ii)} $\rank H_1(M)=4$ implies $M\subset P_1\times P_2$, where
each $P_i$ is a $\theta$-curve in $Y_i$. $($Hence $M\approx
\T^2\#\T^2$ and $M$ is monotonically embedded in $P_1\times
P_2$$)$.
\endproclaim

\noi{\bf Remark.} Note that (i) implies the Borsuk theorem
\cite{Bo3}.

\ms In order to prove Theorem 6D.1 we need the following three
lemmas. The first one will be used in the proof of Lemma 6D.4.

\proclaim {Lemma~6D.2} Let $S_j$ be a closed subset of a graph
$P_j$, $j=1,2$, and let $S_0$ be either the empty set or a circle
in $P_1$ such that $S_0\setminus{}S_1\neq\emptyset$. Suppose
$p_j:M\to{}P_j$ and $q_i:S_i\to{}M$, $i=0,1,2$, are mappings
satisfying the conditions:

{\rm(i)} $H_1(p_1\circ{}q_2)=0$,

{\rm(ii)} $p_j\circ{}q_j:S_j\to{}P_j$, and
$p_1\circ{}q_0:S_0\to{}P_1$ are the inclusions.

\noindent Then there is a monomorphism
$$\varphi:H_1(S_0)\oplus{}H_1(S_1)\oplus{}H_1(S_2)\to{}H_1(M).$$
\endproclaim

\demo{Proof} We will prove that the homomorphism $\varphi$ defined
by
$$\varphi(x_0,x_1,x_2)=H_1(q_0)(x_0)+H_1(q_1)(x_1)+H_1(q_2)(x_2)$$
is a monomorphism.

Since $H_2(P_j,S_j)=0$, by (ii) and the exactness of the homology
sequence of $(P_j,S_j)$ we have

(1) $H_1(p_j\circ{}q_j)$ {\it is a monomorphism}.

First, assume that $S_0=\emptyset$. Note that $H_1(S_0)$ is
trivial. Thus, by (i) and (1), Lemma 4D.1 implies that $\varphi$
is a monomorphism.

Now, assume $S_0$ is a circle in $P_1$ such that
$S_0\setminus{}S_1\neq\emptyset$. Then there is a retraction
$r:P_1\to{}S_0$ such that $r(S_1)$ is a proper subset of $S_0$.
Define $p_0=r\circ{}p_1$. Note that

(1$'$) $H_1(p_0\circ{}q_0)$ {\it is a monomorphism}.

\noindent Observe that
$p_0\circ{}q_1=r\circ{}p_1\circ{}q_1:S_1\to{}S_0$ is homotopic to
a constant map as $(p_1\circ{}q_1)(S_1)=S_1$ and
$r(S_1)\varsubsetneq{}S_0$. Thus $H_1(p_0\circ{}q_1)=0$. Also,
$H_1(p_0\circ{}q_2)=H_1(r\circ{}p_1\circ{}q_2)=0$, by (i). Thus we
have

(2) $H_1(p_i\circ{}q_j)=0$ {\it if} $0\le{i}<j\le2$.

\noindent Hence, by (1), (1$'$) and (2), Lemma 4D.1 implies that
$\varphi$ is a monomorphism in this case as well.\qed
\enddemo

In the following two lemmas 6D.3 and 6D.4 we consider two graphs
$P_1=|K_1|$ and $P_2 =|K_2|$, and a special pseudo 2-manifold $M
\subset P_1\times{}P_2$ as in 6B. We keep the notation of section
6B.

\proclaim{Lemma~6D.3} Suppose $M$ is monotonically embedded in
$P_1 \times P_2$. If either $P_1$ or $P_2$ is the one-point union
of finitely many circles, then both $P_1$ and $P_2$ are circles
and $M=P_1\times P_2$.
\endproclaim

\demo{Proof} It is enough to consider the case where $P_1$ is the
one-point union of circles $S_1,\cdots,S_k$, where $k\ge1$. There
is a 1-cell $\sigma\in K_1$ lying in $S_1$. Then for any 1-cell
$\tau\in K_2$ laying in $P_2(\sigma)$ (this set is well defined
because $p_1$ is surjective) we have $P_1(\tau)=S_1$. (In fact, as
$p_2$ is monotone, by Properties (e$'$) and (f$'$) in 6B, we infer
that $P_1(\tau)$ is a circle. Hence $P_1(\tau)=S_1$, because
$\sigma\subset P_1(\tau)$.)

Therefore $S_1\times{}P_2(\sigma)\subset{}M$. Again, by an
argument as above, $P_2(\sigma)$ is a circle. Hence
$M=S_1\times{}P_2(\sigma)$. Consequently, $P_1=S_1$,
$P_2=P_2(\sigma)$ and $M=P_1\times P_2$, which completes the
proof. \qed
\enddemo

\proclaim{Lemma~6D.4} Suppose $M$ is monotonically embedded in
$P_1 \times P_2$. Then we have:

{\rm(i)} If $\rank{}H_1(M)\le3$ then both $P_1$, $P_2$ are circles
and $M=P_1\times P_2$.

{\rm(ii)} If $\rank{}H_1(M)=4$ then both $P_1$, $P_2$ are
$\theta$-curves.
\endproclaim

\demo{Proof} Let $v$ be a vertex of $K_1 (=K'_1)$ such that
$\rank{}H_1(P_2(v))\ge\rank{}H_1(P_2(v'))$ for each vertex
$v'\in{}K_1$, and let $w$ be a vertex of $K_2 (=K'_2)$ such that
$\rank{}H_1(P_1(w))\ge\rank{}H_1(P_1(w'))$ for each vertex
$w'\in{}K_2$. Put $S_1=P_1(w)$ and $S_2=P_2(v)$.

\ss First, we shall show that

\ss (1) {\it if either} $\rank{}H_1(S_1)=1$ {\it or}
$\rank{}H_1(S_2)=1$ {\it then
 both} $P_1$ {\it and} $P_2$ {\it are circles and} $M=P_1\times{}P_2$.

\ss \noindent Suppose $\rank{}H_1(S_1)=1$. Then
$\rank{}H_1(P_1(w'))=1$, so $P_1(w')$ is a circle, for each vertex
$w'\in{}K'_2$. By Property (h$'$) in 6B, $P_1=p_1(M)$ is the
circle $S^1$. Thus, by Lemma~6D.3, the consequent of implication
(1) holds. By a similar argument, the consequent holds as well if
$\rank{}H_1(S_2)=1$.

\ss Now, we will shall show that

\ss (2) {\it if} $\rank{}H_1(S_1)\ge2$ {\it and}
$\rank{}H_1(S_2)\ge2$ {\it then} $\rank H_1(M)\ge4$.

\ss \noindent We intend to apply Lemma~6D.2. Define
$q_i:S_i\to{}M$ by $q_1(x)=(x,w)$, $q_2(y)=(v,y)$. With this
notation, and for $S=\emptyset$, the hypotheses of Lemma 6D.2 are
fulfilled. Thus we infer that

\ss (3) $\rank{}H_1(M)\ge\rank{}H_1(S_1)+\rank{}H_1(S_2)$.

\ss This implies (2). Observe that (2) and (1) imply (i).

\ss Now assume $\rank H_1(M)=4$. Then, by (1),
$\rank{}H_1(S_i)\ge2$ for each $i$. Consequently, by (3), $\rank
H_1(S_1)=\rank H_1(S_2)=2$. It follows that each $S_i$ is a union
of two different circles which intersect in a connected set. We
shall show that $P_1=S_1$.

By Property (g$'$) in 6B, it suffices to show that
$P_1(\tau)\subset S_1$ for each 1-cell $\tau\in K_2$. Suppose, to
the contrary, $P_1(\tau)\nsubseteq S_1$ for a 1-cell $\tau$. By
our assumptions and Property (f$'$), $S_0=P_1(\tau)$ is a circle.
Define $q_0:S_0\to{}M$ by $q(x)=(x,y_0)$, where $y_0\in\tau$ is
any point. With these new data, the hypotheses of Lemma 6D.2 are
fulfilled again. Hence, by this lemma, $\rank H_1(M)\ge5$,
contrary to $\rank H_1(M)=4$. So, indeed $P_1=S_1$.

Thus $P_1$ is a union of two circles which intersect in a
connected set. As $P_1$ is connected, the intersection is an arc
or a point. The later case is excluded by Lemma 6D.3. So $P_1$ is
a $\theta$-curve.

By similar argument $P_2$ is a $\theta$-curve as well. This ends
the proof of (ii), which ends the proof of our lemma. \qed
\enddemo

\ms \centerline {\it Proof of Theorem {\rm 6D.1}} \ms

\noindent Suppose a special pseudo 2-manifold $M$ is a subset of
the product $Y_1\times{}Y_2$ of two curves. By Theorem~4B.1 we can
assume that $M$ is also monotonically embedded in the product
$P_1\times{}P_2$ of two graphs. Thus, the assumptions about $P_1$,
$P_2$ and $M$ from subsection 6B are fulfilled, and the hypothesis
of 6D.4 is fulfilled as well. From the algebraic assumptions about
$H_1(M)$ and Lemma 6D.4, it follows that both $P_i$ are either
circles or $\theta$-curves. Therefore, by Property (f$'$), any two
sets $P_2(\sigma), \sigma \in K_1$, and any two sets $P_1(\tau),
\tau \in K_2$, intersect. By the same Theorem 4B.1, there exist
two mappings $h_1:P_1\to Y_1$ and $h_2:P_2\to Y_2$ such that their
product mapping $h_1\times h_2:P_1\times P_2\to Y_1 \times Y_2$
coincides with the identity on $M$. Thus, by Property (l) from 6B,
we infer that both mappings $h_1:P_1 \to Y_1$ and $h_2:P_2 \to
Y_2$ are embeddings. Therefore we can think of $P_i$ as a subset
of $Y_i$, which completes the proof.\qed\bs

\proclaim{Theorem 6D.5} {\rm ({\bf R. Cauty})} Any closed surface
whose rank of $H_1$ is $\ge5$ can be embedded in a product of two
graphs.
\endproclaim

\ms\noi {\bf Note 1.} Combining Theorems 6D.1 and 6D.5 with other
known results we can summarize the final status of the problem of
embeddability of 2-manifolds into products of two curves in the
following way.

For closed surfaces we have: Any surface which can be
embedded in the product of two curves can be also embedded in the
product of two graphs (see Corollary 5B.2). Among orientable
surfaces the sphere $\Sb^2$ is the only one not embeddable in any
product of two graphs (see \cite{Bo3} and \cite{Ku}, cf. 6A, the
Remarks following 5D.6, 6D.1, and Theorem 6C.1). Any
non-orientable surface, except the following five: the projective
plane $\P^2$, the Klein bottle $\K=\P^2\# \P^2$, $\P^2\# \P^2\#
\P^2$, $\P^2 \# \P^2 \# \P^2 \# \P^2$ and $\P^2\# \P^2\# \P^2\#
\P^2\# \P^2$, can be embedded in a product of two graphs (see
\cite{C1}). (The above sequence of surfaces could have been
written down in a more concise form using the well known
equivalence $\P^2\# \P^2\# \P^2\approx\T^2\# \P^2$.) The rank of
the 1st homology of these non-orientable surfaces is 0, 1, 2, 3,
4, respectively. Therefore, there are exactly six closed surfaces
not embeddable in any product of two curves: one orientable -
$\Sb^2$ - with genus 0, and five non-orientable with genus
$1,2,3,4$ and $5$. Thus we have exactly six exceptional closed
surfaces.

For surfaces with non-empty boundary we show in section 6F that
each of them can be embedded in the product of the simple triod
$T$ and the unit interval $I$ (the product is often called the
"three-page book" ). (This result is probably known to a number of
topologists.)\ss

\noindent{\bf Note 2.} For each $n\geq 3$ there exist exceptional
closed orientable $n$-manifolds different from $\Sb^n$ (see
Corollary 5D.6).\qed

\bs We are going to prove Theorem 6D.5 by an argument somewhat
different from the original one devised by R. Cauty. Let $M$ be a
closed 2-manifold with $\rank{}H_1(M)\ge5$, we have to show that
$M$ can be embedded in the product of two graphs.

The original proof is based on the following idea: if $M_1$ and
$M_2$ are two closed surfaces and each can be embedded in a
product of two graphs, then the same is true of the connected sum
$M_1\#M_2$ (\cite{C1}, \cite{Ku}). Moreover, $\rank{}H_1(M\#\T^2)=
\rank{}H_1(M)+2$. Since any torus lies in a product of two graphs,
and any closed orientable surface different from $\Sb^2$ can be
obtained from a torus by repeated applications of the operation of
connected sum with $\T^2$, any closed orientable surface can be
embedded in a product of two graphs (\cite{Ku}). On the other
hand, any closed non-orientable surface whose rank of $H_1$ is
$\ge5$ can be obtained by the same construction starting from a
closed non-orientable surface whose rank of $H_1$ is either 5 or
6. So, the proof will be done once we show that either of the two
surfaces embeds in a product of two graphs. The desired embeddings
have been constructed by R. Cauty \cite{C1}.

Our approach is slightly different: Theorem 6D.5 is obtained as a
direct consequence of Theorem 6C.1 (about embeddability of
orientable surfaces) and the following two lemmas (about
embeddability of non-orientable surfaces).

\proclaim{Lemma 6D.6} Any closed non-orientable surface $M$ whose
rank of $H_1$ is odd and $\ge5$ can be embedded in a product of
two graphs.
\endproclaim

\demo{Proof} Suppose $\rank H_1(M)=2k+1$, where $k\ge2$. We shall
show that $M$ can be embedded in the product $P_1\times P_2$ of
two graphs. Define $P_1=S_1\cup S_1'$, where $S_1$ and $S_1'$ are
oriented circles such that their intersection $S_1\cap S_1'$ is a
union of $k$ disjoint arcs $A_1,\cdots,A_k$, and the orientations
induce the same orientation on $A_1$ and opposite orientations on
$A_2$. And define $P_2=S_2\cup S_2'$ to be a $\theta$-curve, with
$S_2$ and $S_2'$ being oriented circles such that the intersection
$S_2\cap S_2'$ is an arc $A$, and the orientations induce opposite
orientations on $A$. It remains to construct a non-orientable
surface $\widetilde{M}$ in $P_1\times P_2$ homeomorphic to $M$.

There are 1-dimensional $CW$ complexes $K_1$ and $K_2$ such that:
$P_i=|K_i|$, each $A_j$ is the carrier of a 1-cell $\sigma_j$ of
$K_1$ oriented coherently with $S_1$, and $A$ is the carrier of a
1-cell $\tau$ of $K_2$ oriented coherently with $S_2$. Suppose
$\partial\tau=w_1-w_0$, where $w_i$ are vertices of $\tau$.

Consider the tori $T=S_1\times{}S_2$, $T'=S_1'\times{}S_2'$. Note
that their intersection
$$T\cap{}T'=(S_1\cap{}S_1')\times(S_2\cap{}S_2')=
(A_1\cup\cdots\cup A_k)\times{}A$$ is the union of $k$ disjoint
2-cells $\sigma_1\times{}\tau$, $\cdots$ $\sigma_k\times{}\tau$ of
$K_1\kw{}K_2$. Put $\widetilde{M}=N\cup{}N'$, where
$N=T\setminus(\kn{\sigma_1}\times{}\kn{\tau}\cup\cdots\cup
\kn{\sigma_k}\times{}\kn{\tau})$ and
$N'=T'\setminus(\kn{\sigma_1}\times{}\kn{\tau}\cup\cdots\cup
\kn{\sigma_k}\times{}\kn{\tau})$. One easily sees that
$\widetilde{M}$ is a (connected) surface. We shall show that
$\widetilde{M}$ is not orientable.

For suppose $\widetilde{M}$ is orientable. Then
$H_2(\widetilde{M})\approx\Z$ and each 2-cell of $K_1\kw{}K_2$
lying in $\widetilde{M}$ can be assigned an orientation such that
the 2-chain $z$ which is the sum of oriented 2-cells with
coefficient 1 is a 2-cycle (representing a generator of
$H_2(\widetilde{M})$). Then $z=c+c'$, where $c$ (respectively
$c'$) is the sum of the oriented 2-cells lying in $N$
(respectively in $N'$). Since $N$ is orientable we may assume that
the orientations of the 2-cells lying in $N$ are induced by the
orientations of $S_1$ and $S_2$. As
$0=\partial{}z=\partial{}c+\partial{}c'$ we have
$\partial{}c=-\partial{}c'$. Since the oriented 1-cell
$\sigma_1\otimes{}w_1$ enters $\partial{}c$ with coefficient 1, it
must enter $\partial{}c'$ with coefficient -1. Let $\tau_1\in K_2$
be the 1-cell with carrier $S_2\setminus \kn{A}$ oriented
coherently with $S_2$. Let $\tau'_1\in K_2$ be the 1-cell with
carrier $S'_2\setminus \kn{A}$ oriented coherently with $S'_2$. It
follows that the 2-cell $\sigma_1\times{}\tau'_1$ lying in $N'$
has been assigned the orientation $\sigma_1\otimes{}\tau'_1$.
Consequently, the orientations of the remaining 2-cells lying in
$N'$ are also induced by the orientations of $S'_1$ and $S_2'$.
But this leads to a contradiction: the oriented 1-cell
$\sigma_2\otimes{}w_1$ enters $\partial{}c$ with coefficient 1,
and $\partial{}c'$ with coefficient 1 as well.

Now we prove that $\rank{}H_1(\widetilde{M})=2k+1$. Note that
$\widetilde{M}$ is the union of two tori with the interiors of $k$
disjoint discs removed, hence the Euler characteristic
$\chi(\widetilde{M})=\chi(N)+\chi(N')-\chi(N\cap{}N')=(-k)+(-k)-0=-2k$.
So, $\rank{}H_1(\widetilde{M})=1-\chi(\widetilde{M})=2k+1$. It
follows that $\widetilde{M}\thickapprox{M}$. \qed
\enddemo\bs

The ideas used in the above proof can be easily generalized to
give a proof of the following much more general result.\ms

\proclaim {\bf Proposition 6D.7} Let $M$ be an $n$-manifold,
$n\ge2$, and let $h:M\to M$ be an involution. Suppose $M=M_0\cup
M_1$ is a union of two bordered $n$-manifolds such that
$h(M_i)=M_{1-i}$, where $i=0,1$, and $M_0\cap M_1=N_0\sqcup N_1$
is the disjoint union of orientable $(n-1)$-manifolds such that
$h(N_i)=N_i$. If $h$ preserves an orientation of $N_0$ and
reverses an orientation of $N_1$, then $M$ is not
orientable.\qed\endproclaim\ms

In the next lemma we use some graphs $P_n$, $n\ge3$. To describe
$P_n$ we first pick some points in the plane $\R^2$. Denote
$o=(0,0)$, and let $v_i=(\cos{2\pi i\over n}, \sin{2\pi i\over
n})$ for $i=0,...,n-1$. Then define
$$P_n=\bigcup^{n-1}_{i=0}(ov_i\cup v_i v_{i+1}),$$ where
$v_n=v_0=(1,0)$. Thus $P_n=|K_n|$, where $|K_n|$ is the
(simplicial) complex composed of segments and their endpoints. Let
$\widetilde{P}_n$ denote the union $\widetilde{P}_n=P_n\cup S$,
where $S$ is a circle such that $S\cap P_n=v_0v_1$.

\proclaim{Lemma 6D.8} Any closed non-orientable surface whose rank
of $H_1$ is equal to $2n$ with $2n\ge6$, can be embedded in a
product of two graphs. Moreover, if $n$ is odd it can be embedded
in $P_n\times P_n$, if $n$ is even it can be embedded in
$\widetilde{P}_{n-1}\times \widetilde{P}_{n-1}$.
\endproclaim

\demo{Proof}We shall prove the stronger version. To this end, we
have to consider two cases.

Case 1: {\it n is odd}. To complete the proof in this case it is
enough to construct a non-orientable surface $M$ in $P_n\times
P_n$ whose rank of $H_1$ is equal to $2n$. To this end, for each
$i=0, ...,n-1$, consider: \ss the oriented circles
$S_i=ov_i\cup{v_iv_{i+1}}\cup{v_{i+1}o}$ in $P_n$, \ss the tori
$T_i=S_i\times{}S_i$ in $P_n\times P_n$, and \ss the 2-cells
$D_i=ov_i\times{ov_i}$ in $T_i$. \ss\noi (Here $v_n=v_0$.) Then
$T_i\cap{}T_{i+1}=D_{i+1}$, where $T_n=T_0$ and $D_n=D_0$. Put
$$M=M_0\cup\cdots\cup M_{n-1},$$ where $M_i=T_i \setminus (\kn
D_i\cup\kn D_{i+1}).$ First we shall show that $M$ is a closed
(connected) surface. Note that $M$ is connected because each $M_i$
is connected and contains $(o,o)$. One easily sees that $M$ is
locally planar at each point different from $(o,o)$. Thus, it
remains to show that $M$ is locally planar at $(o,o)$ too. To this
end, note that $K$ is a pseudo-manifold, and that all the 2-cells
of $K$ containing $(o,o)$ can be arranged into a sequence

\ms \hskip1cm $ ov_0\times ov_1 \rightarrow ov_2\times ov_1
\rightarrow ov_2\times ov_3\rightarrow \cdots \rightarrow
ov_{n-1}\times ov_0\rightarrow$

\ms \hskip1cm $ov_1\times ov_0 \rightarrow ov_1\times ov_2
\rightarrow ov_3\times ov_2\rightarrow \cdots \rightarrow
ov_0\times ov_{n-1}.$ \footnote{There is an action of
$\Z_{2n}\cong \Z_2\bigoplus\Z_n$ on $M$, with $(o,o)$ as a unique
fixed point, such that the sequence represents a cyclic orbit of
the 2-cell $ov_0\times ov_1$, and each cell is followed by its
image under a generating homeomorphism. (The generating
homeomorphism is the restriction to $M$ of the homeomorphism
$h:P_n \times P_n \to P_n \times P_n$ defined by
$h(x,y)=(\rho(y),\rho(x))$, where $\rho:P_n \to P_n$ is the
rotation around $o$ through the angle $\frac{2\pi}{n}$.)}

\bs \noi(Each cell $ov_i\times ov_{i'}$ is followed by
$ov_{i'+1}\times ov_{i+1}$, where indices are reduced modulo $n$.
Note also that both cells of the $i$th column lie in $M_i$ and
have only the point $(o,o)$ in common.) Each cell borders upon the
next one along in turn:
$$\{o\}\times ov_1, ov_2 \times \{o\}, \cdots, ov_{n-1}\times \{o\},
\{o\}\times ov_0,$$ $$ov_1\times\{o\}, \{o\}\times ov_2, \cdots,
\{o\}\times ov_{n-1}, ov_0\times \{o\}.$$ This shows that the
union of the 2-cells is a 2-disc whose interior in $M$ contains
point $(o,o)$. Thus, $M$ is locally planar at that point. Hence
$M$ is a closed surface.

Now we shall show that $\chi(M)=1-2n$. First, one can prove by
induction that the set $N=(T_0\cup\cdots\cup T_{n-2})
\setminus(\kn D_1 \cup\cdots \cup\kn D_{n-2})$ is the connected
sum $T_0 \#\cdots\# T_{n-2}$. Then we note that $M=N'\cup
M_{n-1}$, where $N'=N\setminus(\kn D_0 \cup\kn D_{n-1})$. Also
note that $N'\cap{}M_{n-1}=\partial D_0 \cup
\partial D_{n-1}$ is the one-point union of two circles.
Consequently,
$$\chi(M)=\chi(N')+\chi(M_{n-1})-\chi(N'\cap{}M_{n-1})=[(n-2)(-2)+(-2)]+(-2)-(-1)=-2n+1.$$

It follows that $M$ is not orientable (as the Euler characteristic
of an orientable surface is even). Finally,
$\rank{}H_1(M)=1-\chi(M)=2n$, which concludes the proof in Case 1.

Case 2: {\it n is even}. As in Case 1, it is enough to construct a
non-orientable surface $M'$ in $\widetilde{P}_{n-1}\times
\widetilde{P}_{n-1}$ with rank of $H_1$ equal to $2n$. Since $n-1$
is odd and $\ge3$, by Case 1 there is a closed non-orientable
surface $M$ in $P_{n-1}\times P_{n-1}$ with rank of $H_1$ equal to
$2(n-1)$. Note that $A= $cl$(S\setminus P_{n-1})$ is an arc. It
follows that $M'=(M\cup (S\times S))\setminus (\kn A\times \kn A)$
is homeomorphic to $M\# \T^2$. Thus, $M'$ is a closed
non-orientable surface in $\widetilde{P}_{n-1} \times
\widetilde{P}_{n-1}$ with $\chi(M')=-2(n-1)-1$.Consequently, rank
$ H_1(M')=1+2(n-1)+1=2n$, which concludes the proof in Case 2.
\qed
\enddemo

\bs \centerline{6E. {\it Retracting products onto surfaces}}

\bs In connection with the proof in section 6A, let us make some
observations concerning the existence of retractions from
$Y_1\times Y_2$ onto closed surfaces $M\subset Y_1\times Y_2$,
where $Y_1$ and $Y_2$ are curves. We prove that the only surface
which is such a retract is the torus.\bs

\proclaim{Theorem 6E.1} Let $M$ be a closed surface in the product
$Y_1\times Y_2$ of two curves.

\ss {\rm(A)} If $M$ is a torus then $M=P_1\times P_2$, where each
$P_i$ is a circle in $Y_i$. Consequently, $M$ is a retract of
$Y_1\times Y_2$.

\ss {\rm(B)} If $M$ is not a torus then $M$ is not a retract of
$Y_1\times Y_2$.
\endproclaim

\bs \noi {\bf Remark.} Part (B) implies the Borsuk theorem
\cite{Bo3} (because any topological 2-sphere in a 2-dimensional
space is a retract of that space, see e.g., \cite {Kur, \S53, VI,
Theorem 1, p.354}).\qed

\bs Part (A) directly follows from condition (i) of Theorem 6D.1,
and the fact that any circle in a curve $Y$ is a retract of $Y$
(see e.g., \cite {Kur, \S53, VI, Theorem 1, p.354}). To prove Part
(B) we need the following fact. With no doubt, it is known to many
(all?) topologists. As we are unable to indicate an appropriate
reference, a proof is supplied for completeness.

 \proclaim{Lemma 6E.2} Let $M$ be a closed surface and let $S$
 be a circle in $M$. If $S$ is contractible in $M$ then $S$ bounds
 a disc $D\subset M$, i.e. $S=\partial D$.
\endproclaim

\demo{Proof} Let $p:\widetilde{M}\to M$ be the universal covering
projection. Then $\widetilde{M}$ is homeomorphic to either $\Sb^2$
or $\R^2$. By the homotopy lifting property of $p$, there is a
circle $\widetilde{S}\subset \widetilde{M}$ such that $p$ maps
$\widetilde{S}$ bijectively onto $S$. We shall show that there is
a disc $\widetilde{D}\subset \widetilde{M}$ bounded by
$\widetilde{S}$ such that {\it for every covering transformation}
$\tau:\widetilde{M}\to \widetilde{M}$ {\it different from the
identity we have}

\ss(1) $\tau(\widetilde{D})\cap\widetilde{D}=\emptyset$.

\ss First consider any covering transformation
$\tau:\widetilde{M}\to \widetilde{M}$ different from the identity.
If $\widetilde{M}\approx \Sb^2$ then $\tau$ is unique and is an
involution. Note that

\ss(2) $\tau(\widetilde{S})\cap \widetilde{S}=\emptyset$.

\ss\noi In fact, $\tau(\widetilde{S})\nsubseteq\widetilde{S}$,
otherwise $p$ would not be injective on $\widetilde{S}$ as $\tau$
has no fixed point. Then suppose (2) is false. Therefore, by the
above observation, we get a simple triod $T=oa \cup ob \cup oc$
such that $oa \cup ob \subset \widetilde{S}$ and $oc\subset
\tau(\widetilde{S})$. But this is impossible because $p$ is a
local homeomorphism and $p(T)\subset S$. This proves (2).

Next choose the disc $\widetilde D$ in $\widetilde M$ bounded by
$\widetilde S$. If $\widetilde M \approx \R^2$ then such a disc is
uniquely determined (by the Jordan theorem); if $\widetilde M
\approx \Sb^2$ then take $\widetilde D$ disjoint with
$\tau(\widetilde S)$ (which can be done by (2)). Observe that

\ss (3) $\tau(\widetilde{S})\cap \widetilde{D}=\emptyset.$

\ss\noi According to the choice of the disc, it is enough to prove
(3) only for $\widetilde M \approx \R^2$. Suppose, to the
contrary, that (3) is false. Then $\tau (\widetilde{S})\subset
\widetilde{D}$ by (2). Consequently, $\tau (\widetilde{D})\subset
\widetilde{D}$ because $\tau(\widetilde D)$ is a disc bonded by
$\tau(\widetilde S)$ (as $\tau$ is a homeomorphism of the plane).
Thus, by the Brouwer Fixed Point Theorem, $\tau$ has a fixed
point, a contradiction. This proves (3). Then note that

\ss (4) $\widetilde{D}\nsubseteq \tau(\widetilde{D}).$

\ss\noi Else $\tau^{-1}(\widetilde{D})\subset \widetilde{D}$, a
contradiction (as $\tau^{-1}$ is another covering transformation
different from the identity). Now we are ready to prove (1).

Suppose $\tau(\widetilde{D})\cap \widetilde{D}\neq\emptyset$. By
(4) there is an arc $ab\subset \widetilde{D}$ meeting
$\tau(\widetilde{D})$ only at $b$. It follows from (2) that
$ab\cup \tau(\widetilde{D})$ is a 2-umbrella in $\widetilde{M}$, a
contradiction. This proves (1).

It follows from (1) that $p$ maps $\widetilde{D}$ homeomorphically
onto $p(\widetilde{D})$. Therefore, $D=p(\widetilde{D})$ is the
desired disc in $M$ bounded by $S$, which completes the proof.
\enddemo

 \ms \centerline {\it Proof of Part {\rm(B)} of Theorem {\rm
6E.1}}

\ms Let $M$ be a closed surface different from the torus and lying
in the product $Y_1\times Y_2$ of two curves.

First we discuss the special case where $Y_1=P_1$ and $Y_2=P_2$
are graphs and $M$ is monotonically embedded in $P_1\times{P_2}$.
Let $K_i$ denote a regular complex such that $P_i=|K_i|$. Then we
adopt the notation of Subsection 6B.

By Property (b) in 6B, there are 1-cells $\sigma_0\in{K_1}$ and
$\tau_0\in{K_2}$ such that $\sigma_0\times\tau_0\subset{M}$. By
Property (f$'$) in 6D the sets $P_1(\tau_0)=S_1$ and
$P_2(\sigma_0)=S_2$ are circles. Let
$$M_1=S_1\times\tau_0\cup\sigma_0\times{S_2},{}
M_2=M\setminus(S_1\times\kn\tau_0\cup\kn\sigma_0\times{S_2}),$$
and let
$D_0=(S_1\setminus\kn\sigma_0)\times(S_2\setminus\kn\tau_0)$. Note
that $D_0$ is a disc, $M_1$ and $M_2$ are surfaces with boundary,
$M=M_1\cup{M_2}$ and $M_1\cap{M_2}=
\partial{M_1}=\partial{M_2}=\partial{D_0}$.

Now suppose $M$ is a retract of $P_1\times{P_2}$. Then $S=\partial
D_0$ is contractible in $M$, because $S$ is contractible in $D_0$
and $D_0 \subset P_1\times P_2$. Thus, by Lemma 6E.2 there is a
disc $D\subset M$ bounded by $S$. It follows that either $D=M_1$
or $D=M_2$. However, $M_1=(S_1\times{S_2})\setminus\kn{D_0}$ is
not a disc, hence $D=M_2$. Therefore, $M=M_1\cup D$ is a torus.
This contradicts our assumption, and ends the proof in the special
case.

\ss In the general case, by Theorem 4.3, there exist two graphs
$P_1$, $P_2$, an embedding $g:M\to P_1\times{P_2}$, and a map
$h:P_1\times{P_2}\to{}Y_1\times{Y_2}$ such that

\ss (i) {\it $g(M)$ is monotonically embedded in
$P_1\times{P_2}$}, and

\ss (ii) {\it $h\circ{g}:M\to Y_1\times{Y_2}$ is the inclusion}.

\ss Since Part (B) holds in the special case, it follows from (i)
that $g(M)$ is not a retract of $P_1\times{P_2}$. Suppose $M$ is a
retract of $Y_1\times{Y_2}$ and let $r:Y_1\times{Y_2}\to{M}$ be a
retraction. Then, by (ii), $g\circ{r}\circ{h}(x)=x$ for each
$x\in{g(M)}$. Thus $g\circ{r}\circ{h}:P_1\times{P_2}\to{g(M)}$ is
a retraction, a contradiction. This contradiction completes the
proof of Part (B) in the general case.
 \qed

\bs \centerline {6F. {\it Embedding bordered surfaces in the
"three-page book"}}\ms

By a {\it bordered surface} we mean a connected $2$-manifold with
non-empty boundary. The {\it "three-page book"} stands for the
product $I\times T$, where $T=a_0v\cup a_1v\cup a_2v$ is a simple
triod (i.e. a union of three segments $a_iv$ mutually disjoint
except a common endpoint $v$). In this section we prove the
following theorem.\ms

\proclaim{\bf Theorem 6F.1} Any bordered surface can be embedded
in the "three-page book".\endproclaim

\demo{Proof} To this end, consider a bordered surface $M$ and let
$k$ (${\geq 1}$) be the number of components of $\partial M$. It
is known that\ss

{\rm(1)} $\chi(M)\leq 2-k$.

\ss\noi In order to show that $M$ can be embedded in the
"three-page book" we recall a convenient classic model of $M$ (cf.
\cite{Ma, p. 43}). One gets the model by repeated application of
the following procedure: given a bordered surface $N$ we take a
long rectangular strip and pase both ends of the strip to the
boundary of $N$ so that the ends do not overlap on the boundary of
$N$. This produces a bordered surface whose Euler characteristic
is equal to $\chi (N)-1$. Attaching such a strip to a closed disc
we get either an annulus or the M\"{o}bius band. In the first case
the strip is called "annular" and in the other "twisted".

We begin with a closed disc $D$. Let $S_k$ denote the bordered
surface obtained from $D$ by attaching $k-1$ annular strips as
shown in Figure 6F.1.\ms

(FIGURE 6F.1)

\ms\noi The bordered surface which results is a model of the
$2$-sphere $\Sb^2$ with the interiors of $k$ disjoint discs
removed. Therefore, \ss {\rm(2)} $\chi(S_k)= 2-k$.\ss\noi
Consequently, \ss {\rm(3)} $\chi(M)\leq \chi(S_k)$.\ss \noi To get
the final model of $M$ we shall modify $S_k$ so as to keep the
number of boundary components the same and reduce the Euler
characteristic to $\chi (M)$. It is necessary to distinguish two
essentially different cases.

{\bf Case I}: {\it M is orientable}. Then $\chi(M)=l-k$, where $l$
is even. By (2) and (3) we infer that $\chi(M)=\chi(S_k)-2m_0$,
where $m_0\geq 0$. The desired model of $M$ is obtained by
attaching $m_0$ pairs of "crossed" annular strips to the boundary
of $D$ as shown in Figure 6F.2.\ms

(FIGURE 6F.2)

\ms {\bf Case II}: {\it M is non-orientable}. Then
$\chi(M)=\chi(S_k)- m_1$, where $m_1\geq1$. To get the desired
model of $M$ we attach $m_1$ "twisted" strips to $S_k$ as shown in
Figure 6F.3.\ms

(FIGURE 6F.3)

\ms To complete the proof of our theorem it is enough to embed the
above models in the product $I\times T$, where $T$ is the simple
triod. If $M$ is orientable then one can easily detect a copy of
the model in $I\times T$ as shown in Figure 6F.4. (For the disc
corresponding to $D$ we take the rectangle $I \times a_2v$. The
desired copy is the union of the rectangle, $k-1$ annular strips
lying in $I \times a_0v$, and $m_0$ pairs of "crossed" annular
strips lying in $I \times (a_0v \cup a_1v)$.)\ms

(FIGURE 6F.4)

\ms Now assume $M$ is non-orientable. Then we define a copy of the
model of $M$ in $I\times T$ as shown in Figure 6F.5. (This time
for the disc corresponding to $D$ we take the rectangle $I \times
a_2v$ with $m_1$ "gates". The desired copy is the union of the
rectangle, $k-1$ annular strips lying in $I \times a_0v$, and
$m_1$ "twisted" strips lying in $I \times (a_0v \cup a_1v)$.)\ms

(FIGURE 6F.5)

\enddemo

\bs \centerline {6G. {\it Embedding surfaces in the second
symmetric product of a curve}}\ms

As we have noted in Chapter 4, any $2$-dimensional compactum can
be embedded in the symmetric product $SP^3(\mu)$, where $\mu$
stands for the Menger universal curve. Illanes and Nadler have
asked about embeddability of $\Sb ^2$ in the second symmetric
product of the Menger curve (thus enquiring about possible
extension of the Borsuk Theorem 1.2 to embeddings into symmetric
products of curves). We shall show in a moment that the answer is
negative, but first we note a positive result.

Notice that the Cartesian product $X\times Y$ naturally embeds in
the second symmetric product $SP^2 (X\vee Y)$ (where $X$ and $Y$
are pointed spaces). This combined with Theorems 6C.1 and 6D.5
implies the following corollary.

\proclaim{Theorem 6G.1} Any closed surface embeddable in the product
of two graphs can be also embedded in the second symmetric product
of a graph. In particular, any orientable closed surface different
from $\Sb^2$, and any non-orientable surface $M$ with $\rank
H_1(M)\ge 5$ can be embedded in the second symmetric product of a
graph.\qed
\endproclaim

To get additional information about possible embeddability of the
six exceptional surfaces in the second symmetric product of a
curve, we first formulate a particular case of Corollary 4G.3.

\proclaim{Corollary 6G.2} Neither the $2$-sphere nor the
projective plane, nor the Klein bottle can be embedded $($even up
to shape$)$ in the second symmetric product of a curve.\qed
\endproclaim

This corollary answers in the negative the Illanes and Nadler
question.\ms

The next lemma is needed in the proof of our final theorem.

\proclaim{Lemma 6G.3} Let $X$ be the union of three different
topological $n$-tori $T_1,T_2,T_3$, $n\geq2$, such that the
intersection of any two different tori is an $(n-1)$-torus, and
the intersection of all is non-void and not an $(n-1)$-torus. Then
X is not embeddable in any product of $n$ curves.
\endproclaim

\demo{Proof} Suppose $X$ can be embedded in the product $Y_1
\times\cdots\times Y_n$ of $n$ curves. Without loss of generality
we can assume that $X\subset Y_1 \times\cdots\times Y_n$. Let
$\star=(\star_1,\cdots,\star_n)$ be a point of $T_1 \cap T_2 \cap
T_3$. By Theorem 5D.5, for each $i=1,2,3$, there exist topological
circles $S_{i,1}\subset Y_1,\cdots, S_{i,n} \subset Y_n$ such that
$$T_i=S_{i,1} \times\cdots\times S_{i,n}.$$ Call $S_{i,j}$
the $j$th {\it factor} of $T_i$. The $j$th factors of two tori
either coincide or meet at $\{\star_j\}$ only. In the latter case
all the remaining factors coincide. Indeed, by our hypothesis, for
any two different elements $i_1,i_2$, the intersection
$$T_{i_1} \cap T_{i_2}=
(S_{i_1,1}\cap S_{i_2,1})\times\cdots\times (S_{i_1,n}\cap
S_{i_2,n})$$ is an $(n-1)$-torus. It follows that $S_{i_1,j}\cap
S_{i_2,j}$ is a one-point of the form $\{\star_j\}$ for exactly
one $j$, and $S_{i_1,k}=S_{i_2,k}$ for $k\neq j$. Then we say that
$T_{i_1}$ and $T_{i_2}$ {\it differ at index} $j$. Without loss of
generality we can assume that \ss

\itemitem{(1)} $T_1$ {\it and} $T_2$ {\it differ at index} 1. \ss

\noi Let us notice that

\itemitem{(2)} $T_1$ {\it and} $T_3$ {\it differ at index} $j>1$. \ss

\noi Indeed, otherwise $S_{3,1}\cap S_{1,1}=\{\star_1\}$ and
$S_{3,j}=S_{1,j}$ for $j>1$. If $S_{3,1}=S_{2,1}$ then $T_3=T_2$,
a contradiction. Thus $S_{3,1}\neq S_{2,1}$, hence $S_{3,1}\cap
S_{2,1}=\{\star_1\}$ and $S_{3,j}=S_{2,j}$ for $j>1$. Therefore,
the intersection $T_1\cap T_2\cap T_3$ is an $(n-1)$-torus, a
contradiction. Now, by (1) and (2), we infer that $$T_2 \cap
T_3=\{\star_1\}\times S_{1,2}\times\cdots\times S_{1,j-1}\times
\{\star_j\}\times S_{1,j+1}\times\cdots\times S_{1,n}$$ is an
$(n-2)$-torus, a contradiction. This ends the proof. \qed\enddemo

In our final observation we point out the fact that the symmetric
product produces spaces more complex than the corresponding
Cartesian product does. Already the second symmetric product may
transform a relatively simple graph into a 2-dimensional
polyhedron which is not embeddable in any Cartesian product of two
curves. In fact, we have the following

\proclaim{Theorem 6G.4} The $n$th symmetric product, $n\geq 2$, of
a bouquet of $n+1$ topological circles cannot be embedded in the
product of $n$ curves.
\endproclaim

\demo{Proof} Let $\bigvee_{i=1}^{n+1} S_i$ be a bouquet of $n+1$
topological circles with base point $\ast$. We have to show that
$SP^n(\bigvee_{i=1}^{n+1} S_i)$ can not be embedded in the product
of $n$ curves. To this end, it is enough to point out a subset $Q$
of $SP^n(\bigvee_{i=1}^{n+1} S_i)$ that is not embeddable in any
product of $n$ curves. The space $Q$ will be defined as a subspace
of an $n$-dimensional $CW$ complex. The cells $e_J$ of the complex
correspond to proper sets $J\subset\{1,\cdots,n+1\}$. If
$J=\{n_1,\cdots,n_k\}$, where $n_1<\cdots <n_k$, then $e_J$ is
defined as follows:
$$e_J=q(\{(x_1,\cdots,x_n)\in (\bigvee_{i=1}^{n+1} S_i)^n:
x_j\in S_{n_j}\; for\; j=1,\cdots,k;\;\; x_j=\ast \; for\; j>k
\}),$$ \noi where $q:(\bigvee_{i=1}^{n+1} S_i)^n\to
SP^n(\bigvee_{i=1}^{n+1} S_i)$ is the quotient map. We have
exactly one $0$-cell $e_{\emptyset}= \{q(\star,\cdots,\star)\}$.
One easily sees that $\dim e_J=|J|$ and $e_J\cap e_{J'}=e_{J\cap
J'}$. Now we are ready to define $Q$, put $Q=T_1\cup T_2 \cup
T_3$, where $T_1=e_{\{2,\cdots, n+1\}}$,
$T_2=e_{\{1,3,\cdots,n+1\}}$, $T_3= e_{\{1,2,4,\cdots,n+1\}}$. It
remains to prove that $Q$ is not embeddable in any product of $n$
curves. Notice the following properties of the cells: each $T_j$
is an $n$-torus, the intersection of any two tori $T_j$ is an
$(n-1)$-torus, each of the three tori contains the 0-cell, and the
intersection of all three tori is the $(n-2)$-torus
$e_{\{4,\cdots,n+1\}}$. Thus the hypotheses of Lemma 6G.3 are
fulfilled. It follows from this lemma that $Q$ has the desired
property, which ends the proof. \qed
\enddemo

\bs \centerline {6H. {\it Embedding surface-like continua in
products of two curves}} \ms

Here we apply Cauty's Theorem 6D.5 to derive a result on
embeddability of some special surface-like continua into products
of two curves.

Let $p$ be a prime number. We denote the mapping cylinder of a
$p$-to-1 covering map $f = f_p: \Sb^1 \to \Sb^1$ by $M(p)$. Let
$\partial M(p)$ denote the circle in $M(p)$ corresponding to the
domain of $f$. Now we construct an inverse sequence of
2-dimensional polyhedra $P_i$. Denote by $P_1$ the $2$-sphere with
a triangulation $K_1$, and assume a 2-dimensional polyhedron $P_i$
with a triangulation $K_i$ has been defined. Replace each
$2$-simplex $\sigma$ of $K_i$ by $M(p)_\sigma$, where
$M(p)_\sigma$ is obtained from $M(p)$ by identifying $\partial
M(p)$ with the boundary $\partial \sigma$. In this way we obtain a
polyhedron $P_{i+1}$. Define $\varphi_i:P_{i+1} \to P_i$ to be any
mapping satisfying the conditions: \roster
\item "" (1) $\varphi_i (x) = x$ for each $x\in |K_i^{(1)}|$, and
\item "" (2) ${\varphi_i}{|M(p)_\sigma}$ is a relative homeomorphism of
$(M(p)_\sigma,\partial M(p)_\sigma)$ onto $(\sigma,b(\sigma))$,
where $b(\sigma)$ is the barycenter of $\sigma$.
\endroster

\noi Let $K_{i+1}$ be a triangulation of $P_{i+1}$ such that $
$diam$ [\varphi_j \circ \cdots \circ \varphi_i(\sigma)] <
1/2^{i-j}$ for each 2-simplex $\sigma \in K_{i+1}$. Then we define
the {\it Pontrjagin surface mod} $p$ as the inverse limit:
$$
\Pi_p = \varprojlim \{ P_i,\varphi_i \}.
$$

Then using Theorem 6D.5 and a technique analogous to that in
\cite{Dr$_2$, Lemma 2}, one can prove the following:

\proclaim{Theorem 6H.1} The Pontrjagin surface mod $2$ can be
embedded in a product of two curves.\qed\endproclaim

\bs\centerline{\bf Problems to Chapter 6}\bs

A graph $P$ is said to be a {\it quasi-factor} of a closed surface
$M$ if $M$ can be surjectively embedded in $P\times Q$, where $Q$
is another graph. If $P$ contains no proper quasi-factor of $M$
then it is called {\it minimal}. If $M$ is monotonically embedded
in $P\times Q$ then $P$ is said to be a {\it monotone
quasi-factor} of $M$. As in the case of quasi-factors, we can
define {\it minimal monotone quasi-factors} of $M$. The results
from subsection 6C show that the graph $\Theta_{m+1}$ is a minimal
(monotone) quasi-factor of any orientable surface of genus $m$.

\proclaim {Problem 6C.1} {\it Determine all quasi-factors} $(${\it
surjective monotone quasi-factors}$)$ {\it of a closed surface}.
\endproclaim

\proclaim {Problem 6C.2} {\it Determine all minimal quasi-factors}
$(${\it minimal monotone quasi-factors}$)$ {\it of a closed
surface}.
\endproclaim

\proclaim {Problem 6D.1} {\it Suppose a closed surface} $M$ {\it
can be embedded in the second symmetric product of a curve}. {\it
Can} $M$ {\it be embedded in the second symmetric product of a
graph}? $($Same question for a polyhedron in place of $M$.$)$
\endproclaim

\proclaim{Problem 6G.1} Is it possible to embed the remaining
three exceptional surfaces $($ that is, the non-orientable
surfaces with genus $3, 4$ and $5$$)$ in the second symmetric
product of a graph?
\endproclaim

\proclaim{Problem 6H.1} Can the Pontrjagin surface mod $p$, $p
\neq 2$, be embedded in a product of two curves?
\endproclaim

\proclaim{Problem 6H.2} Is each Pontrjagin surface a quasi
$2$-manifold? \endproclaim

\proclaim{Problem 6H.3} Is it possible to embed any
$3$-dimensional product of two Pontrjagin surfaces in a product of
three curves?
\endproclaim

\bs \centerline{\bf APPENDIX}

\bs The main result of this Appendix is Theorem~A1 below which
expresses certain property of tensor product. All groups under
discussion are Abelian.

\proclaim{Theorem A1} Let If $f:A\to{}B$ be a homomorphism of
groups. If $G$ is a group which is either non-torsion or a direct
sum of cyclic groups, and
$$f\otimes 1_{G}:A\otimes{}G\to{}B\otimes{}G$$
is non-trivial then, for each $k\ge1$, the homomorphism
$$
f\otimes 1_{\bigotimes^kG}:A\otimes{\bigotimes{}^kG}
\to{}B\otimes{\bigotimes{}^kG},
$$
is non-trivial as well.
\endproclaim

To prove this result we need the following Lemmas A2, A3 and A4.
Lemma~A2 shows that Theorem~A1 can be derived from its special
case where $f$ is the inclusion homomorphism.

\proclaim{Lemma A2} Let $f:A\to{}B$ be a homomorphism and let
$i:f(A)\to{}B$ denote the inclusion homomorphism. Then $f\otimes
1_G:A\otimes{}G\to{}B\otimes{}G$ is non-trivial if and only if
$i\otimes 1_G:f(A)\otimes{}G\to{}B\otimes{}G$ is non-trivial, for
any Abelian group $G$.
\endproclaim

\demo{Proof} Let $f':A\to{}f(A)$ be defined by $f'(a)=f(a)$, for
all $a\in{}A$. Then $f=i\circ{}f'$, so $f\otimes 1_G=(i\otimes
1_G)\circ(f'\otimes 1_G)$. As $f'$ is an epimorphism, $f'\otimes
1_G$ is an epimorphism as well, the conclusion follows. \qed
\enddemo

\proclaim{Lemma A3} Let $i:A\to{}B$ be the inclusion of groups,
where $A$ is a $p$-group for some prime $p$. If $G_1$ and $G_2$
are groups such that both homomorphisms
$$i\otimes 1_{G_k}:A\otimes{}G_k\to{}B\otimes{}G_k \ ,$$
$k=1,2$, are non-trivial then
$$i\otimes 1_{G_1\otimes{}G_2}:A\otimes(G_1\otimes{}G_2)\to{}
B\otimes(G_1\otimes{}G_2)$$ is non-trivial as well.
\endproclaim

\demo{Proof} First we prove this lemma under additional assumption
of $B$ being a $p$-group. Let $G_k^{(p)}$ be a $p$-basic subgroup
of $G_k$ (cf. \cite{F, p. 136}), and let $i_k$ denote the
inclusion $G_k^{(p)}\hookrightarrow{}G_k$. By \cite{F, Theorem
61.1, p.261},
$1_C\otimes{}i_k:C\otimes{}G_k^{(p)}\to{}C\otimes{}G_k$ is an
isomorphism for any $p$-group $C$. Consequently, since
$i\otimes1_{G_k}$ is non-trivial and both $A$ and $B$ are
$p$-groups, the homomorphism
$$i\otimes1_{G_k^{(p)}}:
A\otimes{}G_k^{(p)}\to{}B\otimes{}G_k^{(p)}$$ is non-trivial.
Since $G_k^{(p)}$ is a direct sum of cyclic $p$ groups and
infinite cyclic groups, and since tensor product commutes with
direct sums, it follows that there exists a cyclic (infinite or a
$p$-group) subgroup $H_k$ of $G_k^{(p)}$ (a direct summand of
$G_k^{(p)}$), such that

\ms \itemitem{(i)} {\it the homomorphism}
$i\otimes1:A\otimes{}H_k{}\to{}B\otimes{}H_k$ {\it is
non-trivial}. \ms

Since $\Z\otimes{}C\cong{}C\otimes\Z\cong{}C$ for any group $C$,
and $\Z(p^r)\otimes\Z(p^s)\cong\Z(p^t)$, where $t=\min(r,s)$ (cf.
\cite{F, p. 255)}, it follows that $H_1\otimes{}H_2$ is isomorphic
to either $H_1$ or $H_2$. Consequently, by (i), the homomorphism
$$i\otimes1:A\otimes(H_1\otimes{}H_2)\to{}
B\otimes(H_1\otimes{}H_2)$$ is non-trivial. Since
$H_1\otimes{}H_2$ is a direct summand of
$G_1^{(p)}\otimes{}G_2^{(p)}$, it follows that

\ms \itemitem{(ii)} {\it the homomorphism}
$i\otimes1:A\otimes(G_1^{(p)}\otimes{}G_2^{(p)})\to{}
B\otimes(G_1^{(p)}\otimes{}G_2^{(p)})$ {\it is non-trivial}. \ms

Now, by \cite{F, p. 255, property (C)} and \cite{F, Theorem 61.1,
p. 261}, for any $p$-group $C$, we have the following sequence of
(natural) isomorphisms
$$
(C\otimes{}G_1^{(p)})\otimes{}G_2^{(p)}
@>{(1\otimes1)\otimes{}j_2}>> (C\otimes{}G_1^{(p)})\otimes{}G_2
@>{(1\otimes{}j_1)\otimes1}>> (C\otimes{}G_1)\otimes{}G_2 \ ;
$$
so
$$
1\otimes(j_1\otimes{}j_2): C\otimes(G_1^{(p)}\otimes{}G_2^{(p)})
\to C\otimes(G_1\otimes{}G_2)
$$
is a (natural) isomorphism. Thus, by (ii), since $A$ and $B$ are
$p$-groups, the homomorphism $i\otimes1_{G_1\otimes{}G_2}$ is
non-trivial, which ends the proof of the lemma in case $B$ is a
$p$-group.

Now, we prove the lemma in the general case. The $p$-group $A$ is
contained in the $p$-component $B_p$ of the torsion part of $B$.
Let $i_1:A\hookrightarrow{}B_p$ and $i_2:B_p\hookrightarrow{}B$
denote the inclusions homomrphisms. Then $i=i_2\circ{}i_1$.
Since $i\otimes1_{G_k}$ is non-trivial, the homomorphism
$$i_1\otimes1_{G_k}:A\otimes{}G_k\to{}B_p\otimes{}G_k$$
is also non-trivial, for $k=1,2$. Consequently, referring to the
special case, the homomorphism
$$i_1\otimes1_{G_1\otimes{}G_2}:
A\otimes(G_1\otimes{}G_2)\to{}B_p\otimes(G_1\otimes{}G_2)$$ is
non-trivial. \ms

\noindent Since $B_p$ is a pure subgroup of $B$ (see \cite{F, p.
114)}, by \cite{F, Theorem 60.4, p. 259} it follows that
$$i_2\otimes1:
B_p\otimes(G_1\otimes{}G_2)\to{}B\otimes(G_1\otimes{}G_2)$$ is a
monomorphism. Consequently, $i\otimes1_{G_1\otimes{}G_2}=
(i_2\otimes1_{G_1\otimes{}G_2})\circ
(i_1\otimes1_{G_1\otimes{}G_2})$ is non-trivial, which proves the
lemma. \qed \enddemo

\proclaim{Lemma A4} Let $i:A\to{}B$ be an inclusion homomorphism
of groups. If both $A$ and $G$ contain an element of infinite
order then the homomorphism
$$i\otimes1_G:A\otimes{G}\to{}B\otimes{G}$$
is non-trivial.
\endproclaim

\demo{Proof} By the assumption $A$ contains an infinite cyclic
subgroup $C$. Since $G$ contains an element of infinite order, the
quotient $G/T(G)$ is not trivial and torsion-free. Let us consider
the following commutative diagram

$$
\CD C\otimes{G} @>{(i|C)\otimes1_G}>>
  B\otimes{G}\\
   @V{1_C\otimes\pi}VV
     @V{1_B\otimes\pi}VV \\
 C\otimes{(G/T(G))} @>{(i|C)\otimes1_{G/T(G)}}>>
  B\otimes{(G/T(G))} \ ,
\endCD
$$
where $\pi$ denotes the projection $G\to{G/T(G)}$.

Since $i|C:C\hookrightarrow{B}$ is a monomorphism and $G/T(G)$ is
torsion-free, by \cite{F, Theorem 60.6, p. 260},
$(i|C)\otimes1_{G/T(G)}$ is a monomorphism as well. Next
$1_C\otimes\pi$ is not trivial, because $\pi$ is an epimorphism
and $G/T(G)$ is not trivial. Consequently, by the commutativity of
the diagram, $(i|C)\otimes1_{G}$ is non-trivial. Thus
$i\otimes1_{G}$ is non-trivial.
 \qed \enddemo

\ms \centerline{\it Proof of Theorem~A1}

\ms By Lemma~A2 we may assume that $f:A\to{}B$ is the inclusion
homomorphism. Then we divide the proof into three parts.

\ms (I) First, we assume that $A$ is a torsion group (and $G$ is
an arbitrary (Abelian) group). Then $A$ is a direct sum of its
$p$-components $A_p$. By the assumption the homomorphism
$f\otimes{}1_G:A\otimes{G}\to{}B\otimes{G}$ is non-trivial.
Consequently, there is a prime $p$ such that the homomorphism
$$(f|A_p)\otimes1_G:
A_p\otimes{G}\hookrightarrow{B}\otimes{G}$$ is non-trivial. By
Lemma~A3, the homomorphism
$$
(f|A_p)\otimes1_{\bigotimes^kG}:
A_p\otimes{\bigotimes{}^kG}\hookrightarrow
{B}\otimes{\bigotimes{}^kG}
$$
is non-trivial as well. Consequently, the homomorphism
$f\otimes1_{\bigotimes^kG}$ is non-trivial, which ends the proof
in this case.

\ms (II) Now, we assume that both $A$ and $G$ contain an element
of infinite order. Then also $\bigotimes^kG$ contains an element
of infinite order. Thus, in this case, the theorem follows by
Lemma~A4.

Observe that (I) and (II) imply the conclusion of the theorem in
the case $G$ contains an element of infinite order.

\ms (III) Finally, we assume that $G$ is the direct sum of cyclic
groups. Let $G=\bigoplus{G_s}$, where each $G_s$ is a cyclic
group. Then ${\bigotimes^kG}$ is also a direct sum of cyclic
groups, in which an isomorphic copy of each $G_s$ appears at least
once. Consequently, ${\bigotimes^kG}$ has an isomorphic copy of
$G$ as a direct summand. Thus, in this case, the required
conclusion follows since the tensor product commutes with the
direct sum. This ends the proof of the theorem. \qed

\bs Note that the conclusion of Theorem~A1 holds for arbitrary
group $G$ provided the image $f(A)$ is a torsion group. In
general, if $G$ is a torsion group, the conclusion of Theorem~A1
does not hold, even if $G$ is a $p$-group such that
$G\otimes{G}\ne0$.

\proclaim{Example} Let $G=\Z(p^{\infty})\oplus\Z(p)$, where $p$ is
an arbitrary prime. Then $G$ is a $p$-group such that
$G\otimes{G}=\Z(p)$. Let $f:\Z\to\Z$ be a homomorphism defined by
$f(1)=p$. One can observe that $f\otimes1_{G}$ is not trivial
since $f\otimes{1}_{\Z(p^{\infty})}:\Z\otimes\Z(p^{\infty})
\to\Z\otimes\Z(p^{\infty})$ is an epimorphism (and
$Z\otimes\Z(p^{\infty})$ is isomorphic to $Z(p^{\infty})$). On the
other hand $f\otimes1_{G\otimes{G}}$ is trivial since
$f\otimes1_{\Z(p)}$ is trivial.
\endproclaim

%%%%%%%%%% References %%%%%%%%%%%%%%%%%%%%%%%%%%

\enddocument